\documentclass[12pt]{article}
\usepackage{amsmath}
\usepackage{amssymb}
\usepackage{amsthm}
\usepackage{amscd}
\usepackage{amsxtra}
\usepackage{verbatim}
\usepackage{xcolor}
\usepackage{color}
\usepackage{enumerate}
\usepackage{mathrsfs}
\usepackage{stmaryrd}
\usepackage[all]{xy} 

\usepackage{array,arydshln}
\usepackage{bm}

\allowdisplaybreaks

\numberwithin{equation}{section}

\definecolor{dblue}{rgb}{0,0,0.45}
\definecolor{red}{rgb}{0.7,0,0}

 \RequirePackage{geometry}
\newtheorem{theorem}{Theorem}[section]
\newtheorem{lemma}[theorem]{Lemma}
\newtheorem*{lemma*}{Lemma}

\newtheorem{proposition}[theorem]{Proposition}

\theoremstyle{definition}

\newtheorem{remark}[theorem]{Remark}
\newtheorem{definition}[theorem]{Definition}
\newtheorem{example}[theorem]{Example}

\theoremstyle{remark}

\newcommand{\E}{{\mathbb E}}
\newcommand{\N}{{\mathbb N}}
\newcommand{\R}{{\mathbb R}}

\newcommand{\cC}{{\mathcal C}}

\newcommand{\cF}{{\mathcal F}}

\newcommand{\cM}{{\mathcal M}}

\newcommand{\cP}{{\mathcal P}}

\newcommand{\cQ}{{\mathcal Q}}
\newcommand{\cS}{{\mathcal S}}

\newcommand{\cV}{{\mathcal V}}
\newcommand{\cW}{{\mathcal W}}
\newcommand{\cX}{{\mathcal X}}

\newcommand{\la}{\langle}
\newcommand{\ra}{\rangle}
\newcommand{\nn}{\nonumber}
\newcommand{\ve}{\varepsilon}

\newcommand{\vertiii}[1]{{\left\vert\kern-0.25ex\left\vert\kern-0.25ex\left\vert #1 
    \right\vert\kern-0.25ex\right\vert\kern-0.25ex\right\vert}}

\date{}
\begin{document}

\title{
Averaging principle for rough slow-fast systems of level 3
}
\author{   Yuzuru \textsc{Inahama} 
}
\maketitle

\begin{abstract}
The averaging principle for slow-fast systems of 
various kind of stochastic (partial) differential equations 
has been extensively studied.
An analogous result was shown 
for slow-fast systems of rough differential equations 
driven by random rough paths a few years ago
and the study of ``rough slow-fast systems" seems to be 
gaining momentum now.
In all known results, however, the driving rough paths are 
of level 2. 
In this paper we formulate rough slow-fast systems 
driven by random rough paths of level 3 
and prove the strong averaging principle of Khas'minski\u{\i}-type.
\end{abstract}

\section{Introduction}

Let $(w_t)$ and $(b_t)$ be two independent standard
(finite-dimensional) Brownian motions (BMs).
A slow-fast system of (finite-dimensional) stochastic differential 
equations (SDEs) of It\^o-type are given  by
\begin{equation} \nn%
\left\{
\begin{array}{lll}
X^{\ve}_t &=& x_0 + 
\int_0^{t}
 f(X^\varepsilon_s, Y^\varepsilon_s) ds
 +
\int_0^{t} \sigma (X^\varepsilon_s, Y^\varepsilon_s) db_s,
 \\
Y^{\ve}_t &=& y_0 + 
\varepsilon^{-1}\int_0^{t}
 g(X^\varepsilon_s, Y^\varepsilon_s) ds
 +
\varepsilon^{-1/2}\int_0^{t}
 h(X^\varepsilon_s, Y^\varepsilon_s) dw_s,
\end{array}
\right.
\end{equation}
where $0<\ve \ll 1$ is a small parameter.
The processes
$X^\varepsilon$ and $Y^\varepsilon$ are called 
the slow component and the fast component, respectively.
Suitable conditions are imposed on $g$ and $h$ so that 
the following so-called frozen SDE satisfies certain ergodicity
for every $x$:
\[
Y^{x,y}_t =y+ \int_0^t g(x, Y^{x,y}_t)dt 
     +  \int_0^t  h(x, Y^{x,y}_t) dw_t.
\]
An associated unique invariant probability measure 
is denoted by $\mu^x$. We set 
$\bar{f} (x) =\int   f(x,y) \mu^x (dy)$ 
and set $\bar{\sigma} (x)$ in a similar way.
Consider the following averaged SDE:
\[
\bar{X}_t = x_0 + \int_0^{t}
 \bar{f} (\bar{X}_s) ds+
\int_0^{t} \bar\sigma (\bar{X}_s) db_s.
\]
The averaging principle of this type, 
which was initiated by Khas'minski\u{\i} \cite{khas},
 claims that 
$X^\varepsilon$ converges to $\bar{X}$ in an appropriate sense 
as $\ve \searrow 0$.
Early results are well-summarized in \cite[Chapter 7]{fw}.
Since then 
the study of averaging principle has been developed greatly.
Today, so many kinds of slow-fast systems of 
stochastic (partial) differential equations are known to 
satisfy the averaging principle.

Rough differential equations (RDEs) are generalized 
controlled ordinary differential equations, 
which can be regarded  
as ``de-randomization" of usual SDEs.
The generalized controls (i.e. drivers) are called rough paths (RPs).
When the RPs are random, we have
systems quite similar to SDEs
even when there is no (semi)martingale property.
For basic information on rough path theory, see \cite{fh, fvbook} 
among others.

The study of slow-fast systems of random RDEs looks
quite natural and interesting.
The first averaging result was obtained by \cite{pix2},
which adopted fractional calculus approach to RP theory.
It was then followed by \cite{ina_thk}, which used 
the controlled path theory to prove the averaging principle.
Within just two years after that, four papers 
already appeared along this research direction.
A large deviation 
principle associated with the averaging principle 
was shown in \cite{yangxu}.
A convergence rate of the averaging principle was obtained in \cite{yjym}.
The averaging principle for slow-fast systems 
of semilinear rough partial differential equations was proved in \cite{llpx}.
In these three works, the fast component is driven by 
Brownian RP, while the slow component is driven by 
fractional Brownian RP or more general random RP.
In a recent paper \cite{phsx} the authors studied a slow-fast system 
whose fast component is driven by fractional Brownian noise
and proved an almost-sure version of averaging principle.

In the theory of Gaussian RPs, a prominent example 
is fractional Brownian RP (i.e. a canonical lift of fractional Brownian 
motion) with Hurst parameter $H \in (1/4, 1/2]$.
When $H \in (1/4, 1/3]$, the fractional Brownian RP is a 
(random) RP of level 3.
In all the existing works on slow-fast systems of RDEs, however,
the driving RPs are of level 2.
It is therefore quite natural and important to generalize 
the theory of rough slow-fast systems to the case of level 3.
Our main purpose in this paper is 
to construct rough slow-fast systems of level 3
and prove the strong averaging principle,
which can be regarded as a level-3 version of the preceding work \cite{ina_thk}.

The structure of this paper is as follows.
In Section \ref{sec.ass_res} we provide assumptions 
on the coefficients and driving RP
and then state our main theorem (Theorem \ref{thm.main}).
In Section \ref{sec.CPtheory} we carefully explain
controlled path theory of level 3.
Although this is basically known among experts, 
there seem to be few papers which actually elaborate on the subject.
Section \ref{sec.driving} is devoted to constructing random
mixed RPs which drive our slow-fast systems.
Unlike the other sections, this section is specific to the
level 3 case and has no counterpart 
in the preceding work \cite{ina_thk}.
Hence, for those who already understand the level 2 case in \cite{ina_thk}, this section is the most important.
Section \ref{sec.SF} consists of two parts.
In the first half we precisely formulate our rough slow-fast system
in a deterministic way.
The second half discusses probabilistic aspects of the system
including the proof of our main theorem.
However, since the second half is quite similar to a counterpart
in the level 2 case in  \cite{ina_thk}, our exposition is 
a little bit sketchy.


\bigskip\noindent
{\bf Notation:}
\\
\noindent
Before closing Introduction, 
 we introduce the notation which will be used throughout the paper.
We write $\N =\{1,2, \ldots\}$ and $\N_0 :=\N \cup \{0\}$.
We set $\llbracket 1, k\rrbracket :=\{ 1, \ldots, k\}$ for $k\in\N$.
Let $T \in (0, \infty)$ be arbitrary and we work on the time interval
$[0,T]$ unless otherwise specified.
For a subinterval $[a,b]\subset [0,T]$,
we set $\triangle_{[a,b]} =\{ (s,t)\in \R^2 \mid  a\le s\le t \le b\}$.
When $[a,b]= [0,T]$, we simply write $\triangle_{T}$ for this set.

Below, $\cV$, $\cV_i~(i\in\N)$ and $\cW$ are Euclidean spaces.
The set of bounded linear maps from $\cV$ to $\cW$ is denoted 
by $L(\cV, \cW)$, which coincides with the set of all real 
$(\dim \cW) \times (\dim \cV)$ matrices.
For $k \ge 2$, the set of $k$-bounded linear maps from 
$\cV_1\times \ldots\times \cV_k$ to $\cW$ is denoted by 
$L^{(k)}(\cV_1,\ldots, \cV_k; \cW)$.
There are natural identification as follows;
 $L(\cV, \cW)\cong \cV^* \otimes \cW$ and 
 $L^{(k)}(\cV_1,\ldots, \cV_k; \cW)\cong 
L( \cV_1 \otimes\ldots\otimes \cV_k, \cW)$.
As usual, the truncated tensor algebra 
of degree $k~(k\in\N)$ over $\cV$ is defined by 
$T^k (\cV) := \oplus_{i=0}^k \cV^{\otimes i}$, 
where we set $\cV^{\otimes 0}:=\R$.

\begin{itemize} 
\item
The set of all continuous path $\varphi\colon [a,b] \to\cV$
is denoted by $\cC ([a,b], \cV)$. 
With the usual sup-norm $\|\varphi\|_{\infty, [a,b]}$ on the $[a,b]$-interval,
$\cC ([a,b], \cV)$ is a Banach space.
The difference of $\varphi$ is frequently denoted by $\varphi^1$,
that is, $\varphi^1_{s,t} := \varphi_t - \varphi_s$ for $(s,t)\in \triangle_{[a,b]}$.

\item
Let $0< \gamma \le 1$.
For a path $\varphi \colon [a,b] \to \cV$,
the $\gamma$-H\"older seminorm is defined by 
\[
\|\varphi\|_{\gamma,[a, b]} :=\sup _{a \le s<t \le b} 
\frac{\left|\varphi_{t}-\varphi_{s}\right|_{\cV}}{(t-s)^{\gamma}}.
\]
If the right hand side is finite, we say $\varphi$ is 
$\gamma$-H\"older continuous on $[a,b]$.
The space of all $\gamma$-H\"older continuous paths on $[a,b]$ is denoted by 
$\cC^\gamma ([a,b], \cV)$.
The norm on this Banach space is
 $|\varphi_a|_{\cV}+\|\varphi\|_{\gamma,[a, b]}$.

\item
Let $\gamma >0$.
For a continuous map $\eta \colon \triangle_{[a,b]} \to \cV$, we set 
\[
\|\eta\|_{\gamma,[a, b]} :=\sup _{a\le s<t\le b} 
\frac{\left|\eta_{s,t}\right|_{\cV}}{(t-s)^{\gamma}}.
\]
If this is finite, then $\eta$ vanishes on the diagonal.
The set of all such $\eta$ with $\|\eta\|_{\gamma,[a, b]}<\infty$
is denoted by $\cC^\gamma_{(2)} ([a,b], \cV)$,
which is a Banach space with $\|\eta\|_{\gamma,[a, b]}$.

\item
When $[a,b] =[0,T]$, we write
$\cC (\cV)$, $\cC^\gamma (\cV)$, $\cC_{(2)}^\gamma (\cV)$
for these spaces 
and $\|\cdot\|_{\infty}$, $\|\cdot\|_{\gamma}$, $\|\cdot\|_{\gamma}$
for the corresponding (semi)norms for simplicity of notation.
For $z\in \cV$, we set $\cC_z (\cV):=\{ \varphi \in \cC (\cV)\mid \varphi_{0}=z\}$. 
We also set $\cC^\gamma_z (\cV)$ in a similar way.

\item
Let $U$ be an open set of $\cV$.
For $k \in \N_0$,  $C^k (U, \cW)$ stands for the set of 
$C^k$-functions from $U$ to $\cW$.
(When $k=0$, we simply write $C (U, \cW)$ 
instead of $C^0 (U, \cW)$.)
The set of bounded $C^k$-functions $f \colon U\to \cW$
whose derivatives up to order $k$ are all bounded 
is denoted by $C_{{\rm b}}^k (U, \cW)$, which is a Banach space with the norm
$\| f\|_{C_{{\rm b}}^k } := \sum_{i=0}^k \|\nabla^i f\|_{\infty}$.
(Here, $ \|\cdot\|_{\infty}$ stands for the usual sup-norm on $U$.)

\item
Let $\gamma \in (1/2,1]$ and $m\in\N$. 
If $w$ belongs to $\cC_0^\gamma (\cV)$, then we can define 
\[
S(w)^m_{s,t} := \int_{s\le t_1 \le \cdots\le t_m\le t}  dw_{t_1} \otimes\cdots\otimes dw_{t_m},
\qquad 
(s, t) \in \triangle_T
\]
as an iterated Young integral. 
We call $S(w)^m$ the $m$th signature of $w$.
It is well-known that $S_k (w)_{s,t} :=
(1, S(w)^1_{s,t}, \ldots, S(w)^k_{s,t})\in T^k (\cV)$ and 
Chen's relation holds, that is,
\[
S_k (w)_{s,t} = S_k (w)_{s,u}\otimes S_k (w)_{u,t}, \qquad  s\le u \le t.
\]
Here, $\otimes$ stands for the multiplication in $T^k (\cV)$.

\item
Let $\alpha \in (1/4, 1/2]$ and write $k:=\lfloor 1/\alpha\rfloor$. 
We recall the definition of 
$\alpha$-H\"older rough path ($\alpha$-RP or RP).
A continuous map 
$X=(1, X^{1},\ldots,  X^{k})\colon \triangle_{T} \to T^k (\cV)$
is called $\cV$-valued $\alpha$-RP if 
$\|X^i\|_{i \alpha} <\infty$ for all $i\in \llbracket 1,k \rrbracket$ and 
\begin{equation} \label{eq.0711-1}
X_{s, t}=X_{s, u}\otimes X_{u, t},
\qquad
s\le u \le t.
\end{equation}
holds in $T^k (\cV)$. (This is called Chen's relation.)
The set of all $\cV$-valued $\alpha$-RPs is denoted by 
$\Omega_{\alpha} (\cV)$.
With the distance 
$d_{\alpha} (X, \hat{X}) :=\sum_{i=1}^k\|X^i - \hat{X}^i\|_{i \alpha}$,
$\Omega_{\alpha} (\cV)$ is a complete metric space.
The homogeneous norm of $X$ is denoted by 
$\vertiii{X}_\alpha:= \sum_{i=1}^k \|X^i\|_{i\alpha}^{1/i}$.
The dilation by $\delta \in \R$ is defined by 
$\delta X =(1, \delta X^{1},\ldots,  \delta^k X^{k})$. 
It is clear that 
$\vertiii{\delta X}_\alpha=|\delta|\cdot\vertiii{X}_\alpha$.
A typical example of RP is $S_k (w)$ for $w \in \cC_0^\gamma (\cV)$ 
with $\gamma \in (1/2,1]$.
This is called a natural lift of $w$.
We view $S_k$ as a continuous map from $\cC_0^\gamma (\cV)$
to $\Omega_{\alpha} (\cV)$ and call it the lift map.

\item
Let $\alpha \in (1/4, 1/2]$ and write $k:=\lfloor 1/\alpha\rfloor$. 
We define $G\Omega_{\alpha} (\cV)$ to be the $d_\alpha$-closure
of $S_k (\cC_0^1 (\cV))$.
It is called the $\alpha$-H\"older 
geometric RP space over $\cV$
and is a complete and separable metric space.
It also coincides with the $d_\alpha$-closure
of $S_k (\cC_0^\gamma (\cV))$ for any $1/2 <\gamma \le 1$.
A geometric RP $X \in G\Omega_{\alpha} (\cV)$ satisfies another important algebraic property called the Shuffle relations.
To explain it, we identify $\cV =\R^d~(d=\dim \cV)$
and the coordinate of $X^i$'s are denoted by $X^{1,p}$, 
$X^{2,pq}$ and $X^{3,pqr}$ ($p,q,r \in \llbracket 1,d\rrbracket$).
Then we have
\begin{equation} \label{eq.0711-2}
X^{1,p}_{s,t}\,X^{1,q}_{s,t}= X^{2,pq}_{s,t}+X^{2,qp}_{s,t}, 
\quad
X^{1,p}_{s,t}\,X^{2,qr}_{s,t}=X^{3,pqr}_{s,t}+X^{3,qpr}_{s,t}+X^{3,qrp}_{s,t}
\end{equation}
for all $p,q,r \in \llbracket 1,d\rrbracket$ and $(s,t) \in \triangle_{T}$.
(When $1/3 <\alpha \le 1/2$, we only have the first formula.)
For basic information on $\alpha$-H\"older geometric RPs,
the reader is referred to \cite[Chapter 9]{fvbook}. 
\end{itemize}


\section{Assumptions and main result}\label{sec.ass_res}

In this section we first introduce natural assumptions 
on the coefficients and driving random RP
of the following slow-fast system and then state our main theorem.

Our slow-fast system of RDEs is given by 
\begin{equation} \label{def.SFeq}
\left\{
\begin{array}{lll}
X^{\ve}_t &=& x_0 + 
\int_0^{t}
 f(X^\varepsilon_s, Y^\varepsilon_s) ds
 +
\int_0^{t} \sigma (X^\varepsilon_s) dB_s,
 \\
Y^{\ve}_t &=& y_0 + 
\varepsilon^{-1}\int_0^{t}
 g(X^\varepsilon_s, Y^\varepsilon_s) ds
 +
\varepsilon^{-1/2}\int_0^{t}
 h(X^\varepsilon_s, Y^\varepsilon_s) dW_s.
\end{array}
\right.
\end{equation}
The time interval is $[0, T]$.
Here, $0 <\ve \le 1$ is a small parameter
and (the first level path of)
$(X^{\ve}, Y^{\ve})$ takes values in $\R^m \times \R^n$. 
The starting point $(x_0, y_0)$ is always deterministic and arbitrary.
(We will not keep track of the dependence on $T, x_0, y_0$.)
At the first stage, \eqref{def.SFeq} is a deterministic system of RDEs 
driven by an ($d+e$)-dimensional third-level
RP which is denoted by $(B, W)$. 
A precise definition of the system \eqref{def.SFeq}
will be given in Subsection \ref{subsec.deter}.

When we consider this slow-fast system of RDEs,
the following are imposed as our standing assumptions:
\begin{itemize}
\item
$\sigma\in C^4 (\R^m, L(\R^d, \R^m))$ and
$h \in C^4 (\R^m \times \R^n, L(\R^e, \R^n))$,
\item
$f \in C (\R^m \times \R^n, \R^m)$ and 
$g \in C(\R^m \times \R^n, \R^n)$
are locally Lipschitz continuous.
\end{itemize}
These guarantee the existence of a unique 
local solution of \eqref{def.SFeq}.
(This fact is well-known. See also Remark \ref{rem.loc.sol} below.)
Since we show the strong version of the averaging principle in this work, 
we assume that $\sigma$  depends only on the slow component.

We set 
\begin{equation}\label{def.tilde_g}
\tilde{g} (x,y) :=g (x,y) +  \frac12
\left\{ \sum_{i=1}^n \sum_{j=1}^e
\frac{\partial h^{kj} }{\partial y_i} (x,y) h^{ij}(x,y)
\right\}_{1\le k \le n},
\quad
(x,y) \in \R^m\times \R^n.
\end{equation}
Here, we wrote $h =\{h^{ij}\}_{1\le i \le n, 1\le j \le e}$.
The second term on the right hand side above is the standard It\^o-Stratonovich correction term when $x$ is viewed as a parameter.


To formulate our main theorem, we introduce more 
assumptions on these coefficients.

\medskip
\noindent
${\bf (H1)}$~$\sigma$ is of $C^4_{{\rm b}}$.

\medskip
\noindent
${\bf (H2)}$~$f$ is bounded and globally Lipschitz continuous.

\medskip
\noindent
${\bf (H3)}$~$h$ is globally Lipschitz continuous.

\medskip
\noindent
${\bf (H4)}$~$\tilde{g}$ is globally Lipschitz continuous.

\medskip
\noindent
${\bf (H5)}$~There exist constants $\gamma_1 >0$ and
$C >0$ such that, for all $x \in \R^m$ and $y \in \R^n$, 
\[
2\langle y, \tilde{g} (x,y)\rangle +  |h (x,y)|^2
\le -\gamma_1 |y|^2 + C(|x|^2 +1).
\]

\medskip
\noindent
${\bf (H6)}$~There exists a constant $\gamma_2 >0$ such that,
for all $x \in \R^m$ and $y_1, y_2 \in \R^n$, 
\[
2\langle y_1- y_2, \tilde{g} (x, y_1)-\tilde{g} (x,y_2)\rangle
+  | h(x,y_1)- h (x,y_2) |^{2}
 \le  -\gamma_2 |y_1-y_2|^{2}.
 \]

\medskip

Let  $\tfrac14 <\alpha_0 \le \tfrac13$ and let
$(\Omega, \mathcal{F}, {\mathbb P};
 \{\cF_t\}_{0\le t\le T})$ be a filtered probability space
satisfying  the usual condition. 
On this probability space, 
the following two independent random variables 
$w$ and $B=(B^1, B^2, B^3)$  are defined.
The former, 
$w =(w_t)_{0\le t \le T}$, is a standard $e$-dimensional 
$\{\cF_t\}$-BM. 
The Stratonovich RP lift of $w$ is denoted by $W=(W^1, W^2)$.
The latter, $B=\{(B^1_{s,t}, B^2_{s,t}, B^3_{s,t})\}_{(s,t)\in \triangle_T}$, is an 
$G\Omega_{\alpha} (\R^d)$-valued 
random variable (i.e., random RP) for every $\alpha \in (1/4,\alpha_0)$.
Here, $G\Omega_{\alpha} (\R^d)$ is the space of 
$\alpha$-H\"older geometric RPs over $\R^d$.
We assume that $(B^1_{s,t}, B^2_{s,t}, B^3_{s,t})$ is $\cF_t$-measurable 
for every $(s,t)\in \triangle_T$.

We assume the following condition on  the integrability of $B$.
Below, $\vertiii{B}_{\alpha}:= 
\|B^1\|_{\alpha}+ \|B^2\|_{2\alpha}^{1/2}+ \|B^3\|_{3\alpha}^{1/3}$ denotes the $\alpha$-H\"older homogeneous RP norm over the time interval $[0,T]$.

\medskip 
\noindent 
${\bf (A)}$~  For every $\alpha \in (1/4,\alpha_0)$ and 
$p \in [1,\infty)$,  we have $\E [\vertiii{B}_{\alpha}^p  ] <\infty$.

\medskip 

\noindent
Under this assumption, the mixed random RP $(B,W)$
and the slow-fast system \eqref{def.SFeq} of RDEs
driven by it 
can be defined in a natural way
(see Subsection \ref{subsec.probab} for precise definitions).
We will show the averaging principle for \eqref{def.SFeq}
when it is driven by this random RP.
Note that, unlike in the second-level case
in the preceding work \cite{ina_thk}, the ``Brownian component" of 
the mixed random RP is of Stratonovich-type.

Next, we introduce 
the frozen SDE and the averaged RDE 
associated with the slow-fast system 
\eqref{def.SFeq} in the usual way.
The frozen SDE is given as follows:
\[
Y^{x,y}_t =y+ \int_0^t \tilde{g}(x, Y^{x,y}_t)dt 
     +  \int_0^t  h(x, Y^{x,y}_t) d^{{\rm I}}w_t,
\]
Here, $(x,y) \in \R^m \times \R^n$ are deterministic and arbitrary
and $d^{{\rm I}}w_t$ stands for the standard It\^o integral
with respect to a standard $e$-dimensional BM $(w_t)$.
We are only interested in the law of $Y^{x,y}$
and hence any realization of BM will do.
Under the assumptions of Theorem \ref{thm.main} below, 
the Markov semigroup $(P^x_t)_{t\ge 0}$
defined by $P^x_t \varphi (y)= {\mathbb E} [\varphi (Y^{x,y}_t) ]$
for a bounded measurable function $\varphi$
has a unique invariant probability measure,
which is denoted by $\mu^{x}$.
(This fact is well-known. See \cite[Appendix A]{pix1} 
among many others.)

Define the averaged drift by 
$\bar{f} (x) =\int_{\R^n}   f(x,y) \mu^x (dy)$ for $x\in \R^m$.
The averaged RDE is given as follows:
\begin{equation} \label{def.avRDE}
\bar{X}_t = x_0 + 
\int_0^{t}
 \bar{f} (\bar{X}_s) ds
 +
\int_0^{t} \sigma (\bar{X}_s) dB_s
\end{equation}
Here, $x_0 \in \R^m$ is the same as in \eqref{def.SFeq}.
Under the assumptions of Theorem \ref{thm.main} below, 
$\bar{f}$ is again bounded and globally Lipschitz. 
(This fact is also well-known. See  \cite[Lemma A.1]{pix1} 
for example.)
Therefore, this RDE has a unique global solution
for every realization of $B=(B^1, B^2, B^3)$.
(See Propositions \ref{prop.0429} below for details.)

Now we are in a position to state our main result,
whose proof will be provided at the end of Section 5.
It claims that (the first level path of) the slow component 
of the slow-fast system \eqref{def.SFeq} of RDEs
converges to (the first level path of) the 
averaged RDE \eqref{def.avRDE} in $L^p$-sense as $\ve\searrow 0$.
Here, $\|\cdot\|_\beta$ stands for the $\beta$-H\"older (semi)norm of a usual path over the time interval $[0,T]$.
\begin{theorem} \label{thm.main}
Assume ${\bf (A)}$ and ${\bf (H1)}$--${\bf (H 6)}$.
Then, for every $p\in [1,\infty)$ and $\beta \in (\tfrac14, \alpha_0)$, 
we have
\[
\lim_{\ve \searrow 0} \E [\| X^\ve - \bar{X} \|_\beta^p] =0.
\]
\end{theorem}

\begin{remark} \label{rem.probsp}
The law of $X^\ve - \bar{X}$ is uniquely determined by the law of 
$B=(B^1, B^2, B^3)$
and the $e$-dimensional Wiener measure.
In fact, $X^\ve - \bar{X}$ is obtained as a functional of $B$ and $w$.
So, the choice of a filtered probability space that carries 
$B$ and $w$ does not matter. 
(Verifying the existence of such a filtered probability space is easy.)
 \end{remark}

\begin{example} \label{exmpl.B}
A prominent example of $B =(B^1, B^2,B^3)$ 
satisfying Assumption ${\bf (A)}$ in Theorem \ref{thm.main} is 
fractional Brownian RP with Hurst parameter $H \in (\tfrac14, \tfrac13]$.
In this case, $\alpha_0 =H$ in ${\bf (A)}$.
This is a canonical RP lift of $d$-dimensional fractional BM 
with Hurst parameter $H$. For more information, see \cite[Chapter 15]{fvbook}.
Concerning this example, we will provide more explanations in
 Remark \ref{rem.250115}.
\end{example}

\begin{remark} \label{rem.coefficient}
The assumptions on the coefficients $\sigma, h, f, g$ in Theorem \ref{thm.main} are stronger than those in 
the author's preceding work \cite{ina_thk} on the level-two case.
In particular, $\tilde{g}$ is not allowed to be of super-linear growth in the present paper.
This is not because of theoretical limitation, but simply because 
computations become quite involved in the level-three case.
(The assumptions in Theorem \ref{thm.main}  are basically 
similar to those in \cite{pix1, pix2}.)
It could be interesting to relax these assumptions.
 \end{remark}

\section{Controlled path theory of level 3}\label{sec.CPtheory}

In this section we collect basic results 
from the theory of controlled paths of level 3.
It should be noted that they are basically known.
For instance, \cite{boge} studied the theory of controlled paths of 
any level.
For the level 3 case, the unpublished work \cite{har} could be useful.
However, they are few works which really elaborate 
somewhat tedious computations. Moreover, our RDEs 
are slightly more general than standard ones.
We will therefore provide a detailed explanation below.

Throughout this section $T >0$ and 
$\alpha \in (1/4, 1/3]$ and we let $\cV$ and $\cW$ be 
Euclidean spaces.

\subsection{Definition of controlled paths}

First we recall the definition of a controlled path (CP)
with respect to a geometric RP
$X =(1, X^{1}, X^{2}, X^{3})\in G\Omega_{\alpha} (\cV)$.
Let $[a,b]\subset [0,T]$ be a subinterval.
We say that $(Y, Y^{\dagger}, Y^{\dagger\dagger},Y^{\sharp}, Y^{\sharp\sharp})$ is a $\cW$-valued
CP with respect to $X$ on $[a,b]$ if
\begin{align*}
(Y, Y^{\dagger}, Y^{\dagger\dagger},Y^{\sharp}, Y^{\sharp\sharp}) 
&\in 
\cC^\alpha ([a,b], \cW) \times 
\cC^\alpha ([a,b], L(\cV,\cW)) 
\\
&\times \cC^{\alpha} ([a,b], L(\cV, L(\cV,\cW)))
 \times
\cC_{(2)}^{3\alpha} ([a,b],\cW) \times \cC_{(2)}^{2\alpha} ([a,b], L(\cV, \cW))
\end{align*}
and, for all $(s, t) \in \triangle_{[a,b]}$,
\begin{align} 
Y_{t}-Y_{s}
&=Y_{s}^{\dagger} X_{s, t}^{1} +Y_{s}^{\dagger\dagger} X_{s, t}^{2}+Y_{s, t}^{\sharp}, 
\label{def_CP1}\\
Y^{\dagger}_t -Y^{\dagger}_s
&=
Y_{s}^{\dagger\dagger} X_{s, t}^{1}  +Y_{s, t}^{\sharp\sharp}.
\label{def_CP2}
\end{align}
Notice the natural identifications
 $L(\cV, L(\cV,\cW))\cong L^2 (\cV\times \cV; \cW)
\cong L(\cV^{\otimes 2}, \cW)$.
The set of all such CPs with respect to $X$
is denoted by $\cQ^{\alpha}_X ([a,b], \cW)$. 
($X$ is often referred to as a reference RP.)
For simplicity,
$ (Y, Y^{\dagger}, Y^{\dagger\dagger},Y^{\sharp}, Y^{\sharp\sharp})$
will often be written as $(Y, Y^{\dagger}, Y^{\dagger\dagger})$.
Obviously, \eqref{def_CP1} and \eqref{def_CP2} imply that
both $Y^{\sharp}$ and $Y^{\sharp\sharp}$ must vanish on the diagonal.
Note that in fact $X^3$ is not involved in the definition of a CP.

A natural seminorm on $\cQ^{\alpha}_X ([a,b], \cW)$ is defined by
\[
\| (Y, Y^{\dagger}, Y^{\dagger\dagger}) \|_{\cQ^{\alpha}_X, [a,b]}
=\|Y^{\dagger\dagger} \|_{\alpha,[a, b]} 
+ 
\|Y^{\sharp} \|_{3\alpha,[a, b]}
+ 
\|Y^{\sharp\sharp} \|_{2\alpha,[a, b]}
\]
Then, $\cQ^{\alpha}_X ([a,b], \cW)$ is a Banach space with the norm 
\[
 |Y_a|+ |Y^\dagger_a| +|Y^{\dagger\dagger}_a| +
\| (Y, Y^{\dagger}, Y^{\dagger\dagger}) \|_{\cQ^{\alpha}_X, [a,b]}.
\]
(When $[a,b]=[0,T]$, we write 
$\cQ^{\alpha}_X (\cW)$ and $\| \cdot \|_{\cQ^{\alpha}_X}$ 
for simplicity.)
Then, there exist positive constants $C$ and $C^\prime$
depending only on $\alpha$ and $b-a$ such that
\begin{align}  
\|Y^{\dagger} \|_{\alpha,[a, b]} 
&\le
\|Y^{\dagger\dagger} \|_{\infty, [a, b]} \|X^1 \|_{\alpha,[a, b]} +\|Y^{\sharp\sharp} \|_{\alpha,[a, b]}
\nn\\
&\le
\{|Y^{\dagger\dagger}_a| +(b-a)^\alpha \|Y^{\dagger\dagger} \|_{\alpha,[a, b]} \}
\|X^1\|_{\alpha,[a, b]} 
+(b-a)^\alpha\|Y^{\sharp\sharp} \|_{2\alpha,[a, b]}
\nn\\
&\le
C (1+ \|X^1\|_{\alpha} ) (|Y^{\dagger\dagger}_a| +
\| (Y, Y^{\dagger}, Y^{\dagger\dagger}) \|_{\cQ^{\alpha}_X, [a,b]})
\label{ineq.0714-1}
\end{align}
and
\begin{align}  
\|Y\|_{\alpha,[a, b]} 
&\le
\|Y^{\dagger} \|_{\infty, [a, b]} \|X^1 \|_{\alpha,[a, b]} 
+
\|Y^{\dagger\dagger} \|_{\infty, [a, b]} \|X^2 \|_{\alpha,[a, b]} +\|Y^{\sharp} \|_{\alpha,[a, b]}
\nn\\
&\le
\{|Y^{\dagger}_a| +(b-a)^\alpha \|Y^{\dagger} \|_{\alpha,[a, b]} \}
\|X^1\|_{\alpha,[a, b]} 
\nn\\
&\quad
+
\{|Y^{\dagger\dagger}_a| +(b-a)^\alpha \|Y^{\dagger\dagger} \|_{\alpha,[a, b]} \}
(b-a)^{\alpha}\|X^2\|_{2\alpha,[a, b]} 
+(b-a)^{2\alpha} \|Y^{\sharp} \|_{3\alpha,[a, b]}
\nn\\
&\le
C^\prime (1+ \|X^1\|_{\alpha}^2 +\|X^2\|_{2\alpha}  )
 (|Y^{\dagger}_a| +|Y^{\dagger\dagger}_a| +
\| (Y, Y^{\dagger}, Y^{\dagger\dagger}) \|_{\cQ^{\alpha}_X, [a,b]}).
\label{ineq.0714-2}
\end{align}


\begin{example} \label{ex.0715}
Here are a few typical examples of CPs for 
a given RP $X \in G\Omega_{\alpha} (\cV)$. 
(In the first three examples the time interval is $[0,T]$ just for simplicity.
It can be replaced by any subinterval $[a,b]$. )

\begin{enumerate} 
\item
For $\xi \in \cW$, $\sigma\in L(\cV,\cW)$ and
$\eta \in L(\cV, L(\cV,\cW))$,
\[
t \mapsto (\xi + \sigma X^1_{0,t}+ \eta X^2_{0,t}, \,\, 
\sigma+ \eta X^1_{0,t}, \,\,\eta)
\] 
belongs to 
$\cQ^{\alpha}_X (\cW)$.
Note that $\sharp$- and $\sharp\sharp$-components of this CP are zero.
Hence, the $\cQ^{\alpha}_X$-seminorm of this CP is zero, too.

\item
If $\varphi \in \cC^{3\alpha} (\cW)$, 
then obviously
$(\varphi, 0, 0) \in \cQ^{\alpha}_X (\cW)$ with
$\|\varphi\|_{3\alpha} = \| (\varphi, 0, 0) \|_{\cQ^{\alpha}_X}$.
In this way, we have a natural continuous embedding 
$\cC^{3\alpha} (\cW) \hookrightarrow \cQ^{\alpha}_X (\cW)$.

\item
Suppose that $(Y, Y^{\dagger}, Y^{\dagger\dagger})\in \cQ^{\alpha}_X (\cW)$ and 
$g\colon \cW \to \cW'$ is a $C^3$-function from $\cW$ to 
another Euclidean space $\cW'$.
Then $(g(Y), g(Y)^\dagger, g(Y)^{\dagger\dagger}) \in \cQ^{\alpha}_X (\cW')$
if we set
\begin{itemize}
\item
$g(Y)_t:=g(Y_t)$.
\item
$g(Y)^\dagger_t := \nabla g(Y_t) Y^\dagger_t$,
where the right hand side is 
the composition of the two linear maps 
$\nabla g(Y_t)\in L(\cW,\cW')$ and $Y^\dagger_t \in L(\cV,\cW)$.
\item
$g(Y)^{\dagger\dagger}_t := \nabla g(Y_t) Y^{\dagger\dagger}_t 
+\nabla^2 g(Y_t) \la Y^\dagger_t \bullet, Y^\dagger_t \star\ra$,
where the first term on the right hand side is 
the composition of the two linear maps 
$\nabla g(Y_t)\in L(\cW,\cW')$ and $Y^{\dagger\dagger}_t 
\in L(\cV, L(\cV,\cW))$. 
(Note that the second term is symmetric in $\bullet$ and $\star$.)
\end{itemize}

Using Taylor's theorem and \eqref{def_CP1}, we can verify this fact as follows:
\begin{align} 
g(Y_{t})- g(Y_{s}) 
&=\nabla g(Y_{s})\langle Y_{s, t}^{1}\rangle
  + \frac{1}{2}\nabla^2 g(Y_{s})\langle Y_{s, t}^{1}, Y_{s, t}^{1}\rangle
  \nn\\
  &\qquad
   +\frac{1}{2}
   \int_{0}^{1} d \theta \, (1-\theta)^2 
   \, \nabla^{3} g(Y_{s}+\theta Y_{s, t}^{1})
   \langle Y_{s, t}^{1}, Y_{s, t}^{1}, Y_{s, t}^{1}\rangle 
   \nn\\
&=
\nabla g(Y_{s})\langle  Y^{\dagger}_sX_{s, t}^{1} 
 +   Y^{\dagger\dagger}_sX_{s, t}^{2}\rangle
\nn\\
&\qquad
  + \frac{1}{2}\nabla^2 g(Y_{s})\langle Y^{\dagger}_sX_{s, t}^{1}, Y^{\dagger}_sX_{s, t}^{1}\rangle
  + g(Y)_{s, t}^{\sharp},
\label{in.0715-1}
\end{align}
where
\begin{align}
g(Y)_{s, t}^{\sharp}&:=\nabla g (Y_{s}) Y_{s, t}^{\sharp}
+
\nabla^2 g(Y_{s})\langle Y^{\dagger}_sX_{s, t}^{1}, 
  Y^{\dagger\dagger}_s X_{s, t}^{2}\rangle
+
\frac{1}{2}\nabla^2 g(Y_{s})\langle Y^{\dagger\dagger}_sX_{s, t}^{2},
   Y^{\dagger\dagger}_sX_{s, t}^{2}\rangle
   \nn\\
    & \qquad + 
    \nabla^2 g(Y_{s})\langle Y^{\dagger}_sX_{s, t}^{1}+
  Y^{\dagger\dagger}_s X_{s, t}^{2}, \, Y^{\sharp}_{s, t}\rangle
+
\frac{1}{2}\nabla^2 g(Y_{s})\langle Y^{\sharp}_{s, t},
   Y^{\sharp}_{s,t}\rangle
\nn\\
&\qquad
+\frac12
\int_{0}^1 d \theta(1-\theta)^2
\nabla^{3} g(Y_{s}+\theta Y_{s, t}^{1})\langle Y_{s, t}^{1}, Y_{s, t}^{1}, Y_{s, t}^{1}\rangle.
\label{in.0715-2}
\end{align}
It is easy to see from \eqref{in.0715-2}
that $g(Y)^{\sharp}\in \cC_{(2)}^{3\alpha} (\cW')$. 
Due to the shuffle relation \eqref{eq.0711-2} and 
the symmetry of the bilinear mapping,  it holds that
\[
\frac12 \nabla^2 g(Y_{s})\langle Y^{\dagger}_sX_{s, t}^{1}, Y^{\dagger}_sX_{s, t}^{1}\rangle
=\nabla^2 g(Y_t) \la Y^\dagger_t \bullet, Y^\dagger_t \star\ra
\vert_{( \bullet, \star)=X^2_{s,t}}.
\]
Here, the right hand should be understood as $X^2_{s,t}$ being plugged into 
an element of $L(\cV^{\otimes 2}, \cW^\prime)$.
Thus,  the composition $g (Y)$ satisfies \eqref{def_CP1}.

Keeping \eqref{def_CP1} and \eqref{def_CP2} in mind, we can see that 
\begin{align}
g(Y)^\dagger_t - g(Y)^\dagger_s
&= \nabla g(Y_t) Y^\dagger_t -  \nabla g(Y_s) Y^\dagger_s
\nn\\
&=
\nabla g(Y_s) Y^\dagger_t +\nabla^2 g(Y_s) \la Y^1_{s,t}, Y^\dagger_t\ra 
\nn\\
&\quad
+
\int_{0}^1 d \theta(1-\theta)
\nabla^{3} g(Y_{s}+\theta Y_{s, t}^{1})
\langle Y_{s, t}^{1}, Y_{s, t}^{1}, Y^\dagger_t \rangle
-  \nabla g(Y_s) Y^\dagger_s
\nn\\
&=\nabla g(Y_s) \la Y^{\dagger\dagger}_s X_{s, t}^{1}\ra
  +\nabla^2 g(Y_s) \la Y^\dagger_sX^1_{s,t}, Y^\dagger_s\ra  
   + g(Y)_{s, t}^{\sharp\sharp},
\label{in.0715-3}
\end{align}
where 
\begin{align}
g(Y)_{s, t}^{\sharp\sharp}
&:=
\nabla g(Y_s) \la Y_{s, t}^{\sharp\sharp}\ra
\nn\\
&\qquad
+
\nabla^2 g(Y_s) \la Y_{s}^{\dagger\dagger} X_{s, t}^{2}+Y_{s, t}^{\sharp},
     \, Y_{s}^{\dagger} \ra
\nn\\
&\qquad
+
\nabla^2 g(Y_s) \la Y_{s}^{\dagger} X_{s, t}^{1}+
Y_{s}^{\dagger\dagger} X_{s, t}^{2}+Y_{s, t}^{\sharp}, \,
      Y_{s}^{\dagger\dagger} X_{s, t}^{1}  +Y_{s, t}^{\sharp\sharp}
       \ra
\nn\\
&\qquad
+
\int_{0}^1 d \theta(1-\theta)
\nabla^{3} g(Y_{s}+\theta Y_{s, t}^{1})
\langle Y_{s, t}^{1}, Y_{s, t}^{1}, Y^\dagger_t \rangle.
\label{in.0715-4}
\end{align}
It is easy to see from \eqref{in.0715-4}
that $g(Y)^{\sharp\sharp}\in \cC_{(2)}^{2\alpha} (\cW')$. 
Thus, the composition satisfies \eqref{def_CP2}, too.

\item
Concatenation of two CPs is also a CP.
Let $0\le a <b<c \le T$.
For $(Y, Y^{\dagger}, Y^{\dagger\dagger})\in \cQ^{\alpha}_X ([a,b], \cW)$ and
$(\hat{Y}, \hat{Y}^{\dagger}, \hat{Y}^{\dagger\dagger})\in \cQ^{\alpha}_X ([b,c], \cW)$
with $(Y_b, Y^{\dagger}_b, Y^{\dagger\dagger}_b)
 = (\hat{Y}_b, \hat{Y}^{\dagger}_b, \hat{Y}^{\dagger\dagger}_b)$,
their concatenation 
$(Z, Z^\dagger, Z^{\dagger\dagger}):=
(Y*\hat{Y}, Y^{\dagger}*\hat{Y}^{\dagger}, 
Y^{\dagger\dagger}* \hat{Y}^{\dagger\dagger})$ can naturally be defined and 
belongs to $\cQ^{\alpha}_X ([a,c], \cW)$.
(Here, $*$ stands for the usual concatenation operation 
for two continuous paths.)

It is clear that $Z, Z^{\dagger}, Z^{\dagger\dagger}$ are 
$\alpha$-H\"older continuous on $[a,c]$.
To prove that $Z^{\sharp},  Z^{\sharp\sharp}\in\cC_{(2)}^{2\alpha} (\cW)$,
it is sufficient to observe the following: 
For $a\le s \le b \le t \le c$, we have
\begin{align} 
Z_{s, t}^{\sharp\sharp} &=Z^{\dagger,1}_{s, t}
-Z_{s}^{\dagger\dagger} X_{s, t}^{1} 
\nn\\
&=(Z^{\dagger,1}_{s, b} -Z_{s}^{\dagger\dagger} X_{s,b}^{1})
+(Z^{\dagger, 1}_{b, t} -Z_{b}^{\dagger\dagger} X_{b, t}^{1})
+Z_{s, b}^{\dagger\dagger, 1} X_{b, t}^{1}
\nn\\
&=
Y^{\sharp\sharp}_{s, b}
 + \hat{Y}^{\sharp\sharp}_{b, t}
   +Y_{s, b}^{\dagger\dagger, 1} X_{b, t}^{1}.
\label{in.0715-5}
\end{align}
and
\begin{align} 
Z_{s, t}^{\sharp} &=Z^1_{s,t}
-Z_{s}^{\dagger} X_{s, t}^{1} -Z_{s}^{\dagger\dagger} X_{s, t}^{2}
\nn\\
&=
(Z^1_{s,b}+Z^1_{b,t})
-\{ Z_{s}^{\dagger} X_{s, b}^{1}+Z_{b}^{\dagger} X_{b, t}^{1} 
- Z_{s,b}^{\dagger, 1} X_{b, t}^{1}
\}
\nn\\
& \qquad
-\{
Z_{s}^{\dagger\dagger} X_{s, b}^{2}+Z_{b}^{\dagger\dagger} X_{b, t}^{2}
- Z_{s,b}^{\dagger\dagger, 1} X_{b, t}^{2}
+Z_{s}^{\dagger\dagger} X_{s, b}^{1}\otimes X_{b, t}^{1}
\}
\nn\\
&=
Y_{s, b}^{\sharp}+\hat{Y}_{b, t}^{\sharp} +Y_{s,b}^{\dagger\dagger, 1} X_{b, t}^{2}
+ Y_{s, b}^{\sharp\sharp}  X_{b, t}^{1}.
\label{in.0715-6}
\end{align}
Here, we have used \eqref{def_CP1}--\eqref{def_CP2} and Chen's relation for $X$.
The right hand side of \eqref{in.0715-5} and \eqref{in.0715-6}
are clearly dominated by a constant multiple of 
$(t-s)^{2\alpha}$ and of $(t-s)^{3\alpha}$, respectively.
\end{enumerate}
\end{example}

\subsection{Integration of controlled paths against rough paths}

Next we discuss integration of a CP 
$(Y, Y^{\dagger}, Y^{\dagger\dagger})\in \cQ^{\alpha}_X ([a,b], L(\cV, \cW))$
against a reference RP $X \in G\Omega_{\alpha} (\cV)$,
where $[a, b]\subset [0, T]$.
Note that $Y^{\dagger}$ and $Y^{\dagger\dagger}$ take values in 
\[
L ( \cV, L(\cV,\cW))\cong L^{(2)} (\cV \times\cV, \cW)
\cong
L ( \cV^{\otimes 2},\cW),
\]
and
\[
L(\cV, L ( \cV, L(\cV,\cW))) \cong L^{(3)} (\cV\times \cV \times\cV, \cW)
\cong
L ( \cV^{\otimes 3},\cW),
\]
respectively.
Note also that $Y^{\sharp}$ and $Y^{\sharp\sharp}$  take values in 
$L ( \cV,\cW)$ and $L ( \cV^{\otimes 2},\cW)$, respectively.

First, we define 
\[
J_{s, t}=Y_{s} X_{s, t}^{1}+Y_{s}^{\dagger} X_{s, t}^{2}+Y_{s}^{\dagger\dagger} X_{s, t}^{3},
\qquad
(s,t)\in \triangle_{[a,b]}.
\]
From Chen's relation for $X$ and 
\eqref{def_CP1}--\eqref{def_CP2}, it easily follows that
\begin{equation} \label{eq.0411-3}
J_{s, u}+J_{u, t}-J_{s, t}
=Y_{s, u}^{\sharp} X_{u, t}^{1}
+Y_{s, u}^{\sharp\sharp} X_{u, t}^{2}
 +Y_{s,u}^{\dagger\dagger, 1} X_{u, t}^{3},
\qquad 
a\le s\le u \le t \le b,
\end{equation}
where we set $Y_{s,u}^{\dagger\dagger, 1}= Y_{u}^{\dagger\dagger}-Y_{s}^{\dagger\dagger}$.
Note that the right hand side above is dominated 
by a constant multiple of $(t-s)^{4\alpha}$.

Let $\mathcal{P}=\{s=t_{0}<t_{1}<\cdots<t_{N}=t\}$
be a partition of $[s,t] \subset [a,b]$.
Its mesh size is denoted by $|\mathcal{P}|$.
We define 
$J_{s, t}(\mathcal{P})=\sum_{i=1}^{N} J_{t_{i-1}, t_{i}}$ and
\begin{equation} \label{eq.0411-4}
\int_{s}^{t} Y_{u} d X_{u}=\lim_{|\mathcal{P}| \searrow 0} J_{s, t}(\mathcal{P}),
\qquad (s,t)\in \triangle_{[a,b]}.
\end{equation}
The limit above is known to exist, and is called a RP integral.
It will turn out in the next proposition
that an (indefinite) RP integral against $X$ is again a CP with respect to $X$.
By the way it is defined,
this RP integral clearly has additivity with respect to the interval $[s,t]$.

\begin{proposition} \label{prop.gousei}
Let  $\tfrac14 < \alpha \le \tfrac13$ and $[a,b] \subset [0,T]$.
Suppose that $X \in G\Omega_{\alpha} (\cV)$
and $(Y, Y^{\dagger}, Y^{\dagger\dagger})\in \cQ^{\alpha}_X ([a,b], L(\cV, \cW))$. Then, the limit in \eqref{eq.0411-4} exists for all $(s,t)\in \triangle_{[a,b]}$.
Moreover, we have
\begin{equation}
\Bigl(\int_{a}^{\cdot} Y_{u} d X_{u}, \,Y, \,  Y^{\dagger} \Bigr)
\in \cQ^{\alpha}_X ([a,b], \cW)
\label{eq.0413-3}
\end{equation}
with the following estimate: 
\begin{align} 
\lefteqn{
\Bigl|
\int_{s}^{t} Y_{u} d X_{u} -(Y_{s} X_{s, t}^{1}+Y_{s}^{\dagger} X_{s, t}^{2}
 +Y_{s}^{\dagger\dagger} X_{s, t}^{3})
\Bigr|_\cW
}
\nn \\
&\le 
\kappa_{\alpha} (t-s)^{4\alpha} 
(\|Y^{\sharp} \|_{3\alpha,[a, b]} \|X^1\|_{\alpha,[a, b]}
+ 
\|Y^{\sharp\sharp} \|_{2\alpha,[a, b]} \|X^2\|_{2\alpha,[a, b]}
\nn\\
&\qquad\qquad\qquad\qquad \qquad
+ 
\|Y^{\dagger\dagger} \|_{\alpha,[a, b]} \|X^3\|_{3\alpha,[a, b]}),
\qquad (s,t)\in \triangle_{[a,b]}.
\label{in.0413-2}
\end{align}
Here, we set $\kappa_{\alpha}=2^{4\alpha}\zeta (4\alpha)$
with $\zeta$ being the Riemann zeta function.
\end{proposition}

\begin{proof} 
In this proof the norm of $\cW$ is denoted by $|\cdot|$.
First we prove the convergence. 
For $\mathcal{P}$ given as above,
we can find $i ~(1 \le i \le N-1)$ such that $t_{i+1}-t_{i-1}\le 2(t-s)/(N-1)$.
Then, we see that
\begin{align} 
|J_{s, t}(\mathcal{P}) - J_{s, t}(\mathcal{P}\setminus \{t_i\})|
&=
|J_{t_{i-1}, t_{i}} +  J_{t_{i}, t_{i+1}}  -  J_{t_{i-1}, t_{i+1}}|
\nn\\
&=
|Y_{t_{i-1}, t_{i}}^{\sharp} X_{t_{i}, t_{i+1}}^{1}
+Y_{t_{i-1}, t_{i}}^{\sharp\sharp}X_{t_{i}, t_{i+1}}^{2}
+Y_{t_{i-1}, t_{i}}^{\dagger\dagger, 1}X_{t_{i}, t_{i+1}}^{3}
|
\nn\\
&\le
(\|Y^{\sharp} \|_{3\alpha,[a, b]} \|X^1\|_{\alpha,[a, b]}
+ 
\|Y^{\sharp\sharp} \|_{2\alpha,[a, b]} \|X^2\|_{2\alpha,[a, b]}
\nn\\
&\qquad\qquad\qquad
+ 
\|Y^{\dagger\dagger} \|_{\alpha,[a, b]} \|X^3\|_{3\alpha,[a, b]})
 \Bigl(\frac{2(t-s)}{N-1}\Bigr)^{4\alpha}.
 \nn
\end{align}

Extracting points one by one from $\mathcal{P}$ in this way
until the partition  becomes the trivial one $\{s,t\}$, we have
\begin{align} 
|J_{s, t}(\mathcal{P}) - J_{s, t}|
&\le
\kappa_{\alpha} (t-s)^{4\alpha}
(\|Y^{\sharp} \|_{3\alpha,[a, b]} \|X^1\|_{\alpha,[a, b]}
\nn\\
&\qquad
+ 
\|Y^{\sharp\sharp} \|_{2\alpha,[a, b]} \|X^2\|_{2\alpha,[a, b]}
+ 
\|Y^{\dagger\dagger} \|_{\alpha,[a, b]} \|X^3\|_{3\alpha,[a, b]}).
\label{in.0413-1}
\end{align}
Note that the condition $4\alpha >1$ is used here.
By standard argument in RP theory, 
\eqref{in.0413-1} implies that
$\{J_{s, t}(\mathcal{P}) \}_{\mathcal{P}}$ is Cauchy as 
$|\mathcal{P}| \searrow 0$.
Therefore, the limit in \eqref{eq.0411-4} exists and \eqref{in.0413-2} holds.

Since it obviously holds that
\[
|Y_{s}^{\dagger\dagger} X_{s, t}^{3}|
\le 
\{ |Y_{a}^{\dagger\dagger}|+(b-a)^\alpha \|Y^{\dagger\dagger} \|_{\alpha,[a, b]}\}
\|X^3\|_{3\alpha,[a, b]} (t-s)^{3\alpha}
\]
and
\begin{align*}
\lefteqn{
|Y_t-Y_s - Y_{s}^{\dagger} X_{s, t}^{1}|
}
\nn\\
&\le 
|Y_{s}^{\dagger\dagger} X_{s, t}^{2}|+ |Y_{s, t}^{\sharp}|
\nn\\
&\le
\Bigl\{ (
|Y_{a}^{\dagger\dagger}|+(b-a)^\alpha \|Y^{\dagger\dagger} \|_{\alpha,[a, b]}
)
\|X^2\|_{2\alpha,[a, b]} 
+
 (b-a)^{\alpha}
\|Y^{\sharp\sharp} \|_{2\alpha,[a, b]}  \Bigr\} (t-s)^{2\alpha},
\end{align*}
the path
$(\int_{a}^{\cdot} Y_{u} d X_{u}, \,Y, \,  Y^{\dagger})$
satisfies \eqref{def_CP1}--\eqref{def_CP2}.
Hence, \eqref{in.0413-2} implies \eqref{eq.0413-3}.
\end{proof}

\subsection{Rough differential equations
with bounded and globally Lipschitz drift}\label{sec.rde}

Now we discuss RDEs in the framework of controlled path theory.
We basically follow \cite{har, boge}, but our RDE has a drift term.
In this subsection, $\cS$ is a metric space and 
we assume
$\tfrac14 < \beta < \alpha \le \tfrac13$ and
$X \in G\Omega_{\alpha} (\cV)\subset G\Omega_{\beta} (\cV)$.

We set conditions on the coefficients of our RDE. 
Let $\sigma\colon \cW \to L(\cV,\cW)$ be of $C^4_{{\rm b}}$
and let $f\colon \cW \times \cS \to \cW$ be a continuous map
satisfying the following condition:
\begin{equation}\label{cond.0506-1}
 \sup_{y \in \cW, z\in \cS} |f(y,z)|_{\cW}
+\sup_{y,y' \in \cW, y\neq y',z\in \cS} 
\frac{|f(y,z)-f(y',z) |_{\cW}}{|y-y'|_{\cW}} <\infty.
\end{equation}
The first and the second term above will be denoted by 
$\|f\|_{\infty}$ and $L_f$, respectively.

For an $\cS$-valued continuous path $\psi \colon [0,T] \to \cS$, 
we consider the following RDE driven by $X$ with the initial value
 $\xi \in\cW$: For all $t \in [0,T]$,
\begin{align} 
Y_{t}
=\xi +\int_{0}^{t} f(Y_{s}, \psi_{s}) ds
+\int_{0}^{t} \sigma(Y_{s}) d X_{s},
\qquad
Y^\dagger_t = \sigma(Y_t),  
\quad 
Y^{\dagger\dagger}_t = \nabla \sigma(Y_t)Y^\dagger_t.
\label{rde.0413}
\end{align}
For every $(Y, Y^{\dagger}, Y^{\dagger\dagger})\in \cQ^{\beta}_X (\cW)$,
the right hand side of this system of equations also
belongs to $\cQ^{\beta}_X (\cW)$, due to Example \ref{ex.0715} and
Proposition \ref{prop.gousei}.
Therefore, \eqref{rde.0413} should be understood as
 an equality in $\cQ^{\beta}_X (\cW)$.
(Following the preceding works \cite{har, boge}, 
we slightly relax the H\"older topology of the space of CPs for quick proofs.)

One should note that if $(Y, Y^{\dagger}, Y^{\dagger\dagger})\in \cQ^{\beta}_X (\cW)$ solves \eqref{rde.0413},
$Y^\dagger_t = \sigma(Y_t)$ and 
$Y^{\dagger\dagger}_t = \nabla \sigma\cdot \sigma (Y_t)$,
in particular, $Y^\dagger$ and $Y^{\dagger\dagger}$
are completely determined by the first level part $Y$.
Precisely, $\nabla \sigma \cdot \sigma (y)$ above 
stands for the element of $L (\cV\otimes \cV, \cW)$ defined by
\[
\cV\otimes \cV \ni v\otimes v^\prime 
\mapsto
\nabla \sigma(y)\la \sigma(y)v, v^\prime\ra \in \cW.
\]

Here we state our main result in this section.
\begin{proposition} \label{prop.0429}
Let the assumptions be as above.
Then, for every $X\in G\Omega_{\alpha} (\cV)$, $\xi \in \cW$ 
and $\psi$, there exists a unique global solution 
$(Y, Y^{\dagger}, Y^{\dagger\dagger})\in \cQ^{\beta}_X (\cW)$ of RDE \eqref{rde.0413}.
Moreover, it satisfies the following estimate:
there exist positive constants $c$ and $\nu$ 
independent of $X,  \xi, \psi, \sigma, f$ such that 
\[
\|Y\|_{\beta} \le c \{(K+1) (\vertiii{X}_{\alpha} +1)  \}^\nu,
\qquad 
X\in G\Omega_{\alpha} (\cV). 
\]
Here, we set
$K := \| \sigma\|_{C_b^4}\vee \|f\|_{\infty}\vee L_f$. 
\end{proposition}

\begin{proof} 
In this proof the norm of finite-dimensional spaces  
is denoted by $|\cdot|$.
Without loss of generality we may assume $T=1$.
Let $\tau \in (0,1]$ and $\xi \in \cW$. 
For simplicity, we write $\eta (y) := \nabla \sigma \cdot \sigma (y)$.
We denote by $c_i >0$ and $\nu_i\in \N$ $(i=1,2,\ldots)$
 certain constants independent of $X,  \xi, \tilde{\xi}, \psi, \sigma, f,
\tau$ and $(s,t) \in \triangle_\tau$.

We define
$\cM^{\xi}_{[0, \tau]}, \cM^1_{[0, \tau]}, \cM^2_{[0, \tau]} \colon \cQ^{\beta}_X ([0,\tau], \cW) \to 
\cQ^{\beta}_X ([0,\tau], \cW)$ by 
\begin{align} 
\cM^{1}_{[0, \tau]} (Y, Y^\dagger, Y^{\dagger\dagger}) 
&=
\Bigl( \int_{0}^{\cdot} \sigma(Y_{s}) d X_{s}, \,\, \sigma(Y),
\,\, \nabla \sigma(Y) Y^\dagger
\Bigr),
\nn\\
\cM^{2}_{[0, \tau]} (Y, Y^\dagger,Y^{\dagger\dagger}) 
&= \Bigl(\int_{0}^{\cdot} f(Y_{s}, \psi_{s}) ds,
\,\, 0,
\,\, 0  \Bigr),
\nn\\
\cM^\xi_{[0, \tau]} (Y, Y^\dagger,Y^{\dagger\dagger}) 
&=(\xi, 0, 0) 
   +\cM^{1}_{[0, \tau]} (Y, Y^\dagger,Y^{\dagger\dagger}) 
     + \cM^{2}_{[0, \tau]} (Y, Y^\dagger,Y^{\dagger\dagger}) .
\end{align}
If $(Y, Y^\dagger,Y^{\dagger\dagger})$ 
starts at $(\xi, \sigma(\xi), \eta(\xi))$, so does 
$\cM^\xi_{[0, \tau]}(Y, Y^\dagger, Y^{\dagger\dagger})$.
A fixed point of $\cM^\xi_{[0, \tau]}$ is a solution of RDE \eqref{rde.0413}
on the interval $[0,\tau]$.

We also set 
\begin{align*}
B_{[0, \tau]}^\xi
=\{(Y, Y^\dagger, Y^{\dagger\dagger}) \in \cQ^{\beta}_X ([0,\tau], \cW)
&\mid
\,
\|(Y, Y^\dagger, Y^{\dagger\dagger})\|_{ \cQ^{\beta}_X,  [0, \tau]} \le  1,  
\nn\\
&
  \qquad Y_0=\xi, \,\, Y^\dagger_0=\sigma(\xi),  \,\,
  Y^{\dagger\dagger}_0= \eta (\xi)\}.
\end{align*}
This set is something like a ball of radius $1$ centered at 
\[
t \mapsto 
\bigl(
\xi + \sigma (\xi) X^1_{0,t}+ \eta (\xi) X^2_{0,t}, \,\, 
\sigma (\xi)+  \eta (\xi) X^1_{0,t}, \,\,  \eta (\xi)
\bigr)
\]
(see Example \ref{ex.0715}).
Since the initial point $(Y_0, Y^\dagger_0, Y^{\dagger\dagger}_0)$ is fixed, 
$\|\cdot\|_{ \cQ^{\beta}_X,  [0, \tau]}$ defines a distance on this set.


For a while from now, we will work only on $[0,\tau]$ and therefore
omit $[0,\tau]$ from the subscript for notational simplicity.
We will often write  $\varphi_{s, t}^{1}:=\varphi_t -\varphi_s$ for a 
usual path $\varphi$ that takes values in a vector space.

Let $\xi, \tilde{\xi}\in \cW$ and pick 
$(Y, Y^\dagger, Y^{\dagger\dagger})\in B^\xi$ and $(\tilde{Y}, \tilde{Y}^\dagger, \tilde{Y}^{\dagger\dagger}) \in B^{\tilde\xi}$ arbitrarily.
For simplicity we write $\Delta := \|  (Y, Y^\dagger, Y^{\dagger\dagger})-(\tilde{Y},  \tilde{Y}^\dagger, \tilde{Y}^{\dagger\dagger}) \|_{\cQ^{\beta}_X}$.
One can show the following estimates 
for all $(s,t) \in \triangle_\tau$ by repeatedly 
using \eqref{def_CP1} and \eqref{def_CP2}:
\begin{align} 
\|Y^{\dagger\dagger}\|_{\infty} &\le
   |\eta(\xi)|+\sup _{s \leq \tau} |Y_{s}^{\dagger\dagger}-Y_{0}^{\dagger\dagger}| 
\le
 K^2 +\|Y^{\dagger\dagger}\|_{\beta} \tau^{\beta} \le  K^2+1,
\label{sono11}\\
\|Y^{\dagger\dagger}- \tilde{Y}^{\dagger\dagger}\|_{\infty} 
 &\le
  | \eta(\xi)-\eta(\tilde\xi)|
   +\|Y^{\dagger\dagger}- \tilde{Y}^{\dagger\dagger} \|_{\beta} \tau^{\beta}
 \le 
 2K^2 |\xi -\tilde\xi| +\Delta,
\label{sono12}\\
|Y_{s, t}^{\dagger,1}| &\le
   |Y_{s}^{\dagger\dagger} X_{s, t}^{1}|+|Y_{s, t}^{\sharp\sharp}|
 \nn\\ 
  &\le (K^2 +1)\|X^{1}\|_{\alpha} (t-s)^{\alpha}  
   +\|Y^{\sharp\sharp}\|_{2 \beta}(t-s)^{2 \beta}
   \nn\\
     &\le   (K^2 +1)(\|X^{1}\|_{\alpha}+1)  (t-s)^{\alpha},
   \label{sono13}\\
|Y_{s, t}^{\dagger,1}- \tilde{Y}_{s, t}^{\dagger,1}| &\le
   |(Y_{s}^{\dagger\dagger} - \tilde{Y}_{s}^{\dagger\dagger})X_{s, t}^{1}|
     +|Y_{s, t}^{\sharp\sharp}-\tilde{Y}_{s, t}^{\sharp\sharp} |
     \nn\\
       &\le 
        (2K^2 |\xi -\tilde\xi| + 
        \Delta)
        \|X^{1}\|_{\alpha} (t-s)^{\alpha}   
          +\|Y^{\sharp\sharp} -  \tilde{Y}^{\sharp\sharp} \|_{2 \beta}(t-s)^{2 \beta}
            \nn\\
       &\le  \bigl\{ 2K^2 \|X^{1}\|_{\alpha} |\xi -\tilde\xi| 
       +  (\|X^{1}\|_{\alpha} +1) \Delta  
         \bigr\}     (t-s)^{\alpha},
     \label{sono14}
     \\
     \|Y^{\dagger}\|_{\infty} 
      &\le
   |\sigma (\xi)|+\sup _{s \le \tau} |Y_{0, s}^{\dagger, 1}| 
 \nn\\
       &\le
       K+ (K^2 +1)(\|X^{1}\|_{\alpha}+1)
\le  2(K^2 +1)(\|X^{1}\|_{\alpha}+1),
\label{sono15}
\\
     \|Y^{\dagger}- \tilde{Y}^{\dagger}\|_{\infty} 
      &\le
   |\sigma (\xi) - \sigma (\tilde{\xi})|
      +\sup _{s \le \tau} |Y_{0, s}^{\dagger, 1}
           -\tilde{Y}_{0, s}^{\dagger, 1}| 
             \nn\\
             &\le
                (K+ 2K^2 \|X^{1}\|_{\alpha})|\xi-\tilde{\xi}| 
                +( \|X^{1}\|_{\alpha} +1)\Delta.
\label{sono16}
     \end{align}
From these, we also have
\begin{align}  
|Y_{s, t}^{1}| 
&\le
   |Y_{s}^{\dagger} X_{s, t}^{1}|
   + |Y_{s}^{\dagger\dagger} X_{s, t}^{2}|
   +|Y_{s, t}^{\sharp}|
 \nn\\ 
  &\le
  2(K^2 +1)(\|X^{1}\|_{\alpha}+1)\|X^{1}\|_{\alpha} (t-s)^{\alpha} 
   \nn\\
    &\qquad\quad
   +(K^2 +1)\|X^{2}\|_{2\alpha} (t-s)^{2\alpha} 
    +\|Y^{\sharp}\|_{3 \beta}(t-s)^{3 \beta}
     \nn\\
      &\le
        4(K^2 +1) (\vertiii{X}_\alpha^2 +1)   (t-s)^{\alpha},
\label{sono17}
\\
|Y_{s, t}^{1} - \tilde{Y}_{s, t}^{1} | 
&\le
   |(Y_{s}^{\dagger} - \tilde{Y}_{s}^{\dagger}) X_{s, t}^{1}|
   + |(Y_{s}^{\dagger\dagger} -\tilde{Y}_{s}^{\dagger\dagger})X_{s, t}^{2}|
   +|Y_{s, t}^{\sharp}- \tilde{Y}_{s, t}^{\sharp}|
 \nn\\ 
  &\le
     \|Y^{\dagger}- \tilde{Y}^{\dagger}\|_{\infty} 
     \|X^{1}\|_{\alpha} (t-s)^{\alpha} 
      \nn\\
       &\qquad \quad + 
      \|Y^{\dagger\dagger}- \tilde{Y}^{\dagger\dagger}\|_{\infty} 
     \|X^{2}\|_{2\alpha} (t-s)^{2\alpha} 
        + 
        \|Y^{\sharp} -  \tilde{Y}^{\sharp} \|_{3 \beta}(t-s)^{3 \beta}
            \nn\\ 
  &\le
   \bigl\{ 
     (K+2 K^2 \|X^{1}\|_{\alpha})|\xi-\tilde{\xi}| \|X^{1}\|_{\alpha}
                +(1+ \|X^{1}\|_{\alpha} )\Delta \|X^{1}\|_{\alpha}
\nn\\
&\qquad \quad +
 2K^2 |\xi -\tilde\xi| \|X^{2}\|_{2\alpha}
 +
 \Delta \|X^{2}\|_{2\alpha}
  +\Delta
   \bigr\} (t-s)^{\alpha}
     \nn\\
      &\le 
      2\bigl\{
      (K^2+1) (\vertiii{X}_\alpha^2+1)  |\xi-\tilde{\xi}|+
       (\vertiii{X}_\alpha^2+1)\Delta
      \bigr\}
       (t-s)^{\alpha}.
  \label{sono18}
     \end{align}
Note that $Y^{\dagger}$ and $Y$ are in fact $\alpha$-H\"older continuous.
Hence, if $\tau$ is small, the $\beta$-H\"older seminorms of 
$Y^{\dagger}$, $Y^{\dagger}-\tilde{Y}^{\dagger}$,
$Y$ and $Y-\tilde{Y}$ can be made very small. 
From \eqref{sono18} we can easily see that
 \begin{align}
 \| Y- \tilde{Y} \|_\infty 
  &\le  |Y_0- \tilde{Y}_0| 
      + \sup_{s\le \tau} | Y^1_{0,s}- \tilde{Y}^1_{0,s} |
   \nn\\
     &\le 
     \{1 + 2(K^2+1) (\vertiii{X}_\alpha^2+1) \tau^\alpha\}
      |\xi -\tilde\xi|
    +
      2(\vertiii{X}_\alpha^2+1)\Delta 
          \tau^{\alpha}.      \label{sono5}      
     \end{align}
%

Next we calculate
$(\sigma(Y), \sigma(Y)^{\dagger}, \sigma(Y)^{\dagger\dagger})$,
whose explicit form was provided in the third item of Example \ref{ex.0715}.
We can easily see from \eqref{sono11}, \eqref{sono13},
\eqref{sono15} and  \eqref{sono17} that 
\begin{align} 
\lefteqn{
|\sigma (Y)^{\dagger\dagger}_t -\sigma (Y)^{\dagger\dagger}_s|
}
\nn\\
&\le 
|\nabla \sigma (Y_t) Y^{\dagger\dagger}_t -\nabla \sigma (Y_s) Y^{\dagger\dagger}_s|
+
|\nabla^2 \sigma(Y_t) \la Y^\dagger_t \bullet, Y^\dagger_t \star\ra
-\nabla^2 \sigma(Y_s) \la Y^\dagger_s \bullet, Y^\dagger_s \star\ra
|
\nn\\
&\le 
\|\nabla^2 \sigma \|_\infty |Y^1_{s,t}| \| Y^{\dagger\dagger} \|_\infty
+
\|\nabla \sigma \|_\infty |Y^{\dagger\dagger, 1}_{s,t}|
\nn\\
&\qquad\quad
+
\|\nabla^3 \sigma \|_\infty |Y^1_{s,t}| \| Y^{\dagger} \|_\infty^2
+
2\|\nabla^2 \sigma \|_\infty |Y^{\dagger, 1} _{s,t}| \| Y^{\dagger} \|_\infty
\nn\\
&\le 
c_1 (K+1)^7 (\vertiii{X}_\alpha+1)^4 (t-s)^\beta.
\label{ineq.913-1}
\end{align}
From the explicit form of $\sigma (Y)_{s, t}^{\sharp\sharp}$
in \eqref{in.0715-4}, we have
\begin{align}
|\sigma (Y)_{s, t}^{\sharp\sharp}|
&\le
|\nabla \sigma (Y_s) \la Y_{s, t}^{\sharp\sharp}\ra|
\nn\\
&\qquad
+
|\nabla^2 \sigma (Y_s) \la Y_{s}^{\dagger\dagger} X_{s, t}^{2}+Y_{s, t}^{\sharp},
     \, Y_{s}^{\dagger} \ra|
\nn\\
&\qquad
+
|\nabla^2 \sigma (Y_s) \la 
Y_{s}^{\dagger} X_{s, t}^{1}+Y_{s}^{\dagger\dagger} X_{s, t}^{2}+Y_{s, t}^{\sharp}, \,
      Y_{s}^{\dagger\dagger} X_{s, t}^{1}  +Y_{s, t}^{\sharp\sharp}
       \ra|
\nn\\
&\qquad
+
\Bigl|\int_{0}^1 d \theta(1-\theta)
\nabla^{3} \sigma (Y_{s}+\theta Y_{s, t}^{1})
\langle Y_{s, t}^{1}, Y_{s, t}^{1}, Y^\dagger_t \rangle\Bigr|
\nn\\
&\le
c_2  (K+1)^7 (\vertiii{X}_\alpha+1)^5 (t-s)^{2\beta},
\label{ineq.913-2}
\end{align}
where we have used \eqref{sono11}, \eqref{sono13},
\eqref{sono15} and  \eqref{sono17} again.
Similarly, we see from \eqref{in.0715-2} that
\begin{align}
|\sigma (Y)_{s, t}^{\sharp}|
&\le 
|\nabla \sigma (Y_{s}) Y_{s, t}^{\sharp}|
+
|\nabla^2 \sigma (Y_{s})\langle Y^{\dagger}_sX_{s, t}^{1}, 
  Y^{\dagger\dagger}_s X_{s, t}^{2}\rangle|
\nn\\
    & \qquad +
\frac{1}{2}  
|\nabla^2 \sigma (Y_{s})\langle Y^{\dagger\dagger}_sX_{s, t}^{2},
   Y^{\dagger\dagger}_sX_{s, t}^{2}\rangle|
   \nn\\
    & \qquad + 
    |\nabla^2 \sigma (Y_{s})\langle Y^{\dagger}_sX_{s, t}^{1}+
  Y^{\dagger\dagger}_s X_{s, t}^{2}, \, Y^{\sharp}_{s, t}\rangle|
+
\frac{1}{2}  |\nabla^2 \sigma (Y_{s})\langle Y^{\sharp}_{s, t},
   Y^{\sharp}_{s,t}\rangle|
\nn\\
&\qquad
+\frac12
\Bigl|  \int_{0}^1 d \theta(1-\theta)^2
\nabla^{3} \sigma (Y_{s}+\theta Y_{s, t}^{1})\langle Y_{s, t}^{1}, Y_{s, t}^{1}, Y_{s, t}^{1}\rangle \Bigr|
\nn\\
&\le
c_3  (K+1)^7 (\vertiii{X}_\alpha+1)^6 (t-s)^{3\beta},
\label{ineq.913-3}
\end{align}
where we have used \eqref{sono11}, \eqref{sono13},
\eqref{sono15} and  \eqref{sono17} again.
It immediately follows from \eqref{ineq.913-1}--\eqref{ineq.913-3} that
\begin{equation} \label{ineq.913-4}
\|(\sigma(Y), \sigma(Y)^{\dagger}, \sigma(Y)^{\dagger\dagger}) 
\|_{\cQ^{\beta}_X}
\le
c_4  (K+1)^7 (\vertiii{X}_\alpha+1)^6.
\end{equation}

We will check that $\cM^\xi$ maps $B^\xi$ to itself 
if $\tau$ is small enough.
As before, we simply write $\nabla\sigma (y) \la y' \ra$ 
for $\nabla \sigma (y) \la y', \bullet \ra = \nabla_{y'} \sigma (y) \in L(\cV, \cW)$
when $y, y' \in \cW$.
It immediately follows from Example \ref{ex.0715} that
\begin{equation}
\|\cM^{2} (Y, Y^\dagger, Y^{\dagger\dagger}) \|_{\cQ^{\beta}_X}
\le 
\Bigl\|  \int_0^\cdot  f(Y_s, \psi_s)ds \Bigr\|_{3\beta}
\le 
K \tau^{1-3\beta} \le K \tau^{\alpha-\beta}.
\label{ineq.917-1}
\end{equation}
We estimate the norm of $\cM^{1}$. 
The $\dagger\dagger$-component of $\cM^{1}$ satisfies 
\begin{align}  
|\nabla \sigma(Y_t) Y^\dagger_t - \nabla \sigma(Y_s) Y^\dagger_s|
&\le
|\nabla \sigma(Y_t) ||Y^{\dagger}_{t} -  Y^{\dagger}_{s}| 
  +  |\nabla \sigma(Y_t)-\nabla \sigma(Y_s) | |Y^{\dagger}_s|
  \nn\\
    &\le 
    K \{ |Y^{\dagger, 1}_{s,t} | +  |Y^{1}_{s,t}  | \|Y^{\dagger}\|_\infty\}
        \nn\\
    &\le 
       9(K+1)^5 (\vertiii{X}_\alpha+1)^3 (t-s)^\alpha, 
       \label{ineq.917-2}
\end{align}
where we have used \eqref{sono13}, \eqref{sono15} and \eqref{sono17}.
The $\sharp\sharp$-component of $\cM^{1}$ reads:
\begin{align}  
\lefteqn{
\sigma(Y_t)-\sigma(Y_s) - \nabla \sigma(Y_s) 
\la Y^\dagger_sX^1_{s,t}\ra
}\nn\\
&=
\{\sigma(Y_t)-\sigma(Y_s) - \nabla \sigma(Y_s) 
\la Y^1_{s,t}\ra\}
+
\nabla\sigma(Y_s)  
\la Y^1_{s,t} - Y^\dagger_sX^1_{s,t} \ra
\nn\\
&=
\int_0^1d\theta (1-\theta)
\nabla^2 \sigma(Y_s + \theta Y^1_{s,t})
  \la    Y^1_{s,t},  Y^1_{s,t}\ra
  +
  \nabla\sigma(Y_s)  \la Y^{\dagger\dagger}_sX^2_{s,t}+Y^{\sharp}_{s,t}\ra.
\label{ineq.917-3}
\end{align}
This implies that
\begin{align}  
|\sigma(Y_t)-\sigma(Y_s) - \nabla \sigma(Y_s) 
\la Y^\dagger_sX^1_{s,t}\ra|
&\le
9 (K+1)^5 (\vertiii{X}_\alpha+1)^4 (t-s)^{2\alpha},
\label{ineq.917-4}
\end{align}
where we have used \eqref{sono17} and \eqref{sono13}.

By \eqref{in.0413-2} in Proposition \ref{prop.gousei},
we can estimate 
the $\sharp$-component of $\cM^{1}$ as follows:
\begin{align}  
\lefteqn{
\Bigl|
\int_{s}^{t} \sigma(Y_{u}) d X_{u} 
 -\sigma(Y_{s})  X_{s, t}^{1}
   -\sigma(Y)_s^{\dagger} X_{s, t}^{2}
\Bigr|
}
\nn\\
&\le
|\sigma(Y)^{\dagger\dagger}_s X_{s, t}^{3}|
+\Bigl|
\int_{s}^{t} \sigma(Y_{u}) d X_{u} 
  -\{ \sigma(Y_{s}) X_{s, t}^{1} 
   +\sigma(Y)_{s}^{\dagger} X_{s, t}^{2}
 +\sigma(Y)_{s}^{\dagger\dagger} X_{s, t}^{3} \}
\Bigr|
\nn\\
&\le
\{ |\nabla \sigma (Y_s) \la Y^{\dagger\dagger}_s \la\bullet, \star \ra,  \,*\ra |
+
|\nabla^2 \sigma (Y_s) 
 \la Y^\dagger_s \bullet, Y^\dagger_s \star, \,*\ra|  \}  
|X_{s, t}^{3}|
\nn\\
&\quad 
+\kappa_{\beta} (t-s)^{4\beta} 
\{ 
\|\sigma (Y)^{\sharp} \|_{3\beta} \|X^1\|_{\beta}
+ 
\|\sigma (Y)^{\sharp\sharp} \|_{2\beta} \|X^2\|_{2\beta}
+ 
\|\sigma (Y)^{\dagger\dagger} \|_{\beta} \|X^3\|_{3\beta}
\}
\nn\\
&\le
c_5 (K+1)^7 (\vertiii{X}_\alpha+1)^7 (t-s)^{3\alpha},
\label{ineq.917-5}
\end{align}
where we have also used \eqref{sono11}, \eqref{sono15},
\eqref{ineq.913-2} and \eqref{ineq.913-3}. 

Combining  \eqref{ineq.917-2}--\eqref{ineq.917-4}, we obtain 
\begin{align} 
\|\cM^{1} (Y, Y^\dagger, Y^{\dagger\dagger}) \|_{\cQ^{\beta}_X}
&\le 
(c_5\vee 9) 
\Bigl\{
(K+1)^4 (\vertiii{X}_\alpha+1)^2 \tau^{\alpha-\beta}
\nn\\
&+
(K+1)^5 (\vertiii{X}_\alpha+1)^4 \tau^{2(\alpha-\beta)}
 +
(K+1)^7 (\vertiii{X}_\alpha+1)^7 \tau^{3(\alpha - \beta)}
\Bigr\}.
\nn
\end{align}
From this and \eqref{ineq.917-1} we can estimate $\cM^{\xi} (Y, Y^\dagger, Y^{\dagger\dagger})$. 
If 
\[
\tau^{\alpha-\beta}
 \le 
\left[ 4 (c_5 \vee 9) (K+1)^4 (\vertiii{X}_\alpha+1)^3 \right]^{-1},
\]
then we have
\[
\|\cM^{\xi} (Y, Y^\dagger, Y^{\dagger\dagger}) \|_{\cQ^{\beta}_X}
\le 
\|\cM^{1} (Y, Y^\dagger, Y^{\dagger\dagger}) \|_{\cQ^{\beta}_X}
+
\|\cM^{2} (Y, Y^\dagger, Y^{\dagger\dagger}) \|_{\cQ^{\beta}_X}
\le 1,
\]
which is equivalent to $\cM^{\xi} (B^\xi) \subset B^\xi$.


Now we will prove that $\cM^\xi$ is a contraction on
$B^\xi$ if we take $\tau$ smaller.
We can easily see from \eqref{sono5} that 
\begin{align}
\lefteqn{
\|\cM^{2} (Y, Y^\dagger, Y^{\dagger\dagger})
-
 \cM^{2} (\tilde{Y}, \tilde{Y}^\dagger, \tilde{Y}^{\dagger\dagger}) \|_{\cQ^{\beta}_X}
}
\nn\\
&\le 
\Bigl\|  \int_0^\cdot \{ f(Y_s, \psi_s) - f(\tilde{Y}_s, \psi_s)\}ds \Bigr\|_{3\beta}
\nn\\
&\le
K \| Y- \tilde{Y}\|_\infty \tau^{1-3\beta} 
\nn\\
&\le
K [
  \{1 + 2(K^2+1) (\vertiii{X}_\alpha^2+1) \tau^\alpha\}
      |\xi -\tilde\xi|
    +
     2 (\vertiii{X}_\alpha^2+1)\Delta 
          \tau^{\alpha}
          ]
\tau^{\alpha-\beta}
\nn\\
&\le
  3(K^2+1) (\vertiii{X}_\alpha^2+1) |\xi -\tilde\xi|
+
 2K (\vertiii{X}_\alpha^2 +1)  \tau^{2\alpha-\beta} \Delta.
 \label{ineq.918-1}
\end{align}

We will estimate the difference of $\cM^1$. 
To do so,  we will first show that there are constants 
$c_6 >0$ and $\nu_1\in\N$ such that
\begin{align}  
\|(\sigma(Y), \sigma(Y)^{\dagger}, \sigma(Y)^{\dagger\dagger}) 
-
(\sigma(\tilde{Y}), \sigma(\tilde{Y})^{\dagger}, \sigma(\tilde{Y})^{\dagger\dagger}) 
\|_{\cQ^{\beta}_X}
\nn\\
\le
c_6  (K+1)^{\nu_1} (\vertiii{X}_\alpha+1)^{\nu_1}
  (|\xi -\tilde{\xi}|+\Delta )
\label{ineq.919-33}
\end{align}
holds for all $(Y, Y^\dagger, Y^{\dagger\dagger})\in B^\xi$ and $(\tilde{Y}, \tilde{Y}^\dagger, \tilde{Y}^{\dagger\dagger}) \in B^{\tilde\xi}$.  

We can estimate
the difference of the $\dagger\dagger$-components as follows:
\begin{align} 
\lefteqn{
|\{\sigma (Y)^{\dagger\dagger}_t -\sigma (Y)^{\dagger\dagger}_s\}
-
\{\sigma (\tilde{Y})^{\dagger\dagger}_t -\sigma (\tilde{Y})^{\dagger\dagger}_s\}
|
}
\nn\\
&\le 
|\{
\nabla \sigma (Y_t) Y^{\dagger\dagger}_t -\nabla \sigma (Y_s) Y^{\dagger\dagger}_s
\}
-
\{
\nabla \sigma (\tilde{Y}_t) \tilde{Y}^{\dagger\dagger}_t 
-\nabla \sigma (\tilde{Y}_s) \tilde{Y}^{\dagger\dagger}_s
\}
|
\nn\\
&\quad +
|\{
\nabla^2 \sigma(Y_t) \la Y^\dagger_t \bullet, Y^\dagger_t \star\ra
-\nabla^2 \sigma(Y_s) \la Y^\dagger_s \bullet, Y^\dagger_s \star\ra
\}
\nn\\
&\qquad\qquad
-
\{
\nabla^2 \sigma(\tilde{Y}_t) \la \tilde{Y}^\dagger_t \bullet, \tilde{Y}^\dagger_t \star\ra
-\nabla^2 \sigma(\tilde{Y}_s) \la \tilde{Y}^\dagger_s \bullet, \tilde{Y}^\dagger_s \star\ra
\}
|
\nn\\
&\le 
\left|
\int_0^1 d\theta \{ \nabla^2 \sigma (Y_s+ \theta Y^1_{s,t}) 
\la  Y^1_{s,t}, Y^{\dagger\dagger}_t\ra
-
 \nabla^2 \sigma (\tilde{Y}_s+ \theta \tilde{Y}^1_{s,t}) 
\la  \tilde{Y}^1_{s,t}, \tilde{Y}^{\dagger\dagger}_t\ra \}
\right|
\nn\\
&\quad + 
| \nabla\sigma (Y_s) \la Y^{\dagger\dagger, 1}_{s, t} \ra 
 - 
 \nabla\sigma (\tilde{Y}_s) \la \tilde{Y}^{\dagger\dagger, 1}_{s, t} \ra|
 \nn\\
   & \quad +
   \left|
   \int_0^1 d\theta \{ \nabla^3 \sigma (Y_s+ \theta Y^1_{s,t}) 
\la  Y^1_{s,t},  Y^\dagger_t \bullet, Y^\dagger_t \star \ra
 -
 \nabla^3 \sigma (\tilde{Y}_s+ \theta \tilde{Y}^1_{s,t}) 
\la  \tilde{Y}^1_{s,t},  \tilde{Y}^\dagger_t \bullet, \tilde{Y}^\dagger_t \star \ra
  \}  \right|
     \nn\\
   & \quad + 
     |   \nabla^2 \sigma(Y_s) 
      \la Y^{\dagger, 1}_{s,t} \bullet, Y^\dagger_t \star\ra
-
\nabla^2 \sigma(\tilde{Y}_s) 
      \la \tilde{Y}^{\dagger, 1}_{s,t} \bullet, \tilde{Y}^\dagger_t \star\ra |
     \nn\\
   & \quad + 
     |   \nabla^2 \sigma(Y_s) 
      \la Y^{\dagger}_{s} \bullet, Y^{\dagger, 1}_{s,t} \star\ra
-
\nabla^2 \sigma(\tilde{Y}_s) 
      \la \tilde{Y}^{\dagger}_{s} \bullet, \tilde{Y}^{\dagger, 1}_{s,t} \star\ra |
\nn\\
&\le
c_7
(K+1)^{\nu_2} (\vertiii{X}_\alpha+1)^{\nu_2}  (|\xi -\tilde{\xi}|+\Delta )
\}
(t-s)^\beta,
\label{ineq.918-2}
\end{align}
where we have used \eqref{sono11}--\eqref{sono5} many times.
%
%


Using \eqref{in.0715-4}, we can estimate
the difference of the $\sharp\sharp$-components
in a similar way as follows:
\begin{align}
\lefteqn{
|\sigma (Y)_{s, t}^{\sharp\sharp}
-\sigma (\tilde{Y})_{s, t}^{\sharp\sharp} |
}
\nn\\
&\le
|\nabla \sigma (Y_s) \la Y_{s, t}^{\sharp\sharp}\ra
-  
  \nabla \sigma (\tilde{Y}_s) \la \tilde{Y}_{s, t}^{\sharp\sharp}\ra |
\nn\\
&\quad
+
|\nabla^2 \sigma (Y_s) \la Y_{s}^{\dagger\dagger} X_{s, t}^{2}+Y_{s, t}^{\sharp}, \, Y_{s}^{\dagger} \ra
 -
  \nabla^2 \sigma (\tilde{Y}_s) \la \tilde{Y}_{s}^{\dagger\dagger} X_{s, t}^{2}+\tilde{Y}_{s, t}^{\sharp}, \, \tilde{Y}_{s}^{\dagger} \ra |
\nn\\
&\quad
+
|\nabla^2 \sigma (Y_s) \la 
Y_{s}^{\dagger} X_{s, t}^{1} +Y_{s}^{\dagger\dagger} X_{s, t}^{2}+Y_{s, t}^{\sharp}, \,
      Y_{s}^{\dagger\dagger} X_{s, t}^{1}  +Y_{s, t}^{\sharp\sharp}\ra
      \nn\\
&\qquad\qquad       -
        \nabla^2 \sigma (\tilde{Y}_s) \la 
        \tilde{Y}_{s}^{\dagger} X_{s, t}^{1}+\tilde{Y}_{s}^{\dagger\dagger} X_{s, t}^{2}+\tilde{Y}_{s, t}^{\sharp}, \,
      \tilde{Y}_{s}^{\dagger\dagger} X_{s, t}^{1}  +\tilde{Y}_{s, t}^{\sharp\sharp} \ra|
\nn\\
&\quad
+
\Bigl|\int_{0}^1 d \theta(1-\theta)
\{
\nabla^{3} \sigma (Y_{s}+\theta Y_{s, t}^{1})
\langle Y_{s, t}^{1}, Y_{s, t}^{1}, Y^\dagger_t \rangle
-
\nabla^{3} \sigma (\tilde{Y}_{s}+\theta \tilde{Y}_{s, t}^{1})
\langle \tilde{Y}_{s, t}^{1}, \tilde{Y}_{s, t}^{1}, \tilde{Y}^\dagger_t \rangle
\}
\Bigr|
\nn\\
&\le
c_8  (K+1)^{\nu_3} (\vertiii{X}_\alpha+1)^{\nu_3}
  (|\xi -\tilde{\xi}|+\Delta )(t-s)^{2\beta},
\label{ineq.919-1}
\end{align}
where we have used \eqref{sono11}--\eqref{sono5}  again.
Note that  
$\nabla^{5} \sigma$ is not involved in the above estimation.

Similarly, we see from \eqref{in.0715-2} that
\begin{align}
\lefteqn{
|\sigma (Y)_{s, t}^{\sharp} - \sigma (\tilde{Y})_{s, t}^{\sharp} |
}
\nn\\
&\le 
|\nabla \sigma (Y_{s}) Y_{s, t}^{\sharp}
  -  \nabla \sigma (\tilde{Y}_{s}) \tilde{Y}_{s, t}^{\sharp}|
\nn\\
&\quad +
|\nabla^2 \sigma (Y_{s})\langle Y^{\dagger}_sX_{s, t}^{1}, 
  Y^{\dagger\dagger}_s X_{s, t}^{2}\rangle
    - 
      \nabla^2 \sigma (\tilde{Y}_{s})\langle \tilde{Y}^{\dagger}_sX_{s, t}^{1}, 
  \tilde{Y}^{\dagger\dagger}_s X_{s, t}^{2}\rangle|
\nn\\
    & \quad +
\frac{1}{2}  
|\nabla^2 \sigma (Y_{s})\langle Y^{\dagger\dagger}_sX_{s, t}^{2},
   Y^{\dagger\dagger}_sX_{s, t}^{2}\rangle
    -
     \nabla^2 \sigma (\tilde{Y}_{s})\langle \tilde{Y}^{\dagger\dagger}_sX_{s, t}^{2},
   \tilde{Y}^{\dagger\dagger}_sX_{s, t}^{2}\rangle |
   \nn\\
    & \quad + 
    |\nabla^2 \sigma (Y_{s})\langle Y^{\dagger}_sX_{s, t}^{1}+
  Y^{\dagger\dagger}_s X_{s, t}^{2}, \, Y^{\sharp}_{s, t}\rangle
   - 
   \nabla^2 \sigma (\tilde{Y}_{s})\langle \tilde{Y}^{\dagger}_sX_{s, t}^{1}+
  \tilde{Y}^{\dagger\dagger}_s X_{s, t}^{2}, \, \tilde{Y}^{\sharp}_{s, t}\rangle |
   \nn\\
    & \quad +
\frac{1}{2}  |\nabla^2 \sigma (Y_{s})\langle Y^{\sharp}_{s, t},
   Y^{\sharp}_{s,t}\rangle
    -  
      \nabla^2 \sigma (\tilde{Y}_{s})\langle \tilde{Y}^{\sharp}_{s, t},
   \tilde{Y}^{\sharp}_{s,t}\rangle |
\nn\\
&\quad
+\frac12
\Bigl|  \int_{0}^1 d \theta(1-\theta)^2 
\{
\nabla^{3} \sigma (Y_{s}+\theta Y_{s, t}^{1})\langle Y_{s, t}^{1}, Y_{s, t}^{1}, Y_{s, t}^{1}\rangle 
-
\nabla^{3} \sigma (\tilde{Y}_{s}+\theta \tilde{Y}_{s, t}^{1})
\langle \tilde{Y}_{s, t}^{1}, \tilde{Y}_{s, t}^{1}, \tilde{Y}_{s, t}^{1}\rangle 
  \}\Bigr|
\nn\\
&\le
c_9  (K+1)^{\nu_4} (\vertiii{X}_\alpha+1)^{\nu_4}
  (|\xi -\tilde{\xi}|+\Delta ) (t-s)^{3\beta},
\label{ineq.919-2}
\end{align}
where we have used \eqref{sono11}--\eqref{sono5} again.
Note that  
$\nabla^{5} \sigma$ is not involved in the above estimation.
Combining \eqref{ineq.918-2}--\eqref{ineq.919-2}, 
we  obtain \eqref{ineq.919-33}.

We estimate the difference of $\cM^{1}$'s. 
Let us first calculate the $\dagger\dagger$-component.
\begin{align}  
\lefteqn{
|\{
\nabla \sigma(Y_t) Y^\dagger_t - \nabla \sigma(Y_s) Y^\dagger_s|
\}
-
\{
\nabla \sigma(\tilde{Y}_t) \tilde{Y}^\dagger_t - \nabla \sigma(\tilde{Y}_s) \tilde{Y}^\dagger_s
\}|
}
\nn\\
&\le
\left|
\int_0^1 d\theta \{
\nabla^2 \sigma (Y_s +\theta Y^{1}_{s,t}) 
\la Y^{1}_{s,t}, Y^{\dagger}_t \ra
\}
- 
\nabla^2 \sigma (\tilde{Y}_s +\theta \tilde{Y}^{1}_{s,t}) 
\la \tilde{Y}^{1}_{s,t}, \tilde{Y}^{\dagger}_t \ra
\right|
\nn\\
&\quad  
  + | \nabla \sigma (Y_s) \la Y^{\dagger, 1}_{s,t}\ra 
    - \nabla \sigma (\tilde{Y}_s) \la \tilde{Y}^{\dagger, 1}_{s,t}\ra |
        \nn\\
    &\le 
      c_{10} (K+1)^{\nu_5} (\vertiii{X}_\alpha+1)^{\nu_5} 
       (|\xi -\tilde{\xi}|+\Delta )(t-s)^\alpha.
       \label{ineq.919-3}
\end{align}
Note that the right hand side of \eqref{ineq.919-3} has $(t-s)^\alpha$, not $(t-s)^\beta$.

Recall the explicit expression of  
the $\sharp\sharp$-component of $\cM^{1}$ in \eqref{ineq.917-3}. 
Using it, we can estimate the difference as follows:
\begin{align}  
\lefteqn{
|  \{
\sigma(Y_t)-\sigma(Y_s) - \nabla \sigma(Y_s) 
\la Y^\dagger_sX^1_{s,t}\ra \}
-
\{
\sigma(\tilde{Y}_t)-\sigma(\tilde{Y}_s) - \nabla \sigma(\tilde{Y}_s) 
\la \tilde{Y}^\dagger_sX^1_{s,t}\ra\}  |
}
\nn\\
&\le
\left| \int_0^1d\theta (1-\theta) 
\{
\nabla^2 \sigma(Y_s + \theta Y^1_{s,t})
  \la    Y^1_{s,t},  Y^1_{s,t}\ra
   - 
  \nabla^2 \sigma(\tilde{Y}_s + \theta \tilde{Y}^1_{s,t})
  \la    \tilde{Y}^1_{s,t},  \tilde{Y}^1_{s,t}\ra
   \} \right|
 \nn\\
 &\qquad +
  | \nabla\sigma(Y_s)  \la Y^{\dagger\dagger}_sX^2_{s,t}+Y^{\sharp}_{s,t}\ra
  - 
  \nabla\sigma(\tilde{Y}_s)  \la \tilde{Y}^{\dagger\dagger}_sX^2_{s,t}+\tilde{Y}^{\sharp}_{s,t}\ra|
   \nn\\
&\le
 c_{11} (K+1)^{\nu_6} (\vertiii{X}_\alpha+1)^{\nu_6} 
  (|\xi -\tilde{\xi}|+\Delta ) (t-s)^{2\alpha}.
\label{ineq.919-4}
\end{align}
Note that the right hand side of \eqref{ineq.919-4} has $(t-s)^{2\alpha}$, not $(t-s)^{2\beta}$.

We now turn to  $\sharp$-component. 
Since the rough path integration is linear when 
the driving RP is fixed,
our calculation is quite similar to that in \eqref{ineq.917-5}.
\begin{align}  
\lefteqn{
\Bigl|
\int_{s}^{t} \{\sigma(Y_{u})- \sigma(\tilde{Y}_{u})\} d X_{u} 
 - \{ \sigma(Y_{s})  -  \sigma(\tilde{Y}_{s})\}X_{s, t}^{1}
   -\{\sigma(Y)_s^{\dagger} -\sigma(\tilde{Y})_s^{\dagger}\}X_{s, t}^{2}
\Bigr|
}
\nn\\
&\le
|\{ \sigma(Y)^{\dagger\dagger}_s - \sigma(\tilde{Y})^{\dagger\dagger}_s\} X_{s, t}^{3}|
\nn\\
&\quad
+\Bigl|
\int_{s}^{t} \{ \sigma(Y_{u}) -  \sigma(\tilde{Y}_{u})\} d X_{u} 
\nn\\
 &\qquad 
  -\left[ \{ \sigma(Y_{s}) - \sigma(\tilde{Y}_{s})\} X_{s, t}^{1} 
   +\{ \sigma(Y)_{s}^{\dagger} - \sigma(\tilde{Y})_{s}^{\dagger}\}X_{s, t}^{2}
 +\{ \sigma(Y)_{s}^{\dagger\dagger} 
   -\sigma(\tilde{Y})_{s}^{\dagger\dagger}\}X_{s, t}^{3} \}
\right] \Bigr|
\nn\\
&\le
|\nabla \sigma (Y_s) \la Y^{\dagger\dagger}_s \la\bullet, \star \ra,  \,*\ra
- 
 \nabla \sigma (\tilde{Y}_s) \la \tilde{Y}^{\dagger\dagger}_s \la\bullet, \star \ra,  \,*\ra|
\cdot |X_{s, t}^{3}|
\nn\\
&\quad +
|\nabla^2 \sigma (Y_s) 
 \la Y^\dagger_s \bullet, Y^\dagger_s \star, \,*\ra
  -  
  \nabla^2 \sigma (\tilde{Y}_s) 
 \la \tilde{Y}^\dagger_s \bullet, \tilde{Y}^\dagger_s \star, \,*\ra  |
   \cdot |X_{s, t}^{3}|
\nn\\
&\quad 
+\kappa_{\beta} (t-s)^{4\beta} 
\{ 
\|Y^{\sharp}- \tilde{Y}^{\sharp} \|_{3\beta} \|X^1\|_{\beta}
+ 
\|Y^{\sharp\sharp}- \tilde{Y}^{\sharp\sharp}\|_{2\beta} \|X^2\|_{2\beta}
+ 
\|Y^{\dagger\dagger}- \tilde{Y}^{\dagger\dagger}\|_{\beta} \|X^3\|_{3\beta}
\}
\nn\\
&\le
c_{12} (K+1)^{\nu_7} (\vertiii{X}_\alpha+1)^{\nu_7}
 (|\xi -\tilde{\xi}|+\Delta ) (t-s)^{3\alpha}.
\label{ineq.919-5}
\end{align}
Note that the right hand side of \eqref{ineq.919-5} has $(t-s)^{3\alpha}$, not $(t-s)^{3\beta}$.
Combining \eqref{ineq.919-3}--\eqref{ineq.919-5}, we easily obtain 
\begin{align} 
\|\cM^{1} (Y, Y^\dagger, Y^{\dagger\dagger})
&- \cM^{1} (\tilde{Y}, \tilde{Y}^\dagger, \tilde{Y}^{\dagger\dagger})
 \|_{\cQ^{\beta}_X}
\nn\\
&\le 
c_{13}(K+1)^{\nu_8} (\vertiii{X}_\alpha+1)^{\nu_8}
 (|\xi -\tilde{\xi}|+\Delta ) \tau^{\alpha -\beta}.
 \label{ineq.919-6}
\end{align}


From \eqref{ineq.919-6} and \eqref{ineq.918-1} we can estimate 
$\cM^{\xi} (Y, Y^\dagger, Y^{\dagger\dagger})$. 
There are constants $c_{14} \ge 4$ and $\nu_9 \ge 2$ such that
if
\begin{equation}\label{def.925-1}
\tau \le \lambda
\qquad \mbox{with} \quad
 \lambda := 
\left[ c_{14} (K+1)^{\nu_9} (\vertiii{X}_\alpha+1)^{\nu_9}
 \right]^{-1/ (\alpha-\beta)}  \, \in (0,1),
\end{equation}
then we have $\cM^{\xi} (B^\xi) \subset B^\xi$ for every $\xi$ and 
\begin{equation}\label{key1}
\|\cM^{\xi} (Y, Y^\dagger, Y^{\dagger\dagger})
-
 \cM^{\xi} (\tilde{Y}, \tilde{Y}^\dagger, \tilde{Y}^{\dagger\dagger})
\|_{\cQ^{\beta}_X}
\le \frac12(|\xi - \tilde{\xi}|+ \Delta).
\end{equation}
In particular, $\cM^{\xi}$ is a contraction on 
$B^\xi = B^\xi_{[0, \tau]}$ for every $\xi$
and therefore has a unique fixed point in this ball. 
Thus, we have obtained
a local solution of RDE \eqref{rde.0413} on $[0,\lambda]$.
Note that $\lambda$ is determined by $\vertiii{X}_{\alpha}$ and $K$, 
but independent of $\xi$ and  $\psi$.

Next,  we do the same thing on the second interval 
$[\lambda, (2\lambda)\wedge 1]$ with the initial condition $\xi$ 
at $t=0$ being replaced by 
$Y_{\lambda}$ at $t=\lambda$.
Since all the estimates above is independent of $\xi$ and $\psi$,
$(Y_s, Y^\dagger_s, Y^{\dagger\dagger}_s)_{s \in [\lambda, (2\lambda)\wedge 1]}$
satisfies the same estimates as those for $(Y_s, Y^\dagger_s, Y^{\dagger\dagger}_s)_{s \in [0,\lambda]}$.
By concatenating them as in the fourth item of Example \ref{ex.0715}, 
we obtain  a solution on $[0, (2\lambda)\wedge 1]$.

We can continue this procedure to obtain a global 
$(Y_s, Y^\dagger_s, Y^{\dagger\dagger}_s)_{s \in [0, 1]}$.
There are (at most) $\lfloor \lambda^{-1}\rfloor +1$ subintervals,
where $\lfloor \cdot \rfloor$ stands for the integer part.
Except (perhaps) the last one, the length of each interval equals 
$\lambda$.
On each subinterval, $(Y, Y^\dagger, Y^{\dagger\dagger})$ satisfies the same estimates.
In particular, Inequality \eqref{sono17} implies that
$\beta$-H\"older norm of $Y$ on each subinterval
 is dominated by $1$.
By H\"older's inequality for finite sums, we can easily see that 
\begin{align*}
\|Y\|_{\beta, [0,1]} \le 
(\lfloor\lambda^{-1}\rfloor +1)^{1-\beta}
\le
c_{15}  \{(K+1) (\vertiii{X}_{\alpha} +1)  \}^{\nu_{10}},
\end{align*}
which is the desired estimate for a global solution.

We now check that any solution 
$(Y, Y^{\dagger}, Y^{\dagger\dagger})
=(Y, \sigma (Y), \nabla \sigma\cdot \sigma (Y))$
of RDE \eqref{rde.0413} defined on $[0, \tau']$ ($0 <\tau' \le 1$)
actually belongs to $\cQ^{\alpha}_X ([0,\tau'], \cW)$. 
By the estimate \eqref{eq.0413-3}, we have
\[
\Bigl|
\int_{s}^{t} \sigma(Y_{u}) d X_{u} 
- (\sigma(Y_{s}) X_{s, t}^{1}+\sigma(Y)^{\dagger}_s X_{s, t}^{2}
 +\sigma(Y)_s^{\dagger\dagger} X_{s, t}^{3})
\Bigr|_\cW
\le 
C_{\beta} (t-s)^{4\beta},
\]
where the constants $C_\beta >0$ depends only on 
$\beta$-H\"older RP norm of $X$ and $\beta$-H\"older seminorm
of $(Y, Y^{\dagger}, Y^{\dagger\dagger})$.
Since $X$ is $\alpha$-H\"older RP and $\sigma(Y)^{\dagger}=
\nabla \sigma\cdot \sigma (Y)$, 
the above inequality implies that $Y^{\sharp}$ is of $3\alpha$-H\"older.
Likewise,
$Y$ is $\alpha$-H\"older continuous and so are $Y^{\dagger\dagger}$.
We can compute $Y^{\sharp\sharp}$ in a similar way as before as follows:
\begin{align*}  
\lefteqn{
Y^{\dagger}_t -Y^{\dagger}_s -Y^{\dagger\dagger}_s X^1_{s,t}
}
\nn\\
&=
\sigma (Y_t)- \sigma (Y_s)- \nabla \sigma(Y_s) \la Y^1_{s,t} \ra
+
\nabla \sigma(Y_s) \la  Y^1_{s,t}-\sigma (Y_s) X^1_{s,t}  \ra
\nn\\
&=
\int_0^1 d\theta \nabla^2 \sigma(Y_s +\theta Y^1_{s,t}) \la Y^1_{s,t},Y^1_{s,t} \ra
+
\nabla \sigma(Y_s) \la  Y^{\dagger\dagger}_sX^2_{s,t} + Y^{\sharp}_{s,t}\ra.
\end{align*}
The right hand side is dominated by a constant multiple of $(t-s)^{2\alpha}$,
which implies that $Y^{\sharp\sharp}$ is of $2\alpha$-H\"older.
Thus, we have seen that any solution is actually $\alpha$-H\"older CP.

Finally, we show the uniqueness of solution.
The uniqueness is a time-local issue, so it suffices to prove 
that any two solutions,
 $(Y, \sigma (Y), \nabla \sigma\cdot \sigma (Y))$ and $(\tilde{Y}, \sigma (\tilde{Y}), \nabla \sigma\cdot \sigma (\tilde{Y}))$,
of RDE \eqref{rde.0413}  must coincide near $t=0$. 
Since they are $\alpha$-H\"older CPs and $\beta <\alpha$,
they both belong to $B^{\xi}_{[0,\tau']}$ for sufficiently small $\tau' >0$.
Recall that we already know that there is only one fixed point 
of $\cM^{\xi}_{[0,\tau']}$ in this ball. 
Hence, these two local solutions 
must be identically equal on $[0,\tau']$. 
 \end{proof}

 \begin{remark} \label{rem.loc.sol}
 {\rm (i)}~By examining the proof of  Proposition \ref{prop.0429},
 one naturally realize the following:
  Just to prove the existence of a unique global solution
 RDE \eqref{rde.0413} for any given $\psi$, $X$ and $\xi$, 
 it suffices to assume that $\sigma$ is of $C^4_{{\rm b}}$ and $f$
 satisfies that
 \[
 \sup_{y \in \cW, t \in [0,T]} |f(y,\psi_t)|_{\cW} + 
\sup_{y,y' \in \cW, y\neq y', t \in [0,T]} 
\frac{|f(y,\psi_t)-f(y',\psi_t) |_{\cW}}{|y-y'|_{\cW}} <\infty.
 \]
 \noindent
 {\rm (ii)}~By a standard cut-off argument, 
 it immediately follows from {\rm (i)} above that if
 $\sigma$ is of $C^4$ and $f$ is locally Lipschitz 
 continuous in the following sense
 \[
\sup_{|y| \vee |y'| \le N, y\neq y', t \in [0,T]} 
\frac{|f(y,\psi_t)-f(y',\psi_t) |_{\cW}}{|y-y'|_{\cW}} <\infty,
\quad
N \in \N,
\]
then RDE \eqref{rde.0413} has a unique
local solution for any given $\psi$, $X$ and $\xi$.
Hence, a unique solution exists up to either the explosion time 
or the time horizon $T$.
  \end{remark}


Together with RDE \eqref{rde.0413}, we also consider the following
RDE on $[0, T]$ with the same $X$, $\sigma$ and $\xi$:
\begin{equation} \label{rde.0506}
\tilde{Y}_{t}
=\xi +\int_{0}^{t} \tilde{f} (\tilde{Y}_{s}, \tilde{\psi}_{s}) ds
+\int_{0}^{t} \sigma (\tilde{Y}_{s}) d X_{s},
\qquad
\tilde{Y}^\dagger_t = \sigma(\tilde{Y}_t),  
\quad 
\tilde{Y}^{\dagger\dagger}_t = \nabla \sigma(\tilde{Y}_t)\tilde{Y}^\dagger_t.
\end{equation}
We assume that 
$ \tilde{f}\colon \cW \times \cS \to \cW$ is also 
continous and satisfies Condition \eqref{cond.0506-1}.
Let $\tilde{\psi}\colon [0,T] \to \cS$ be another continuous path in $\cS$.

\begin{proposition} \label{prop.0506}
Let $\sigma, f, \tilde{f}, \xi$ be as above.
For $X\in  G\Omega_{\alpha} (\cV)$, $\xi \in \cW$ 
and $\psi, \tilde{\psi}$, 
denote by $(Y, \sigma(Y), \nabla \sigma\cdot \sigma (Y))$ and $(\tilde{Y}, \sigma(\tilde{Y}), \nabla \sigma\cdot \sigma (\tilde{Y}))$ 
the unique solutions of RDEs \eqref{rde.0413} and \eqref{rde.0506}
on $[0,T]$, respectively.
For a bounded, globally Lipschitz map
$g \colon \cW \to \cW$,  set
\begin{equation} \label{eq.0506-1}
M_t :=
(Y_t - \tilde{Y}_{t})
-
\int_{0}^{t} \{g(Y_{s})  - g(\tilde{Y}_{s}) \}ds
-
\int_{0}^{t} \{\sigma (Y_{s}) - \sigma (\tilde{Y}_{s}) \}d X_{s},
 \quad t \in [0,T].
 \end{equation}
 Then, $M \in \cC^{1} (\cW)$ and the following estimate holds
  for every $\beta \in (\tfrac14, \alpha)$: there exist positive constants $c$ and $\nu$ such that   %
 \begin{equation} \label{in.0506-2}
 \|Y - \tilde{Y} \|_{\beta} \le c
  \exp \bigl[c (K'+1)^\nu (\vertiii{X}_{\alpha} +1)^\nu \bigr]
     \| M\|_{3\beta}.
 \end{equation}
 Here, we set $K'= \max\{\| \sigma\|_{C_b^4}, \|f\|_{\infty}, L_f, 
 \|\tilde{f}\|_{\infty}, L_{\tilde{f}}, \|g\|_{\infty}, L_g\}$ and
the constants $c$ and $\nu$ are independent of 
 $X,  \xi, \psi, \tilde{\psi},  \sigma, f, \tilde{f}, g, M$. 
\end{proposition}

\begin{proof} 
Without loss of generality we assume $T=1$.
For simplicity we write $(Y, Y^\dagger, Y^{\dagger\dagger})$ and $(\tilde{Y}, \tilde{Y}^\dagger, \tilde{Y}^{\dagger\dagger})$
for $(Y, \sigma(Y), \nabla \sigma\cdot \sigma (Y))$ and $(\tilde{Y}, \sigma(\tilde{Y}), \nabla \sigma\cdot \sigma (\tilde{Y}))$, respectively.
It is easy to see that $M \in \cC^{1} (\cW)$.
Hence, \eqref{eq.0506-1} is in fact an equality in $\cQ^{\beta}_X (\cW)$
with the $\dagger$-parts and $\dagger\dagger$-parts
being clearly equal.

We set 
$\lambda' :=\left[ c_{14} (K'+1)^{\nu_9} (\vertiii{X}_\alpha+1)^{\nu_9}
 \right]^{-1/ (\alpha-\beta)}
$
by just replacing $K$ by $K'$ in \eqref{def.925-1}.
Set $s_j : = j \lambda'$ for $0 \le j \le \lfloor 1/\lambda'\rfloor$
and $s_N :=1$ with $N:= \lfloor 1/\lambda'\rfloor +1$.
Then, on each subinterval $[s_{j-1}, s_j]$,
$(Y, Y^\dagger, Y^{\dagger\dagger}) \in B_{[s_{j-1}, s_j]}^{\xi_j}$,  
 $(\tilde{Y}, \tilde{Y}^\dagger, \tilde{Y}^{\dagger\dagger}) \in B_{[s_{j-1}, s_j]}^{\tilde{\xi}_j}$
 and the 
estimates in the proof of Proposition \ref{prop.0429} are available
(with $K$ being replaced by $K'$). 
Here, we set $\xi_j := Y_{s_j}$ and $\tilde{\xi}_j := \tilde{Y}_{s_j}$.
From \eqref{key1} we have for all $j$ that
\begin{align}
\Bigl\|
\int_{s_{j-1}}^{\cdot} \{g(Y_{s})  - g(\tilde{Y}_{s}) \}ds
-
\int_{s_{j-1}}^{\cdot} \{\sigma (Y_{s}) - \sigma (\tilde{Y}_{s}) \}d X_{s}
\Bigr\|_{\cQ^{\beta}_X,  [s_{j-1}, s_j]}
\nn\\
\le 
\frac12 |\xi_{j-1} - \tilde\xi_{j-1}|+
\frac12
 \| (Y, Y^\dagger, Y^{\dagger\dagger})-(\tilde{Y},  \tilde{Y}^\dagger, \tilde{Y}^{\dagger\dagger}) \|_{\cQ^{\beta}_X, [s_{j-1}, s_j]}.
 \nn
\end{align}
Taking the seminorms of both sides of \eqref{eq.0506-1} on 
each subinterval, 
we can easily see that
\begin{equation}\label{key3}
\| (Y, Y^\dagger, Y^{\dagger\dagger})
 -(\tilde{Y},  \tilde{Y}^\dagger, \tilde{Y}^{\dagger\dagger}) \|_{\cQ^{\beta}_X, [s_{j-1}, s_j]}
\le
2\|M\|_{3\beta} + |\xi_{j-1} - \tilde\xi_{j-1}|.
\end{equation}

Plugging \eqref{key3} into \eqref{sono18} , we obtain for all $j$ that
\begin{align}
|Y_{s, t}^{1}- \tilde{Y}_{s, t}^{1}| 
&\le
     2\bigl\{
      (K'+1)^2 (\vertiii{X}_\alpha^2+1)  |\xi_{j-1}-\tilde{\xi}_{j-1}|
      \nn\\
       &\qquad \quad +
       (\vertiii{X}_\alpha^2+1) 
       (2\|M\|_{3\beta} + |\xi_{j-1} - \tilde\xi_{j-1}|)
      \bigr\}
         (t-s)^{\alpha}
          \nn\\
           &\le   |\xi_{j-1} -\tilde\xi_{j-1}| 
              + \|M\|_{3\beta},
                                        \qquad (s,t) \in \triangle_{[s_{j-1}, s_j]} 
                                         \label{in.0509-2}               
                                         \end{align}
               and in particular
               \begin{equation}\nn
               | \xi_{j} -\tilde\xi_{j} | 
                \le 2 |\xi_{j-1} -\tilde\xi_{j-1}| 
                    + \|M\|_{3\beta}.               
                   \end{equation}
               We have also used that $t-s \le \lambda'$ in \eqref{in.0509-2}.   
                By mathematical induction, we have
               \[
               | \xi_{j} -\tilde\xi_{j} | \le (1+2^1 +\cdots + 2^{j-1}) \|M\|_{3\beta}
                 = (2^{j}-1)\|M\|_{3\beta}, 
                  \quad 1\le j \le N.                          
                   \]
               Then, it follows from \eqref{in.0509-2}  that 
               \begin{align} 
               \|Y^{1}- \tilde{Y}^{1}\|_{\beta, [s_{j-1}, s_j]} 
               &\le
    2\bigl\{
      (K'+1)^2 (\vertiii{X}_\alpha^2+1)  |\xi_{j-1}-\tilde{\xi}_{j-1}|
      \nn\\
       &\qquad \quad +
       (\vertiii{X}_\alpha^2+1) 
       (2\|M\|_{3\beta} + |\xi_{j-1} - \tilde\xi_{j-1}|)
      \bigr\}
         (\lambda')^{\alpha -\beta}
         \nn\\
          & \le  
          2^N 
          \|M\|_{3\beta}, 
                  \quad 1\le j \le N.                                 
               \label{in.0510-1}
\end{align}
By H\"older's inequality for finite sums
and the trivial inequality $N^{1-\beta} 2^N \le 2^{2N}$, 
we see that 
\begin{align} \|Y^{1}- \tilde{Y}^{1}\|_{\beta, [0,1]}
&\le 
          N^{1-\beta} 2^N
           \|M\|_{3\beta}
                    \nn\\
           &\le 
             \exp [ 
             2N (\log 2)  
              ]       \|M\|_{3\beta}
               \nn\\
                 &\le 
                  \exp [ 2( \lfloor 1/\lambda'\rfloor +1) (\log 2)]  \|M\|_{3\beta}.
                        \label{in.0510-2}
\end{align}
By adjusting the positive constants $c$ and $\nu$,
we can easily obtain \eqref{in.0506-2} from \eqref{in.0510-2}.
               \end{proof}

\section{Driving rough path of rough slow-fast system}\label{sec.driving}

In this section we construct the driving RP of our slow-fast system of RDEs.
Unlike the second level case in the previous works \cite{pix2, ina_thk}, 
this part is not so easy. 
We use (a very special case of) the anisotropic RP theory 
developed in \cite{gyur}, in particular, 
an anisotropic version of Lyons' extension theorem.

\subsection{Anisotropic rough paths}

We  write $\cV_1 = \R^d$ and $\cV_2 = \R^e$
and set $\cV := \cV_1\oplus \cV_2 =\R^{d+e}$.
Denote by $\pi_i\colon \cV \to \cV_i~(i=1,2)$ the 
natural projection onto the $i$th component.
Then, $\cV^{\otimes 2} = \oplus_{i, j =1,2} \cV_i\otimes \cV_j$
and 
$\cV^{\otimes 3} = \oplus_{i, j, k =1,2} \cV_i\otimes \cV_j \otimes \cV_k$.
Denote by $\pi_{ij}\colon \cV^{\otimes 2} \to \cV_i\otimes \cV_j
~(i,j=1,2)$ the 
natural projection onto the $(i, j)$th component.
Similarly, $\pi_{ijk}\colon \cV^{\otimes 3} \to \cV_i\otimes \cV_j\otimes \cV_k
~(i,j,k=1,2)$ is defined.

For a continuous map $(\Xi^1, \Xi^2, B^3)\colon \triangle_T 
\to \cV \oplus \cV^{\otimes 2} \oplus \cV_1^{\otimes 3}$, 
we write 
$B^1 := \pi_1 \la \Xi^1 \ra$,  $W^1 := \pi_2 \la \Xi^1\ra$, 
$B^2 := \pi_{11}\la \Xi^2\ra$,  $W^2 := \pi_{22} \la \Xi^2\ra$,  
$I [B, W] := \pi_{12} \la \Xi^2\ra$ and $I [W, B] := \pi_{21} \la \Xi^2\ra$.

Let $\alpha \in (\tfrac14, \tfrac13]$ and $\gamma \in (\tfrac13, \tfrac12]$
with $2\alpha +\gamma >1$.
We say that a continuous map 
\[
(\Xi^1, \Xi^2, B^3)\colon \triangle_T 
\to \cV \oplus \cV^{\otimes 2} \oplus \cV_1^{\otimes 3}
\]
is an anisotropic rough path (ARP) of roughness 
$(\alpha, \gamma)$ if it satisfies the following two conditions:
\begin{enumerate} 
\item[(1)](H\"older regularity) 
\[
\max\{
\|B^1\|_{\alpha}, \|B^2\|_{2\alpha}, \|B^3\|_{3\alpha},
\|W^1\|_{\gamma}, \|W^2\|_{2\gamma}, 
\|I[B, W]\|_{\alpha +\gamma}, \|I[W, B]\|_{\alpha +\gamma} \} <\infty.
\]
\item[(2)](Chen's relation) For all $0\le s \le u \le t \le T$, it holds that 
\begin{align}  
\Xi^1_{s,t} &= \Xi^1_{s,u} + \Xi^1_{u,t}, 
\qquad 
\Xi^2_{s,t} = \Xi^2_{s,u} + \Xi^2_{u,t} + \Xi^1_{s,u} \otimes \Xi^1_{u,t},
\nn\\
B^3_{s,t} &= B^3_{s,u} + B^3_{u,t} +
B^1_{s,u} \otimes B^2_{u,t} + B^2_{s,u} \otimes B^1_{u,t}.
\nn
\end{align}
\end{enumerate}
We denote by $\hat{\Omega}_{\alpha, \gamma} (\cV)$
the set of all ARP of roughness $(\alpha, \gamma)$.
This set naturally becomes a complete metric space equipped
with the H\"older norms.
For $z= (b, w) \in \cC^1_0 (\cV)$, 
we can define its natural lift $\hat{S} (z) \in\hat{\Omega}_{\alpha, \gamma} (\cV)$ by using the 
Riemann-Stieltjes integration in a similar way as in the
usual RP theory (for instance, 
$
\pi_{21} \la \hat{S} (z)_{s,t} \ra:= \int_s^t \int_s^{u_2} dw_{u_1} \otimes db_{u_2}$).
We define $G\hat{\Omega}_{\alpha, \gamma} (\cV)$ to be 
the closure of $\{\hat{S} (z) \mid z =(b, w) \in \cC^1_0 (\cV)\}$, 
which is a complete separable metric space by definition.
An element of $G\hat{\Omega}_{\alpha, \gamma} (\cV)$ 
is called a geometric ARP of roughness $(\alpha, \gamma)$ .

In fact, there is a canonical continuous injection 
$G\hat{\Omega}_{\alpha, \gamma} (\cV) \hookrightarrow 
G\Omega_{\alpha} (\cV)$.
We will explain it in what follows.
Let $(\Xi^1, \Xi^2, B^3) \in \hat{\Omega}_{\alpha, \gamma} (\cV)$.
For  
$(i,j,k) \in \{1,2\}^3$ with $(i,j,k)\neq (1,1,1)$ and $(s,t)\in\triangle_T$,
we define
\begin{equation} \label{def.Lext}
\Xi_{s,t}^{3, [ijk]} 
:=
\lim_{|\cP| \searrow 0}
\sum_{l=1}^N   \pi_{ijk} \bigl\la 
\Xi^1_{s, t_{l-1}}\otimes \Xi^2_{t_{l-1},t_{l}}  
  +  \Xi^2_{s, t_{l-1}}\otimes \Xi^1_{t_{l-1},t_{l}} 
  \bigr\ra,
\end{equation}
where $\cP :=\{s =t_0 <t_1 <\cdots < t_N=t\}$ is a partition of $[s,t]$.
As we will see, the limit on the right hand side exists.
We set $\Xi_{s,t}^3$ by $\pi_{ijk} \la \Xi_{s,t}^3\ra = \Xi_{s,t}^{3, [ijk]}$
if $(i,j,k)\neq (1,1,1)$ and $\pi_{111} \la \Xi_{s,t}^3 \ra= B_{s,t}^3$.
We can easily see that
if $(\Xi^1, \Xi^2, B^3)=\hat{S} (z)$ for some 
$z= (b, w) \in \cC^1_0 (\cV)$, then $(\Xi^1, \Xi^2, \Xi^3)= S_3 (z)$, i.e. the natural third-level 
lift of $z$.

The following two propositions are key results in 
constructing our random RP.
They are actually a special case of an anisotropic version of 
Lyons' extension thereom (see \cite[Theorem 2.6]{gyur}).
Since our situation is quite simple, we will provide direct proofs below.

\begin{proposition} \label{prop.gyur1}
For every $(\Xi^1, \Xi^2, B^3) \in \hat{\Omega}_{\alpha, \gamma} (\cV)$, 
$\Xi_{s,t}^{3, [ijk]}$ in \eqref{def.Lext} is well-defined.
If 
\begin{equation} \label{ass.1002-1}
\max\{
\|B^1\|_{\alpha}, \|B^2\|_{2\alpha}, \|B^3\|_{3\alpha},
\|W^1\|_{\gamma}, \|W^2\|_{2\gamma}, 
\|I[B, W]\|_{\alpha +\gamma}, \|I[W, B]\|_{\alpha +\gamma} \} \le M
\end{equation}
holds for a constant $M>0$,
then there is a constant $C=C_{\alpha, \gamma} >0$ such that
 \begin{align}  \label{est.1002-2}
\|\Xi^{3, [ijk]}\|_{\delta}  \le C M^2 \qquad 
 \mbox{where $\delta:= (6-i-j-k)\alpha+ (i+j+k -3)\gamma$.}
\end{align}
Moreover, $(s,t) \mapsto (\Xi_{s,t}^1, \Xi_{s,t}^2, \Xi_{s,t}^3)$ satisfies Chen's relation.
In particular, $(\Xi^1, \Xi^2, \Xi^3)\in \Omega_{\alpha} (\cV)$.
\end{proposition}

\begin{proof} 
Pick any $(s,t) \in \triangle_T$ with $s<t$.
For $\cP =\{s =t_0 <t_1 <\cdots < t_N=t\}$, we set 
\[
A_{s,t} (\cP) :=\sum_{l=1}^N  \bigl(
\Xi^1_{s, t_{l-1}}\otimes \Xi^2_{t_{l-1},t_{l}}  
  +  \Xi^2_{s, t_{l-1}}\otimes \Xi^1_{t_{l-1},t_{l}} 
  \bigr).
\]
We throw away $t_m~(1 \le m \le N-1)$ from $\cP$.
By Chen's relation, we have
\begin{equation} \label{eq.1002-3}
A_{s,t} (\cP) -A_{s,t} (\cP\setminus \{t_m\})
= 
\Xi^1_{ t_{m-1},  t_{m}}\otimes \Xi^2_{t_{m},t_{m+1}}  
  +\Xi^2_{ t_{m-1},  t_{m}}\otimes \Xi^1_{t_{m},t_{m+1}}.
\end{equation}

We only prove the case $(i,j,k) =(1,2,1)$ for brevity.
(The other cases can be done in the same way.)
For $\mathcal{P}$ given as above,
we can find $m ~(1 \le m \le N-1)$ such that $t_{m+1}-t_{m-1}\le 2(t-s)/(N-1)$. Then, we see that
\begin{align}  \label{ineq.1002-4}
\lefteqn{
|\pi_{121} \la A_{s,t} (\cP)\ra 
- \pi_{121} \la A_{s,t} (\cP\setminus \{t_m\})\ra|
}
\nn\\
&\le 
|B^1_{ t_{m-1},  t_{m}}\otimes I[W, B]_{t_{m},t_{m+1}} | 
  +| I[B, W]_{ t_{m-1},  t_{m}}\otimes B^1_{t_{m},t_{m+1}}|
  \nn\\
&\le 
M^2 (t_{m+1}-t_{m-1})^{2\alpha +\gamma}
\nn\\
&\le 
M^2 \left( \frac{2}{N-1}\right)^{2\alpha +\gamma} (t-s)^{2\alpha +\gamma}.
\end{align}
 Repeating the above estimate, we have
\begin{align}  \label{ineq.1002-5}
|\pi_{121} \la A_{s,t} (\cP)\ra | 
\le 
M^2 2^{2\alpha +\gamma} 
\zeta (2\alpha +\gamma)
(t-s)^{2\alpha +\gamma}.
\end{align}
Here, we have used $2\alpha +\gamma >1$.  
Hence, if $\Xi_{s,t}^{3, [121]} :=\lim_{ |\cP|\searrow 0}\pi_{121} \la A_{s,t} (\cP)\ra$ exists, it has the desired H\"older regularity.

Now we show that $\{\pi_{121} \la A_{s,t} (\cP)\ra\}_{\cP}$ is Cauchy in 
$\cP$ as $|\cP|\searrow 0$.
First, consider the case $\cQ \supset \cP$, that is, 
$\cQ$ is a refinement of $\cP$.
On each subinterval $[t_{l-1}, t_l]$ ($1 \le l \le N$), 
we throw away points of $\cQ \cap [t_{l-1}, t_l]$ one by one 
in the same way as above. 
Then, we see from \eqref{ineq.1002-4} that
\begin{align*}
|\pi_{121} \la A_{s,t} (\cQ)\ra 
- \pi_{121} \la A_{s,t} (\cP)\ra|
&\le 
M^2 2^{2\alpha +\gamma} 
\zeta (2\alpha +\gamma)
\sum_{l=1}^N (t_{l}- t_{l-1})^{2\alpha +\gamma}
\\
&\le
M^2 2^{2\alpha +\gamma} 
\zeta (2\alpha +\gamma) T |\cP|^{2\alpha +\gamma -1}.
\end{align*}
For general $\cP$ and $\cQ$, note that $\cP \cup \cQ$ is a 
refinement of both  $\cP$ and $\cQ$. 
Hence, we have 
\begin{align*}
\lefteqn{
|\pi_{121} \la A_{s,t} (\cQ)\ra - \pi_{121} \la A_{s,t} (\cP)\ra|
}
\nn
\\
&\le
|\pi_{121} \la A_{s,t} (\cQ)\ra - \pi_{121} \la A_{s,t} (\cP \cup \cQ)\ra|
+
|\pi_{121} \la A_{s,t} (\cP \cup \cQ)\ra - \pi_{121} \la A_{s,t} (\cP)\ra|
\nn\\
&\le
2M^2 2^{2\alpha +\gamma} 
\zeta (2\alpha +\gamma) T  (|\cP| \vee |\cQ|)^{2\alpha +\gamma -1}.
\end{align*}
This estimate implies the desired Cauchy condition.

Next, we show Chen's relation. Pick $0\le s< u <t \le T$ arbitrarily.
Let $\cP$ and $\cQ$ be a partition of $[s, u]$ and $[u, t]$, respectively.
Then, $\cP \cup \cQ$ is a partition of $[s, t]$ and
 $|\cP \cup \cQ|= |\cP| \vee |\cQ|$.
By Chen's relation for $(\Xi^1, \Xi^2)$, we have
\[
A_{s,t} (\cP\cup \cQ) =A_{s,u} (\cP) +A_{u,t} (\cQ) + 
\Xi^1_{s,u} \otimes \Xi^2_{u, t} +\Xi^2_{s,u} \otimes \Xi^1_{u, t}.
\]
Applying $\pi_{121}$ to both sides and letting $|\cP \cup \cQ|\searrow 0$,
we have
\[
\Xi_{s,t}^{3, [121]} =\Xi_{s,u}^{3, [121]} + \Xi_{u, t}^{3, [121]}
+
B^1_{s,u} \otimes I[W, B]_{u, t} + I[B, W]_{s,u} \otimes B^1_{u, t}.
\]
From this and the H\"older estimates,  we can also show that 
$\triangle \ni (s,t) \mapsto
\Xi_{s,t}^{3, [121]} \in \cV_1\otimes\cV_2\otimes \cV_1$ is continuous. 
\end{proof}

We write $(\Xi^1, \Xi^2, \Xi^3)={\bf Ext}(\Xi^1, \Xi^2, B^3)$. 
As we have just seen, ${\bf Ext}$ is an injection from 
$\hat{\Omega}_{\alpha, \gamma} (\cV)$ to $\Omega_{\alpha} (\cV)$.
Below we will check that this injection is 
locally Lipschitz continuous, i.e. Lipschitz continous on every bounded 
subset of $\hat{\Omega}_{\alpha, \gamma} (\cV)$.

\begin{proposition} \label{prop.gyur2}
Let  $(\Xi^1, \Xi^2, B^3),  (\tilde{\Xi}^1, \tilde{\Xi}^2,  \tilde{B}^3) \in \hat{\Omega}_{\alpha, \gamma} (\cV)$.
Suppose that \eqref{ass.1002-1} holds for both 
$(\Xi^1, \Xi^2, B^3)$ and $(\tilde{\Xi}^1, \tilde{\Xi}^2,  \tilde{B}^3)$
with a common constant $M>0$ and 
\begin{align*} 
\max\{
\|B^1- \tilde{B}^1\|_{\alpha}, \|B^2- \tilde{B}^2\|_{2\alpha}, \|B^3- \tilde{B}^3\|_{3\alpha},
\|W^1- \tilde{W}^1\|_{\gamma}, 
\|W^2- \tilde{W}^2\|_{2\gamma},
\\ 
\|I[B, W] - I[\tilde{B}, \tilde{W}]\|_{\alpha +\gamma},  \,
\|I[W, B] - I[\tilde{W}, \tilde{B}]\|_{\alpha +\gamma} \} \le \ve
\end{align*}
holds for a constant $\ve >0$,
then there is a constant $C'=C' (M, \alpha, \gamma) >0$ 
(independent of $\ve$) such that
 \begin{align}  \label{est.1002-6}
\|\Xi^{3, [ijk]} -\tilde{\Xi}^{3, [ijk]}\|_{\delta}  \le C' \ve
\qquad 
 \mbox{where $\delta:= (6-i-j-k)\alpha+ (i+j+k -3)\gamma$.}
\end{align}
In particular, ${\bf Ext}\colon  G\hat{\Omega}_{\alpha, \gamma} (\cV)\hookrightarrow G\Omega_{\alpha} (\cV)$ is locally Lipschitz continuous.
\end{proposition}

\begin{proof} 
We use the same notation as in the proof of Proposition \ref{prop.gyur1}.
We see from \eqref{eq.1002-3} that 
\begin{align}
\lefteqn{
\{A_{s,t} (\cP) - \tilde{A}_{s,t} (\cP) \}
- \{ 
A_{s,t} (\cP\setminus \{t_m\}) -\tilde{A}_{s,t} (\cP\setminus \{t_m\})
\}
}
\nn\\
&= 
\{
\Xi^1_{ t_{m-1},  t_{m}}\otimes \Xi^2_{t_{m},t_{m+1}}  
-
\tilde{\Xi}^1_{ t_{m-1},  t_{m}}\otimes \tilde{\Xi}^2_{t_{m},t_{m+1}}  
\}
\nn\\
&\qquad  +\{
  \Xi^2_{ t_{m-1},  t_{m}}\otimes \Xi^1_{t_{m},t_{m+1}}
  - \tilde{\Xi}^2_{ t_{m-1},  t_{m}}\otimes \tilde{\Xi}^1_{t_{m},t_{m+1}}
  \}.
\nn
\end{align}
If $m ~(1 \le m \le N-1)$ is such that $t_{m+1}-t_{m-1}\le 2(t-s)/(N-1)$, then 
\begin{align}
\lefteqn{
| \pi_{121} \la A_{s,t} (\cP)  - \tilde{A}_{s,t} (\cP) \ra 
- 
\pi_{121} \la
A_{s,t} (\cP\setminus \{t_m\}) -\tilde{A}_{s,t} (\cP\setminus \{t_m\})
\ra|
}
\nn\\
&= 
|B^1_{ t_{m-1},  t_{m}}\otimes I[W, B]_{t_{m},t_{m+1}} 
-  \tilde{B}^1_{ t_{m-1},  t_{m}}\otimes I[\tilde{W}, \tilde{B}]_{t_{m},t_{m+1}} | 
\nn\\
&\qquad  +
| I[B, W]_{ t_{m-1},  t_{m}}\otimes B^1_{t_{m},t_{m+1}}
-  I[\tilde{B}, \tilde{W}]_{ t_{m-1},  t_{m}}\otimes \tilde{B}^1_{t_{m},t_{m+1}}|
  \nn\\
&\le 
4M\ve (t_{m+1}-t_{m-1})^{2\alpha +\gamma}
\nn\\
&\le 
4M\ve \left( \frac{2}{N-1}\right)^{2\alpha +\gamma} (t-s)^{2\alpha +\gamma}.
\nn
\end{align}
Using the same argument as in the proof of Proposition \ref{prop.gyur1},
we have
\[
|\pi_{121} \la A_{s,t} (\cP) - \pi_{121} \la \tilde{A}_{s,t} (\cP)\ra | 
\le 
4M\ve 2^{2\alpha +\gamma} 
\zeta (2\alpha +\gamma) (t-s)^{2\alpha +\gamma}.
\]
Letting $|\cP| \searrow 0$, we obtain \eqref{est.1002-6} for 
$(i,j,k) =(1,2,1)$.
(The other cases can be shown in the same way.)
Thus, we have proved the local Lipschitz continuity of ${\bf Ext}\colon\hat{\Omega}_{\alpha, \gamma} (\cV) \to \Omega_{\alpha} (\cV)$. 
Combining this with ${\bf Ext}\circ \hat{S} =S_3$, we see that
${\bf Ext} ( G\hat{\Omega}_{\alpha, \gamma} (\cV)) \subset G\Omega_{\alpha} (\cV)$.
\end{proof}

\subsection{An anisotropic version of 
Kolmogorov's continuity criterion for random rough paths}

In this subsection, we prove a Kolmogorov-type 
continuity criterion for random ARPs of second level. 
A standard version of this criterion in the RP setting
 is found in \cite[Theorem 3.1]{fh}.
We will slightly modify the proof of that theorem.

Let $(Q, R, I)=\{(Q_{s,t}, R_{s,t}, I_{s,t})\}_{(s,t) \in \triangle_T}$ be 
an $\R^3$-valued (two-parameter) stochastic process 
defined on a certain probability space $(\Omega, \cF, {\mathbb P})$
with the following two properties: 
\begin{itemize} 
\item
Chen's relation in a probabilistic sense, that is, 
for every $0 \le s \le u \le t \le T$, 
\begin{align}  \label{cond.1006-1}
{\mathbb P} \left(
Q_{s,t} = Q_{s,u}+Q_{u,t}, \, R_{s,t} = R_{s,u}+R_{u, t}, \,
I_{s,t} = I_{s,u}+I_{u,t} +Q_{s,u}R_{u,t}  
\right) =1.
\end{align}
\item
An $L^q$-version, $q \in [2,\infty)$, of the
 H\"older condition with exponent $\beta_1, \beta_2 \in (0,1]$, that is, 
there exists $C>0$ such that for every $(s, t) \in \triangle_T$
\begin{align}  \label{cond.1006-2}
\|Q_{s,t}\|_{L^q} \le C (t-s)^{\beta_1}, \,\,
 \|R_{s,t}\|_{L^q}\le C (t-s)^{\beta_2}, \,\,
  \|I_{s,t}\|_{L^{q/2}} \le C (t-s)^{\beta_1 +\beta_2}.
\end{align}
\end{itemize}


\begin{proposition} \label{prop.1006-1}
Let $(Q, R, I)$ be as above
and let $q \in [2,\infty)$ and $\beta_1, \beta_2 \in (1/q, 1]$.
Assume  that \eqref{cond.1006-1} and \eqref{cond.1006-2} hold.
Then, there exists a continuous modification of $(Q, R, I)$
 (denoted by the same symbol again) with the following two properties: 
 \begin{itemize} 
\item 
Chen's relation holds almost surely, that is, 
\begin{align}  
{\mathbb P} (
Q_{s,t} = Q_{s,u}+Q_{u,t}, \,\, R_{s,t} = R_{s,u}+R_{u, t}, \,\,
I_{s,t} = I_{s,u}+I_{u,t} +Q_{s,u}R_{u,t}
\nn\\
\qquad \mbox{{\rm for all} $0 \le s \le u \le t \le T$} )=1
\label{cond.1006-3}
\end{align}
\item
For every $\alpha_1 \in (0, \beta_1 -1/q)$ and $\alpha_2 \in (0, \beta_2 -1/q)$,
there exist non-negative random variables 
$M^Q, M^R\in L^q$ and $M^I\in L^{q/2}$ such that
\begin{align}  
{\mathbb P} (
|Q_{s,t}| \le M^Q (t-s)^{\alpha_1}, \, 
&|R_{s,t}| \le M^R (t-s)^{\alpha_2}, 
\nn\\
|I_{s,t}| &\le M^I (t-s)^{\alpha_1 +\alpha_2}
\quad \mbox{{\rm for all} $(s, t) \in \triangle_T$} )=1.
\label{cond.1006-4}
\end{align}
\end{itemize}
\end{proposition}

\begin{proof} 
Without loss of generality, we may assume $T=1$.
Set $D_n =\{k2^{-n}\mid 0\le k\le 2^n\}$ for $n\in \N_0$ and 
set $D = \cup_{n =0}^\infty D_n$.
By \eqref{cond.1006-1} there exists $\Omega_1 \in \cF$
 with ${\mathbb P} (\Omega_1)=1$ on which 
\[
Q_{s,t} = Q_{s,u}+Q_{u,t}, \,\, R_{s,t} = R_{s,u}+R_{u, t}, \,\,
I_{s,t} = I_{s,u}+I_{u,t} +Q_{s,u}R_{u,t}
\]
holds for all $(s, u, t) \in D^3$ with $s\le u \le t$.
In particular, it holds on $\Omega_1$ that 
$Q_{s,s}=R_{s,s}=I_{s,s}=0$ for all $s\in D$.
From here we will work on $\Omega_1$.

We set  
\begin{equation}\label{def.1007-3}
K^Q_n =\sup_{1 \le i \le 2^n} |Q_{ (i-1)2^{-n}, i2^{-n}}|, \qquad n\in \N_0.
\end{equation}
We also set $K^R_n$ and $K^I_n$ in the same way.
Next, define $M^Q$, $M^R$ and $M^I$ by 
\[
M^Q:= 2\sum_{n=0}^\infty 2^{n \alpha_1} K^Q_n,
\quad
M^R:= 2\sum_{n=0}^\infty 2^{n \alpha_2} K^R_n,
\quad
M^I:= M^Q M^R + 2\sum_{n=0}^\infty 2^{n (\alpha_1 +\alpha_2 )} K^I_n.
\]

We can easily show 
\[
\E [|K^Q_n|^q] \le 
\sum_{i=1}^{2^n}E [|Q_{ (i-1)2^{-n}, i2^{-n}}|^q ]
\le 
2^n C^q (2^{-n})^{\beta_1 q} \le C^q (2^{-n})^{\beta_1 q -1}
 \]
and, similarly, $\E [|K^R_n|^q] \le C^q (2^{-n})^{\beta_2 q -1}$
and $\E [|K^I_n|^{q/2}] \le C^{q/2} (2^{-n})^{(\beta_1 +\beta_2) q/2 -1}$.

Then, we have
\[
\|M^Q\|_{L^q} \le 2\sum_{n=0}^\infty 2^{n\alpha_1} \|K^Q_n\|_{L^q} 
\le 
2C \sum_{n=0}^\infty (2^{-n})^{-\alpha_1+\beta_1 -1/q } <\infty,
\]
where we have used $\beta_1 -1/q >\alpha_1$.
In the same way, we have $\|M^R\|_{L^q}<\infty$.
From Schwarz' inequality, 
$\|M^Q M^R\|_{L^{q/2}} 
\le\|M^Q\|_{L^q}^{1/2}\|M^R\|_{L^q}^{1/2} <\infty$ immediately follows.
By a similar argument as above we also have
\[
\|M^I -M^QM^R\|_{L^{q/2}} \le 
2\sum_{n=0}^\infty 2^{n (\alpha_1 +\alpha_2)} \|K^I_n\|_{L^{q/2}} 
\le 
2C \sum_{n=0}^\infty (2^{-n})^{-\alpha_1-\alpha_2+\beta_1 +\beta_2-2/q }<\infty.
\]
Thus, we have shown $\|M^I\|_{L^{q/2}}<\infty$.

Take any $(s,t) \in D^2$ with $s<t$.
Then, there uniquely exists $m\in \N_0$ such that $2^{- (m+1)}< t-s \le 2^{-m}$.
Moreover, there exists a partition 
$\{s=\tau_0 <\tau_1 \cdots <\tau_N =t\}$ with $\tau_i \in D~(0\le i \le N)$ 
such that (1)~$(\tau_{i-1}, \tau_{i}) \in D_n\times D_n$ for some $n \ge m+1$
$(1\le i \le N)$ and (2)~at most two of $[\tau_{i-1}, \tau_{i}]$'s have
length $2^{-n}$ for every $n \ge m+1$.

Then we have
\[
|Q_{s,t}| \le \max_{1 \le i \le N} |Q_{s, \tau_i}|
      \le \sum_{i=1}^{N} |Q_{\tau_{i-1}, \tau_i}| 
                  \le 2 \sum_{n =m+1}^\infty K^Q_n
\]
and 
\[
\frac{|Q_{s,t}|}{(t-s)^{\alpha_1}} 
  \le 2^{\alpha_1 (m+1)} 2\sum_{n =m+1}^\infty  K^Q_n
    \le 2\sum_{n =m+1}^\infty 2^{\alpha_1n} K^Q_n \le M^Q.
\]
In the same way, we have $|R_{s,t}|/(t-s)^{\alpha_2} \le M^R$, too.

Next we estimate $I_{s,t}$. By Chen's relation, we have 
\begin{equation}\label{eq.1007-6}
I_{s,t} =\sum_{i=1}^{N} ( I_{\tau_{i-1}, \tau_i} + Q_{s, \tau_{i-1}}R_{\tau_{i-1}, \tau_i})
\end{equation}
and therefore 
\begin{align}  
|I_{s,t}| 
&\le
\sum_{i=1}^{N}  |I_{\tau_{i-1}, \tau_i}| 
+ \max_{1 \le i \le N} |Q_{s, \tau_i}|\sum_{i=1}^{N}  |R_{\tau_{i-1}, \tau_i}| 
\nn\\
&\le 
2 \sum_{n =m+1}^\infty K^I_n 
+\Bigl(2 \sum_{n =m+1}^\infty K^Q_n\Bigr)
  \Bigl(2 \sum_{n =m+1}^\infty K^R_n\Bigr).
  \nn
\end{align}
Repeating a similar argument as above, we can easily see that
\[
\frac{|I_{s,t}|}{(t-s)^{\alpha_1 +\alpha_2}} 
  \le
  2\sum_{n =m+1}^\infty 2^{(\alpha_1+\alpha_2)n} K^I_n 
  +
  \Bigl(2 \sum_{n =m+1}^\infty  2^{\alpha_1n}K^Q_n\Bigr)
  \Bigl(2 \sum_{n =m+1}^\infty 2^{\alpha_2n}K^R_n\Bigr) \le M^I.
\]

From here, we construct a continuous modification.
Set 
\[
\Omega_2 =\{\omega \in \Omega_1 \mid M^Q (\omega) 
\vee M^R (\omega)\vee M^I (\omega)<\infty\}.
\]
Clearly, $\mathbb{P} (\Omega_2)=1$.   
For $\omega \in \Omega_2$, $(Q (\omega), R(\omega), I(\omega))$ satisfies Chen's relation 
on $\triangle_1 \cap D^2$ and the following estimates:
\begin{align}  \label{ineq.1006-9}
|Q_{s,t}(\omega)| \le M^Q (\omega) (t-s)^{\alpha_1}, \quad 
|R_{s,t}(\omega)| \le M^R (\omega)(t-s)^{\alpha_2}, 
\nn\\
|I_{s,t}(\omega)| \le M^I(\omega) (t-s)^{\alpha_1 +\alpha_2}
\quad
\mbox{for every $(s,t) \in \triangle_1 \cap D^2$.}
\end{align}
This implies that the map
\[
\triangle_1 \cap D^2 \ni \,\, (s,t) \mapsto 
(Q_{s,t} (\omega), R_{s,t}(\omega), I_{s,t}(\omega)) \,\, \in \R^3
\]
is uniformly continuous.
Since $\triangle_1 \cap D^2$ is dense in $\triangle_1$, 
 this map extends to a uniformly continuous map from 
$\triangle_1$ to $\R^3$, which is denoted by $(\tilde{Q} (\omega), \tilde{R}(\omega), \tilde{I}(\omega))$.
More concretely, we take a sequence $\{(s_n, t_n)\}_{n\in\N}\subset \triangle_1 \cap D^2$ which converges to $(s,t)\in \triangle_1$
and define 
\[
(\tilde{Q}_{s,t} (\omega), \tilde{R}_{s,t}(\omega), \tilde{I}_{s,t}(\omega))
=\lim_{n\to\infty}
(Q_{s_n, t_n} (\omega), R_{s_n, t_n}(\omega), I_{s_n, t_n}(\omega)).
\]
Of course, the limit does not depend on the choice of $\{(s_n, t_n)\}$.
Obviously, this extended map satisfies Chen's relation 
and the inequalities in \eqref{ineq.1006-9} for every $(s,t)\in\triangle_1$.

Finally, we prove that $(\tilde{Q}, \tilde{R}, \tilde{I})$ is a modification
of $(Q, R, I)$. 
Pick any $(s,t)\in \triangle$ and take a sequence $\{(s_n, t_n)\}_{n\in\N}\subset \triangle_1 \cap D^2$ which converges to $(s,t)\in \triangle_1$.
We have just seen that $(Q_{s_n, t_n}, R_{s_n, t_n}, I_{s_n, t_n})$
converges to $(\tilde{Q}_{s,t}, \tilde{R}_{s,t}, \tilde{I}_{s,t})$ a.s. as $n\to\infty$.  
On the other hand, it follows from
 \eqref{cond.1006-1} and \eqref{cond.1006-2} 
that $(Q_{s_n, t_n}, R_{s_n, t_n}, I_{s_n, t_n})$
converges to $(Q_{s,t}, R_{s,t}, I_{s,t})$ in $L^1$ as $n\to\infty$. 
Thus, we have shown that
$(\tilde{Q}_{s,t}, \tilde{R}_{s,t}, \tilde{I}_{s,t})=(Q_{s,t}, R_{s,t}, I_{s,t})$ holds a.s.
for every fixed $(s,t)\in \triangle_1$.
\end{proof}


Next we estimate the difference of two
$\R^3$-valued (two-parameter) stochastic processes
$(Q, R, I)=\{(Q_{s,t}, R_{s,t}, I_{s,t})\}_{(s,t) \in \triangle_T}$ 
and 
$(\hat{Q}, \hat{R}, \hat{I})
=\{(\hat{Q}_{s,t}, \hat{R}_{s,t}, \hat{I}_{s,t})\}_{(s,t) \in \triangle_T}$.

\begin{proposition} \label{prop.1007-2}
Let $q \in [2,\infty)$ and $\beta_1, \beta_2 \in (1/q, 1]$.
Assume that both $(Q, R, I)$ and $(\hat{Q}, \hat{R}, \hat{I})$ 
satisfy \eqref{cond.1006-1}--\eqref{cond.1006-2} with a common
constant $C>0$.
(Their continuous modifications as in Proposition \ref{prop.1006-1}
are denoted by the same symbols again.)
Assume further that 
\begin{align}  \label{cond.1007-1}
\|Q_{s,t} -\hat{Q}_{s,t} \|_{L^q} \le \ve (t-s)^{\beta_1}, \quad
 \|R_{s,t}-\hat{R}_{s,t} \|_{L^q}\le \ve (t-s)^{\beta_2}, 
\nn\\
  \|I_{s,t}-\hat{I}_{s,t} \|_{L^{q/2}} \le \ve (t-s)^{\beta_1 +\beta_2}
  \qquad \mbox{{\rm for all} $(s, t) \in \triangle_T$}.
\end{align}
holds for a constant $\ve >0$.
Then, for every $\alpha_1 \in (0, \beta_1 -1/q)$ and $\alpha_2 \in (0, \beta_2 -1/q)$, there exists a constant $C' >0$ such that
\begin{align}  \label{cond.1007-2}
\bigl\|
\|Q -\hat{Q}\|_{\alpha_1} \bigr\|_{L^q}
+
\bigl\| \|R -\hat{R}\|_{\alpha_2} \bigr\|_{L^q}
+
\bigl\|  \|I -\hat{I}\|_{\alpha_1 + \alpha_2} \bigr\|_{L^{q/2}}
\le C'\ve.
\end{align}
Here, $C'$ depends on $q, \beta_1, \beta_2, \alpha_1, \alpha_2, C, T$, but not on $\ve$.
\end{proposition}

\begin{proof} 
We assume $T=1$ again. We argue in the same way
and use the same symbols as in the proof of Proposition \ref{prop.1006-1}.

By replacing $Q$ by $Q -\hat{Q}$ in \eqref{def.1007-3},
we set $K^{Q-\hat{Q}}_n$. 
We set $K^{R-\hat{R}}_n$ and $K^{I-\hat{I}}_n$ in a similar way.
Then, we have
\begin{align}  
\E [|K^{Q-\hat{Q}}_n|^q] 
\le 
\sum_{i=1}^{2^n}E [|Q_{ (i-1)2^{-n}, i2^{-n}} - \hat{Q}_{ (i-1)2^{-n}, i2^{-n}}|^q ]
\le 
\ve^q (2^{-n})^{\beta_1 q -1}.
\label{ineq.1007-4}
\end{align}
We also have
 $\E [|K^{R-\hat{R}}_n|^q]\le \ve^q (2^{-n})^{\beta_2 q -1}$
and 
$\E [|K^{I-\hat{I}}_n|^ {q/2}]\le \ve^{q/2}  (2^{-n})^{(\beta_1 +\beta_2) q/2 -1}$.

Take any $(s,t) \in D^2$ with $s<t$
and let $\{s=\tau_0 <\tau_1 \cdots <\tau_N =t\}$ be the partition 
as in the proof of Proposition \ref{prop.1006-1}. 
From  \eqref{ineq.1007-4} we can easily see that
\[
|Q_{s,t}- \hat{Q}_{s,t}| \le \max_{1 \le i \le N} |Q_{s, \tau_i}- \hat{Q}_{s, \tau_i}|
      \le \sum_{i=1}^{N} |Q_{\tau_{i-1}, \tau_i} - \hat{Q}_{\tau_{i-1}, \tau_i} | 
                  \le 2 \sum_{n =m+1}^\infty K^{Q-\hat{Q}}_n
\]
and 
\[
\frac{|Q_{s,t} - \hat{Q}_{s,t}|}{(t-s)^{\alpha_1}} 
  \le 2^{\alpha_1 (m+1)} 2\sum_{n =m+1}^\infty  K^{Q-\hat{Q}}_n
    \le 2\sum_{n =m+1}^\infty 2^{\alpha_1n} K^{Q-\hat{Q}}_n 
\]
and therefore
\[
\bigl\|\|Q -\hat{Q}\|_{\alpha_1} \bigr\|_{L^q}
\le
2\sum_{n =1}^\infty 2^{\alpha_1n} \| K^{Q-\hat{Q}}_n \|_{L^q}
\le
2 \ve \sum_{n=m+1}^\infty (2^{-n})^{-\alpha_1+\beta_1 -1/q }= C'_1 \ve, 
\]
where we set $C'_1 := 2 \sum_{n=0}^\infty (2^{-n})^{-\alpha_1+\beta_1 -1/q } <\infty$.
In essentially the same way, we can also show
$\bigl\|\|R -\hat{R}\|_{\alpha_2} \bigr\|_{L^q} \le C'_2 \ve$
for some constant $C'_2>0$.


Next we turn to the $I$-component. 
Since both $I$ and $\hat{I}$ satisfy \eqref{eq.1007-6},
we see that 
\begin{align}  
|I_{s,t} - \hat{I}_{s,t}|
&\le 
\sum_{i=1}^{N} | I_{\tau_{i-1}, \tau_i} - \hat{I}_{\tau_{i-1}, \tau_i} |
\nn\\
&\quad 
+ \max_{1\le i \le N} | Q_{s, \tau_{i}} - \hat{Q}_{s, \tau_{i}}|
\sum_{i=1}^{N}  |R_{\tau_{i-1}, \tau_i} |
\nn\\
&\quad 
+ \max_{1\le i \le N} |\hat{Q}_{s, \tau_{i}}|
\sum_{i=1}^{N}  |R_{\tau_{i-1}, \tau_i} -\hat{R}_{\tau_{i-1}, \tau_i}|.
\nn\\
&\le
2 \sum_{n =m+1}^\infty K^{I-\hat{I}}_n 
+\Bigl(2 \sum_{n =m+1}^\infty K^{Q-\hat{Q}}_n\Bigr)
  \Bigl(2 \sum_{n =m+1}^\infty K^R_n\Bigr)
\nn\\
&\qquad  +\Bigl(2 \sum_{n =m+1}^\infty K^Q_n\Bigr)
  \Bigl(2 \sum_{n =m+1}^\infty K^{R-\hat{R}}_n\Bigr).
  \nn
\end{align}
By using a similar argument as above, we can easily
see from the above 
estimate that $\bigl\|\|I -\hat{I}\|_{\alpha_1 +\alpha_2} \bigr\|_{L^{q/2}} \le C'_3 \ve$
holds for some constant $C'_3>0$. This completes the proof.
\end{proof}

\subsection{Lemmas for the geometric property 
of random anisotropic rough paths}

The Kolmogorov-type continuity criterion alone is not 
sufficient for proving the geometric property of random ARPs
which drives our slow-fast system.
In this subsection we will provide a few lemmas for that purpose.

Let $w =(w_t)_{t \in [0,T]}$ be a one-dimensional Brownian motion
defined on a probability space $(\Omega, \cF, \mathbb{P})$.
For $m\in\N$, $w(m)=(w(m)_t)_{t \in [0,T]}$ be the $m$th
dyadic piecewise linear approximation of $w$,
that is the piecewise linear approximation associated with 
the partition $\{ kT2^{-m}\mid 0\le k\le 2^m\}$ of $[0,T]$.
Let $g\in \cC^\alpha_0 (\R)$ with $\alpha \in (0,1)$. 
(Note that $g$ is not random).
We set, for $(s,t) \in \triangle_T$,
\begin{align}  
I[g, w]_{s,t} &=\int_s^t (g_u -g_s) d^{{\rm I}}w_u,
\label{def.1009-1}
\\
I[g, w(m)]_{s,t} &=\int_s^t (g_u -g_s) dw(m)_u.
\label{def.1009-2}
\end{align}
Here, the integral in \eqref{def.1009-1} is an It\^o integral,
while that in \eqref{def.1009-2} is a (random) Riemann-Stieltjes integral.

\begin{lemma} \label{lem.1009-1}
For every $q \in [1,\infty)$ and $\alpha \in (0,1)$, 
there exists a constant $C=C_{q, \alpha}>0$  such that 
\[
\bigl\| I[g, w]_{s,t}\bigr\|_{L^q} \vee 
\sup_{m\in\N}\bigl\| I[g, w(m)]_{s,t}\bigr\|_{L^q} 
\le 
C (t-s)^{\alpha + \tfrac12} \|g\|_\alpha,
\quad (s,t) \in \triangle_T, \,\, g\in \cC^\alpha_0 (\R).
\]
Here, $C$ does not depend on $(s,t)$ or $g$.
\end{lemma}

\begin{proof} 
In this proof, $C_j~(0 \le j \le 4)$ are certain positive constants
which depend on $q$ and $\alpha$ only.
Without loss of generality we may assume $T=1$.

By Burkholder's inequality we have 
\[
\E\left[ |I[g, w]_{s,t}|^q\right] 
\le
C_0 \left( \int_s^t |g_u -g_s|^2du\right)^{\frac{q}{2}}
 \le 
 C_0 \|g\|_\alpha^q (2\alpha +1)^{-q/2} (t-s)^{(2\alpha +1) q/2},
\]
which is the desired estimate for $I[g, w]$.

We calculate $I[g, w (m)]$ for given $m\in\N$.
We write $t^m_k := k2^{-m}~(0\le k\le 2^m)$ and 
$\Delta^m_k w := w_{t^m_k} - w_{t^m_{k-1}}$ for simplicity.
Note that the law of $\Delta^m_k w/ 2^{m/2}$ is the standard normal 
distribution.
We will write $g^1_{s,t} = g_t-g_s$ as usual.
Set 
\[
[g]^m_k := 2^m \int^{t^m_k}_{t^m_{k-1}} g_u du, 
\qquad m\in\N, \, 1 \le k \le 2^m. 
\]

First, consider the case that $s$ and $t$ belong to the same subinterval,
that is, there exists $k$ such that $s, t\in [t^m_{k-1},t^m_{k}]$.
Then, we can easily see that
\begin{align}  
|I[g, w (m)]_{s,t}|
&=
\left|
\int_s^t g^1_{s,u} \frac{\Delta^m_k w}{2^{-m}} du
\right|
\nn\\
&\le 
\|g\|_\alpha \frac{(t-s)^{\alpha +1}}{\alpha +1} 
 \left|\frac{\Delta^m_k w}{2^{-m}}\right|
\le 
\|g\|_\alpha   \frac{(t-s)^{\alpha + \tfrac12 }  }{\alpha +1}
\left|\frac{\Delta^m_k w}{2^{-m/2}}\right|
\nn
\end{align}
and therefore 
\begin{equation} \label{ineq.1009-3}
\bigl\| I[g, w(m)]_{s,t}\bigr\|_{L^q}
\le 
C_1 \|g\|_\alpha (t-s)^{\alpha + \tfrac12 }
\end{equation}
for some constant $C_1 >0$.

Next, consider the case $s \in [t^m_{k-1},t^m_{k}]$
and $t \in [t^m_{l},t^m_{l+1}]$ for some $k, l$ with $k\le l$.
By Chen's relation we have
\begin{align} 
I[g, w(m)]_{s,t}
&= 
 I[g, w(m)]_{s,t^m_{k}}
 +
  I[g, w(m)]_{t^m_{k}, t^m_{l}}
  +
   I[g, w(m)]_{t^m_{l},t}
\nn\\ 
&\qquad +g^1_{s,t^m_{k}} w(m)^1_{t^m_{k}, t^m_{l}} 
+g^1_{s,t^m_{k}} w(m)^1_{t^m_{l}, t} 
+g^1_{t^m_{k},  t^m_{l}}w(m)^1_{t^m_{l}, t} 
\nn\\ 
&=: A_1 +\cdots +A_6.
\label{eq.1009-4}
\end{align}
Note that $A_1$ and $A_3$ were already estimated in \eqref{ineq.1009-3}.
Since $w(m)^1_{t^m_{k}, t^m_{l}}=w^1_{t^m_{k}, t^m_{l}}$,
we see that $\|A_4\|_{L^q} \le C_2\|g\|_\alpha (t-s)^{\alpha +1/2}$.
Since 
$
w(m)^1_{t^m_{l}, t} =[\Delta^m_{l+1} w/ 2^{m/2}] \cdot [(t-t^m_{l})/ 2^{m/2}]
$,
the variance of $w(m)^1_{t^m_{l}, t}$  is dominated by $t-t^m_{l}$.
Then, in the same way as above, we see that
$\|A_5\|_{L^q}+\|A_6\|_{L^q} \le C_3\|g\|_\alpha (t-s)^{\alpha +1/2}$.

Finally, we estimate $A_2$. We may assume $k<l$, here.
We can easily see that
\[
A_2 = \sum_{j=k+1}^l [g_{\cdot} - g_{t^m_{k}}]^m_j (\Delta^m_j w)
=
\int_{t^m_{k}}^{t^m_{l}}  \sum_{j=k+1}^l
[g_{\cdot} - g_{t^m_{k}}]^m_j  \mathbf{1}_{[t^m_{j-1},t^m_{j} ]} (u)\,dw_u.
\]
Noting that $|[g_{\cdot} - g_{t^m_{k}}]^m_j| \le
 \|g\|_\alpha (t^m_{l}-t^m_{k})^{\alpha}$, we see from Burkholder's 
 inequality that
 \begin{align*}
 \E[ |A_2|^q]^{\frac{1}{q}}
 \le
C_4 \left(
  \int_{t^m_{k}}^{t^m_{l}} 
  \left| \sum_{j=k+1}^l
[g_{\cdot} - g_{t^m_{k}}]^m_j  \mathbf{1}_{[t^m_{j-1},t^m_{j} ]} (u) 
\right|^2 du
 \right)^{\frac12} 
 \le 
 C_4 \|g\|_{\alpha} (t^m_{l}-t^m_{k})^{\alpha +1/2},
 \end{align*}
which completes the proof of the lemma.
\end{proof}

\begin{lemma} \label{lem.1009-2}
For every $q \in [1,\infty)$, $\alpha \in (0,1)$ and $\kappa \in (0,\alpha)$, 
there exists a 
sequence of positive numbers $\{\ve_m\}_{m\in\N}$ converging to $0$
as $m\to\infty$ such that 
\[
\bigl\| I[g, w(m) ]_{s,t}- I[g, w]_{s,t}\bigr\|_{L^q} 
\le 
 \ve_m  \|g\|_\alpha (t-s)^{\alpha + \tfrac12 -\kappa}
\]
for all $m\in\N$, $(s,t) \in \triangle_T$ and $g\in \cC^\alpha_0 (\R)$.
Here, $\{\ve_m\}$ does not depend on $(s,t)$ or $g$.
\end{lemma}

\begin{proof} 
We assume $T=1$ again and use the same notation 
as in the proof of Lemma \ref{lem.1009-1}.
We will estimate $I[g, w (m)]-I[g, w]$ for given $m\in\N$.
Below, $C_i~(1 \le i \le 4)$ are certain positive constants
which depends only on $q, \alpha, \kappa$.

First, consider the case that $s$ and $t$ belong to the same subinterval,
that is, there exists $k$ such that $s, t\in [t^m_{k-1},t^m_{k}]$.
Then, we see from \eqref{ineq.1009-3} that
\begin{align}  
\bigl\| I[g, w(m)]_{s,t} -I[g, w]_{s,t} \bigr\|_{L^q}
&\le
 \bigl\| I[g, w(m)]_{s,t} \bigr\|_{L^q}+  \bigl\| I[g, w]_{s,t} \bigr\|_{L^q}
\nn\\
&\le
C_1 \|g\|_\alpha (t-s)^{\alpha + \tfrac12}
\nn\\
&\le
C_1 (2^{-m})^{\kappa} \|g\|_\alpha (t-s)^{\alpha + \tfrac12 -\kappa}.
\label{ineq.1011-1}
\end{align}

Next, consider the case $s \in [t^m_{k-1},t^m_{k}]$
and $t \in [t^m_{l},t^m_{l+1}]$ with $k\le l$.
By Chen's relation we have
\begin{align} 
I[g, w(m)]_{s,t}-I[g, w]_{s,t}
&= 
\{ I[g, w(m)]_{s,t^m_{k}} - I[g, w]_{s,t^m_{k}} \}
 \nn\\
 &\quad 
 +
  \{ I[g, w(m)]_{t^m_{k}, t^m_{l}} -  I[g, w]_{t^m_{k}, t^m_{l}} \}
   \nn\\
 &\quad 
  +
   \{ I[g, w(m)]_{t^m_{l},t} -  I[g, w]_{t^m_{l},t} \}
\nn\\ 
&\quad 
+g^1_{s,t^m_{k}} \{ w(m)^1_{t^m_{k}, t^m_{l}}   -w^1_{t^m_{k}, t^m_{l}} \}
\nn\\ 
&\quad 
+g^1_{s,t^m_{k}} \{w(m)^1_{t^m_{l}, t}  -w^1_{t^m_{l}, t} \}
\nn\\ 
&\quad 
+g^1_{t^m_{k},  t^m_{l}}  \{ w(m)^1_{t^m_{l}, t} - w^1_{t^m_{l}, t}\}
=: 
B_1 +\cdots +B_6.
\label{eq.1011-2}
\end{align}
Note that $B_4 =0$ and $B_1$ and $B_3$ were already 
estimated in \eqref{ineq.1011-1}.
We can easily check that
\begin{align*}
\|B_6\|_{L^q} 
&\le
\|g\|_\alpha (t^m_{l}- t^m_{k})^\alpha 
  \{ \|w(m)^1_{t^m_{l}, t}\|_{L^q} + \|w^1_{t^m_{l}, t}\|_{L^q}\}
  \nn\\
   &\le
    C_2 \|g\|_\alpha (t^m_{l}-t^m_{k})^\alpha  (t- t^m_{l})^{\tfrac12}
    \nn\\
     &\le
    C_2 (2^{-m})^{\kappa}\|g\|_\alpha  (t- s)^{\alpha +\tfrac12 -\kappa}.
\end{align*}
Obviously, $B_5$ satisfies the same estimate, too.

It remains to estimate $B_2$ when $k<l$. First, we should note that 
\begin{align*}  
B_2=
\int_{t^m_{k}}^{t^m_{l}}  \left[
\sum_{j=k+1}^l
[g_{\cdot} - g_{t^m_{k}}]^m_j  \mathbf{1}_{[t^m_{j-1},t^m_{j} ]} (u)
- 
(g_{u} - g_{t^m_{k}})
\right] dw_u.
\end{align*}
The absolute value of the integrand is dominated by
$\|g\|_\alpha (2^{-m})^\alpha$. 
By Burkholder's inequality, we have
\[
\|B_2\|_{L^q}  
\le 
C_3 \|g\|_\alpha (2^{-m})^\alpha (t^m_{l}-t^m_{k})^{\tfrac12}
\le 
 C_3  (2^{-m})^\kappa \|g\|_\alpha(t- s)^{\alpha +\tfrac12 -\kappa}
\]
since $0< \kappa<\alpha$ and $t-s \ge 2^{-m}$.

Hence, taking $\ve_m :=C_4 (2^{-m})^\kappa$ for 
a suitable constant $C_4 >0$, we finish the proof of the lemma.
\end{proof}

\begin{lemma} \label{lem.1009-3}
Let $\alpha \in (0,1)$
and $\{g_m\}_{m\in\N} \subset  \cC^\alpha_0 (\R)$ be a 
sequence which converges to $g \in  \cC^\alpha_0 (\R)$ as $m\to\infty$
in the $\alpha$-H\"older norm. 
Then, 
for every $q \in [1,\infty)$ and $\kappa \in (0,\alpha)$, 
there exists a 
sequence of positive numbers $\{\tilde\ve_m\}_{m\in\N}$ converging to $0$
as $m\to\infty$ such that 
\[
\bigl\| I[g_m, w(m) ]_{s,t}- I[g, w]_{s,t}\bigr\|_{L^q} 
\le 
 C' \tilde\ve_m 
 (t-s)^{\alpha + \tfrac12 -\kappa},
\quad (s,t) \in \triangle_T, \,\, m\in\N.
\]
Here, we set $C':= \|g\|_\alpha \vee \sup_{m\ge 1}\|g_m\|_\alpha$ and
$\{\tilde\ve_m\}$ does not depend on $(s,t)$, $g$, $\{g_m\}$.
\end{lemma}

\begin{proof} 
We may assume $T=1$ as before.
The left hand side of the desired estimate 
is dominated by
\[
\bigl\| I[g_m-g, w(m) ]_{s,t}\bigr\|_{L^q} 
+
\bigl\| I[g, w(m) ]_{s,t}- I[g, w]_{s,t}\bigr\|_{L^q}.
\]
By Lemma \ref{lem.1009-1},
the first term is dominated by $C \| g_m-g\|_\alpha (t-s)^{\alpha +1/2}$.  
By Lemma \ref{lem.1009-2},
the second term is dominated by $\ve_m  \|g\|_\alpha (t-s)^{\alpha + \tfrac12 -\kappa}$.
By setting $\tilde\ve_m :=2C + \ve_m$, we complete the proof.
\end{proof}


\subsection{Construction of driving rough path}

Let $\tfrac14 <\alpha_0 \le \tfrac13$ and 
$(\Omega, \mathcal{F}, {\mathbb P})$ be a probability space. 
Let $w =(w_t)_{0\le t \le T}$ be a standard $e$-dimensional Brownian
motion 
and let $B=\{(B^1_{s,t}, B^2_{s,t}, B^3_{s,t})\}_{0\le s\le t\le T}$ be an 
$G\Omega_{\alpha} (\R^d)$-valued 
random variable (i.e., random RP)
defined on $(\Omega, \mathcal{F}, {\mathbb P})$
for every $\alpha \in (1/4,\alpha_0)$.
We assume that $w$ and $B$ are independent.
As for the integrability of $B$, 
Assumption {\bf (A)} is assumed, that is,
$\vertiii{B}_{\alpha}$ has moments of all orders.
Let $\{\cF_t\}_{0\le t\le T}$ be a filtration satisfying the usual condition
as well as the following two conditions:
 {\rm (i)} $w$ is an $\{\cF_t\}$-BM and 
{\rm (ii)} $t \mapsto (B^1_{0,t}, B^2_{0,t}, B^3_{0,t})$ is $\{\cF_t\}$-adapted.
We set $W =(W^1, W^2)$ as the 
Stratonovich-type Brownian RP, that is,
\[
W^1_{s,t} = w_t-w_s, \qquad  W^2_{s,t} = \int_s^t (w_u -w_s)\otimes 
\circ d w_u,
\qquad (s,t)\in \triangle_T,
\]
where $\circ dw$ stands for the Stratonovich integral.
It is well-known that $W \in G\Omega_{\gamma} (\R^e)$ a.s. 
and that $\vertiii{W}_{\gamma}$ has moments of all orders for every $\gamma \in (1/3, 1/2)$.

We set
\begin{align*}
I [B,W]_{s,t} &:= \int_s^t B^1_{s,u} \otimes d^{{\rm I}}w_u, 
\nn\\
I [W, B]_{s,t} &:= W^1_{s,t}\otimes B^1_{s,t} -  \int_s^t (d^{{\rm I}}w_u)\otimes  B^1_{s,u}.
\end{align*}
for every $(s,t)\in \triangle_T$.
Here, $d^{{\rm I}}w_u$ stands for the It\^o integration.
We can easily see that 
\begin{align*}
I [B,W]_{s,t}&=I [B,W]_{s,u}+I [B,W]_{u,t}+ B^1_{s,u} \otimes W^1_{u,t},
\\
I [W, B]_{s,t}&=I [W, B]_{s,u}+I [W, B]_{u,t}+ W^1_{s,u} \otimes B^1_{u,t}
\end{align*}
hold a.s. for each fixed $0\le s \le u \le t \le T$.

\begin{lemma}\label{lem.1014-1}
Let the situation be as above.
\\
{\rm (1)}~For every $(s,t)\in \triangle_T$, we have 
\[
I [B,W]_{s,t}= \lim_{|\cP|\searrow 0} \sum_{i=1}^N 
 B^1_{s, t_{i-1}} \otimes W^1_{t_{i-1}, t_i},
 \qquad
 I [W, B]_{s,t}= \lim_{|\cP|\searrow 0} \sum_{i=1}^N 
 W^1_{s, t_{i-1}} \otimes B^1_{t_{i-1}, t_i}.
\]
Here, the limits are in $L^2 ({\mathbb P})$ and 
$\cP =\{s=t_0 <t_1 <\cdots < t_N =t\}$ is a partition of $[s,t]$.
\\
{\rm (2)}~For every $\alpha \in (\tfrac14, \alpha_0)$
 and $q \in [1,\infty)$, 
 there exists a positive constant $C>0$ independent of $(s,t)$ such that
\[
\bigl\| I[B, W]_{s,t}\bigr\|_{L^q} \vee \bigl\| I[W, B]_{s,t}\bigr\|_{L^q} 
\le 
C (t-s)^{\alpha + \tfrac12},
\qquad (s,t) \in \triangle_T.
\]
\end{lemma}

\begin{proof} 
We can take any $k\in \llbracket 1, d\rrbracket$ and
 $l \in \llbracket 1, e\rrbracket$ and compute each $(k,l)$-component.
Hence, it is enough to prove the lemma for $d=e=1$.

We prove the first assertion.
The left one is almost obvious from the definition 
of It\^o integral. For the right one, note that 
\[
\sum_{i=1}^N  B^1_{s, t_{i-1}} W^1_{t_{i-1}, t_i}
+
\sum_{i=1}^N 
 W^1_{s, t_{i-1}} B^1_{t_{i-1}, t_i}
=
W^1_{s, t}B^1_{s, t} 
-
\sum_{i=1}^N 
 W^1_{t_{i-1}, t_i}  B^1_{t_{i-1}, t_i}.
\]
By the independence, it holds for $i<j$ that
\[
\E [ W^1_{t_{i-1}, t_i}  B^1_{t_{i-1}, t_i}
W^1_{t_{j-1}, t_j}  B^1_{t_{j-1}, t_j}]
=
\E[ B^1_{t_{i-1}, t_i} B^1_{t_{j-1}, t_j} ] 
\E[  W^1_{t_{i-1}, t_i} ]\E[  W^1_{t_{j-1}, t_j} ] =0.
\]
This implies that
\begin{align}  
\E \left[ \left|\sum_{i=1}^N 
 W^1_{t_{i-1}, t_i}  B^1_{t_{i-1}, t_i} \right|^2 \right]
 &=
 \sum_{i=1}^N \E [(W^1_{t_{i-1}, t_i}  B^1_{t_{i-1}, t_i})^2 ]
 \nn\\
 &\le 
  \sum_{i=1}^N \E [(W^1_{t_{i-1}, t_i})^4]^{1/2}  \E[(B^1_{t_{i-1}, t_i})^4 ]^{1/2}
   \nn\\
 &\le 
  c_1 \sum_{i=1}^N (t_i-t_{i-1})^{1+2\alpha} \le c_2 |\cP|^{2\alpha}.
  \nn
\end{align}
Here, we used Assumption  {\bf (A)} and 
$c_1$ and $c_2$ are certain positive constants.
The right hand side tends to zero as $|\cP| \searrow 0$.
Thus, we have shown (1).

From Burkholder's inequality and  Assumption  {\bf (A)}, 
the second assertion immediately follows.
\end{proof}

In what follows, we assume
 $\alpha \in (\tfrac14, \alpha_0)$ and $\gamma \in  (\tfrac13, \tfrac12)$
with $2\alpha +\gamma >1$.
By Proposition \ref{prop.1006-1} and Lemma \ref{lem.1014-1} (2), 
a continuous modification of $I[B, W]$ and hence that of $I[W, B]$ exist,
which will be denoted by the same symbols.
Moreover, $\|I[B, W]\|_{\alpha +\gamma}$ and $\|I[W, B]\|_{\alpha +\gamma}$
have moments of all orders.
Thus, we have seen that 
$(\Xi^1, \Xi^2, B^3) \in \hat{\Omega}_{\alpha, \gamma} (\cV)$, a.s.

\begin{definition} \label{def.random_uRP}
We write
\[
\Xi^1_{s,t} := \begin{pmatrix}
B^1_{s,t}  \\
W^1_{s,t} \\
\end{pmatrix},
\qquad
\Xi^2_{s,t}
:=\begin{pmatrix}
B^2_{s,t} & I[B, W]_{s,t} \\
I[W, B]_{s,t} & W^2_{s,t} \\
\end{pmatrix}
\]
and set $(\Xi^1, \Xi^2, \Xi^3):={\bf Ext}(\Xi^1, \Xi^2, B^3)$. 
(We also write $\Xi =(\Xi^1, \Xi^2, \Xi^3)$ for simplicity.)
As we will see below, $(\Xi^1, \Xi^2, B^3)\in 
G\hat{\Omega}_{\alpha, \gamma} (\cV)$ a.s.
Therefore, due to Proposition \ref{prop.gyur2},  $\Xi=(\Xi^1, \Xi^2, \Xi^3) \in G\Omega_{\alpha} (\cV)$ a.s.
\end{definition}

This random RP $\Xi$ is the driving RP of our slow-fast system of RDEs.
As we have seen in Proposition \ref{prop.gyur1}, 
$\|\Xi^{3, [ijk]}\|_{\delta}$, where $\delta:= (6-i-j-k)\alpha+ (i+j+k -3)\gamma$, have moments of all orders.

We write the components of $\Xi^3$ as follows:
$\Xi^{3, [111]} =B^3$ and $\Xi^{3, [222]}=W^3$.
For $(i,j,k) \neq (1,1,1), (2,2,2)$, we write 
$\Xi^{3, [121]} = I[B,W,B]$, $\Xi^{3, [221]} = I[W,W,B]$, etc.
Except when $(i,j,k) =(1,1,1)$, the H\"older regularity of 
$\Xi^{3, [ijk]}$ is larger than $1$.
Hence, although these seven components are involved 
in the Riemann-type sum for a RP integral along $\Xi$,
they actually make no contribution to the RP integral.

It remains to show that the random anisotropic RP $(\Xi^1, \Xi^2, B^3)$ 
is geometric.

\begin{lemma} \label{lem.250115}
Let the notation be as above. Then, $(\Xi^1, \Xi^2, B^3)\in 
G\hat{\Omega}_{\alpha, \gamma} (\cV)$, a.s.
\end{lemma}

\begin{proof} 
Clearly, $(\Xi^1, \Xi^2, B^3)$ is a functional of $B=(B^1, B^2, B^3)$ and $w$.
Since $B$ and $w$ are independent, 
we may write ${\mathbb P} = {\mathbb P}^B \times {\mathbb P}^W$
and ${\mathbb E} = {\mathbb E}^B \times {\mathbb E}^W$.
By Fubini's theorem, it is enough to show that, for every fixed 
realization of $B$, $(\Xi^1, \Xi^2, B^3)\in 
G\hat{\Omega}_{\alpha, \gamma} (\cV)$, ${\mathbb P}^W$-a.s.
Hence, we will assume below
that $B$ is an arbitrary (non-random) element 
of $G\Omega_{\alpha} (\cV_1)$.
For simplicity, we assume $T=1$ again.

Let $\{b (m)\}_{m\in\N} \subset \cC_0^1 (\cV_1)$ be a sequence
such that $\lim_{m\to\infty} S_3 (b(m)) =B$ in $G\Omega_{\alpha} (\cV_1)$. 
For the $m$th dyadic approximation $w(m)$, we write 
$(W(m)^1, W(m)^2)$ for $S_3 (w(m))$. 
It is well-known that for every $\gamma\in (1/3, 1/2)$ and 
$q\in [1,\infty)$, there are constants $c_1, c_2 >0$ and $r \in (0,1)$ such that
\begin{align*}  
\E^W [ \| W^i \|_{i\gamma}^{q/i} ]  \vee
\E^W [ \| W(m)^i \|_{i\gamma}^{q/i} ] \le c_1,
\qquad
\E^W [ \| W^i - W(m)^i\|_{i\gamma}^{q/i} ] \le c_2 r^m
\end{align*}
holds for all $m\in\N$ and $i=1,2$.
Here, $c_1, c_2, r$ are independent of $m$.

By Lemma \ref{lem.1009-1},  \ref{lem.1009-3} and 
Proposition \ref{prop.1007-2}, we have 
\[
\lim_{m\to\infty} \E^W \left[
\|I [b(m), w(m)]- I[b,w]\|_{\alpha + \gamma -2\kappa}
\right]=0
\]
for every sufficiently small $\kappa >0$.
From this we can easily see that 
\[
\lim_{m\to\infty} \E^W \left[
\|I [w(m), b(m)]- I[w, b]\|_{\alpha + \gamma -2\kappa}
\right]=0.
\]
Hence, a subsequence of $\{\hat{S} (b(m), w(m)) \}_{m\in\N}$
converges to $(\Xi^1, \Xi^2, B^3)$ in 
$\hat{\Omega}_{\alpha-\kappa, \gamma-\kappa} (\cV)$
and hence
$(\Xi^1, \Xi^2, B^3)\in G\hat{\Omega}_{\alpha-\kappa, \gamma-\kappa} (\cV)$, 
${\mathbb P}^W$-a.s.
Noting that $\alpha-\kappa$ and $\gamma-\kappa$ can be arbitrarily 
close to $\alpha_0$ and $1/2$, respectively, 
we complete the proof.
\end{proof}

In our construction of the random RP $\Xi$, It\^o integration is used. 
So, it is not a priori obvious whether $\Xi$ can be obtained 
via the piecewise linear approximation.
At least, in the case of fractional Brownian motion with 
Hurst parameter $H\in (1/4, 1/3]$, we can prove it.

\begin{remark} \label{rem.250115}
Let $b^H =(b^H_t)_{t \in [0,T]}$
and $w =(w_t)_{t \in [0,T]}$ 
be a $d$-dimensional fractional Brownian motion
with Hurst parameter $H\in (1/4, 1/3]$
and an $e$-dimensional Brownian motion, respectively, 
which are assumed to be independent.
Their $m$th
dyadic piecewise linear approximation ($m\in\N$) are denoted by 
$b^H (m)$ and $w(m)$, respectively.

According to \cite[Theorem 15.42 and Proposition 15.5]{fvbook},
the following limits of three sequences of random RPs exist
a.s and in $L^p~(1 \le p <\infty)$: 
\begin{enumerate} 
\item[(i)]~$\lim_{m\to\infty}S_3 ((b^H (m), w(m)))$ in 
$G\Omega_{\alpha} (\cV_1\oplus \cV_2)$ with $\alpha \in (1/4, H)$.
\item[(ii)]~$\lim_{m\to\infty} S_3 (b^H (m))$ in 
$G\Omega_{\alpha} (\cV_1)$ with $\alpha \in (1/4, H)$.
\item[(iii)]~$\lim_{m\to\infty} S_2 (w(m))$ in 
$G\Omega_{\gamma} (\cV_2)$ with $\gamma \in (1/3, 1/2)$.
\end{enumerate}
(Concerning the lift of Gaussian processes of this kind, 
a recent work \cite{gk} should also be referred to.)

We call $B^H = \lim_{m\to\infty} b^H (m)$ the fractional Brownian RP,
i.e. a canonical lift of $b^H$.
In this case we may take $\alpha_0 =H$.
We can show that the mixed RP 
$\Xi$ as in Definition \ref{def.random_uRP}
(with $B$ being replaced by $B^H$) coincides with 
$\lim_{m\to\infty} S_3 ((b^H (m), w(m)))$ as expected.
We will check this in the next paragraph.

Obviously from {\rm (ii)} and  {\rm (iii)} above, 
5 (out of 7) components of 
$\hat{S} ((b^H (m), w(m)))\in G\hat{\Omega}_{\alpha, \gamma} (\cV)$ converge a.s. 
In essentially the same way as in the proof of Lemma \ref{lem.250115},
\[
\lim_{m\to\infty} \E^W \left[
\|I [b^{H, j} (m), w^k (m)]- I[b^{H,j},w^k]\|_{\alpha + \gamma}
\right]=0
\]
and therefore 
\[
\lim_{m\to\infty} \E^W \left[
\|I [w^k(m), b^{H,j} (m)]- I[w^k, b^{H,j}]\|_{\alpha + \gamma}
\right]=0
\]
hold for almost all fixed $b^H$, 
where the superscripts $j~(1\le j \le d)$ and $k~(1\le k \le e)$ stands for the coordinates of $\cV_1=\R^d$ and $\cV_2 =\R^e$, respectively.
If we take a subsequence which may depend on $b^H$,
we have 
$\lim_{m\to\infty}\hat{S} ((b^H (m), w(m)))=(\Xi^1, \Xi^2, (B^{H})^{3})$ for almost all $w$.
Since ${\bf Ext}\circ \hat{S} =S_3$ and ${\bf Ext}$ is continuous, we 
have
\[
(\Xi^1, \Xi^2, \Xi^{3})={\bf Ext} ((\Xi^1, \Xi^2, (B^{H})^{3}))=
\lim_{m\to\infty} S_3 ((b^H (m), w(m)))
\quad 
\mbox{for a.a. $w$.}
\]
Note that the right hand side above converges even if 
we do not take subsequence.
Hence, we have $(\Xi^1, \Xi^2, \Xi^{3})=\lim_{m\to\infty} S_3 ((b^H (m), w(m)))$ for almost all $(b^H, w)$ by Fubini's theorem.
\end{remark}

\section{Slow-fast system of rough differential equations}
\label{sec.SF}

In this section we define the slow-fast system \eqref{def.SFeq}
of RDEs precisely and study its deterministic and probabilistic aspects. 
In this and the next sections, 
we always assume that $\sigma$ and $h$ are of $C_{{\rm b}}^4$ and 
$f$ and $g$ are locally Lipschitz continuous. 
The regularity parameters satisfy
that  $ \tfrac14 <\beta <\alpha <\alpha_0 \le \tfrac13$
and $\tfrac13 <\gamma <\tfrac12$ with $2\alpha +\gamma >1$. 
These assumptions guarantee that
 our slow-fast system of RDEs has a unique time-local solution
 up to the explosion time.
As before, we write
$\cV_1 =\R^d$, $\cV_2 =\R^e$ and $\cV=\cV_1 \oplus \cV_2$.

\subsection{Deterministic aspects}\label{subsec.deter}
First, we discuss deterministic aspects of our 
slow-fast system of RDEs.
As in Definition \ref{def.random_uRP}, let $(\Xi^1, \Xi^2, B^3)\in 
G\hat{\Omega}_{\alpha, \gamma} (\cV)$ and write
$\Xi={\bf Ext}(\Xi^1, \Xi^2, B^3)$. 
We will often write $\Xi:=(\Xi^1, \Xi^2, \Xi^3)\in G\Omega_{\alpha} (\cV)$.
This is the driving RP of our slow-fast system.
For a CP 
 $(Z, Z^\dagger, Z^{\dagger\dagger})\in \cQ^{\beta}_\Xi (\R^{m+n})$
with respect to $\Xi \in G\Omega_\alpha (\cV)$,
we often write $Z =(X,Y)$.

Respecting the direct sum decomposition 
$\R^{m+n}= \R^{m}\oplus \R^{n}$, 
a generic element of $\R^{m+n}$ is denoted by $z =(x,y)$.
The (partial) gradient operators with respect to $z$, $x$ and $y$ are 
denoted by $\nabla_z$, $\nabla_x$ and $\nabla_y$, respectively.
Hence, $\nabla_z =(\nabla_x, \nabla_y)$ at least formally
when it acts on nice functions on $\R^{m+n}$.
We will write the canonical projections as 
$\rho_1 \la z\ra =x$ and $\rho_2 \la z\ra =y$.
We set 
\[
F_{\varepsilon}(x, y)=\left(\begin{array}{c}f(x, y) 
\\ 
\varepsilon^{-1} g(x, y)\end{array}\right), 
\qquad
 \Sigma_{\varepsilon} (x, y)=\left(\begin{array}{cc}\sigma(x) & O 
 \\ 
 O & \varepsilon^{-1/2} h(x, y)\end{array}\right).
\]
Then, $F_{\varepsilon} \colon \R^{m+n}\to \R^{m+n}$ and 
$\Sigma_{\varepsilon} \colon  \R^{m+n}\to L(\cV_1 \oplus \cV_2, \R^{m+n})$.
In this subsection, $\varepsilon \in (0, 1]$ is arbitrary but fixed.

The precise definition of the slow-fast system \eqref{def.SFeq}
of RDEs (in the deterministic sense) is given as follows:  
\begin{align}\label{rde.SF}
Z_{t}^{\varepsilon} 
&=z_0 +\int_{0}^{t} F_{\varepsilon} (Z^{\varepsilon}_{s}) ds
+\int_{0}^{t} \Sigma_{\varepsilon} (Z^{\varepsilon}_{s}) d \Xi_{s},
\\
(Z^{\varepsilon})^\dagger_t &= \Sigma_\ve (Z^{\varepsilon}_t),
\qquad
(Z^{\varepsilon})^{\dagger\dagger}_t 
= (\nabla_z \Sigma_\ve \cdot \Sigma_\ve ) (Z^{\varepsilon}_t),
\qquad
t \in [0,T].
\nn
\end{align}
(Note that $(Z^{\varepsilon})^\dagger$ and 
$(Z^{\varepsilon})^{\dagger\dagger}$ take value in 
$L(\cV, \R^{m+n})$ and $L(\cV\otimes \cV, \R^{m+n})$, respectively.)
We consider this RDE in the $\beta$-H\"older topology.

Let 
$(Z^{\varepsilon}, \Sigma_{\varepsilon} (Z^{\varepsilon}),
(\nabla_z \Sigma_\ve \cdot \Sigma_\ve ) (Z^{\varepsilon}) )$
be a (necessarily unique) solution of RDE \eqref{rde.SF}
on a certain time interval $[0,\tau]$, where $\tau \in (0,T]$.
Then, the summand of
the modifies Riemann sum that approximates the RP integral 
on the right hand side of \eqref{rde.SF} is given by
\begin{align}
K_{s, t}
&:=
 \Sigma_{\varepsilon} (Z^{\varepsilon}_{s})  \Xi_{s, t}^{1}
+\{  \Sigma_{\varepsilon} (Z^{\varepsilon}) \}^{\dagger}_s \Xi_{s, t}^{2}
+ \{  \Sigma_{\varepsilon} (Z^{\varepsilon}) \}^{\dagger\dagger}_s
\Xi_{s, t}^{3},
\quad
(s,t)\in \triangle_{[0, \tau]}.
\label{eq.1029-1}
\end{align}
Here, 
\begin{align}
\{  \Sigma_{\varepsilon} (Z^{\varepsilon}) \}^{\dagger}_s
&= (\nabla_z \Sigma_\ve \cdot \Sigma_\ve ) (Z^{\varepsilon}_s)
=
\nabla_z \Sigma_{\varepsilon} (Z^{\varepsilon}_s)
\la  \Sigma_{\varepsilon} (Z^{\varepsilon}_s) \bullet, \star \ra
\in L(\cV^{\otimes 2}, \R^{m+n}),
\nn\\
\{  \Sigma_{\varepsilon} (Z^{\varepsilon}) \}^{\dagger\dagger}_s
&=
 \nabla_z \Sigma_{\varepsilon} (Z^{\varepsilon}_s)
  \bigl\la
   (\nabla_z \Sigma_\ve \cdot \Sigma_\ve ) (Z^{\varepsilon}_s) 
   \la  \bullet, \star\ra, \,*
   \bigr\ra
\nn\\
&\qquad 
+
\nabla_z^2  \Sigma_{\varepsilon} (Z^{\varepsilon}_s)
\la \Sigma_\ve (Z^{\varepsilon}_s)\bullet, \Sigma_\ve (Z^{\varepsilon}_s)
\star, \,*\ra
\in L( \cV^{\otimes 3}, \R^{m+n}).
\nn
\end{align}
Here, we used the third item of Example \ref{ex.0715}, again.

\begin{remark}
For the rest of this section, we use the following notation.
Let $\tau \in (0,T]$ and let $\cX$ be a Euclidean space.
If a continuous map $\eta \colon \triangle_{[0, \tau]} \to \cX$ 
belongs to $\cC^\delta_{(2)} ([0, \tau], \cX)$ for $\delta>0$,
we simply write $O ((t-s)^\delta)$ for $\eta_{s,t}$. 
It should be noted:
\begin{itemize}
\item
When we use this ``big $O$" symbol, we do not assume that 
$t-s$ is small.
\item
When $\eta$ depends on a driving RP $\Xi$ or $(\Xi^1, \Xi^2, B^3)$,
$\|\eta\|_{\delta, [0, \tau]}$ may depend on the RP (and the parameter $\ve$).
In other words, this symbol only means RP-wise estimates. 
\end{itemize}
\end{remark}

\begin{remark}
An element of a direct sum space  is denoted by
both a ``column vector" and a``row vector." 
These are not precisely distinguished.
\end{remark}

\begin{lemma}\label{lem.rde.Xve}
Let $0<\tau \le T$ and let the situation be as above. Suppose that 
\[
(Z^{\varepsilon}, \Sigma_{\varepsilon} (Z^{\varepsilon}),
(\nabla_z \Sigma_\ve \cdot \Sigma_\ve ) (Z^{\varepsilon}) )
\in \cQ^{\beta}_\Xi ([0, \tau], \R^{m+n}),
\]
is a unique solution of RDE \eqref{rde.SF} on $[0, \tau]$
and write $Z^{\varepsilon} =(X^{\varepsilon}, Y^{\varepsilon})$.

Then, $(X^{\varepsilon}, \sigma (X^{\varepsilon}), 
 (\nabla_x \sigma_\ve \cdot \sigma_\ve ) (X^{\varepsilon}))$ belongs to
$\cQ^{\beta}_B ([0, \tau], \R^{m})$
 and is a unique local solution of the following RDE driven by 
$B =(B^1, B^2, B^3)$ on $[0, \tau]$:
\begin{align}\label{rde.Fast}
X_{t}^{\varepsilon} 
&=x_0 +\int_{0}^{t} f (X^{\varepsilon}_{s}, Y^{\varepsilon}_{s}) ds
+\int_{0}^{t} \sigma (X^{\varepsilon}_{s}) dB_{s},
\\
(X^{\varepsilon})^\dagger_t &= \sigma(X^{\varepsilon}_t), 
\qquad 
(X^{\varepsilon})^{\dagger\dagger}_t 
= (\nabla_x \sigma \cdot \sigma ) (X^{\varepsilon}_t),
 \qquad t \in [0, \tau].
 \nn
\end{align}
Recall that an RDE of this type was introduced in \eqref{rde.0413}
and  discussed in Subsection \ref{sec.rde}.
\end{lemma}

\begin{proof}
By Proposition \ref{prop.gousei}, we have
\begin{equation}
X^{\ve}_t - X^{\ve}_s = O ((t-s)^1)+\rho_1 \la K_{s,t} \ra, 
\qquad (s,t)\in \triangle_{[0, \tau]}
\label{eq.1031-1}
\end{equation}
since $4\beta >1$.

It is clear that 
\begin{equation}
\rho_1 \la \Sigma_{\varepsilon} (Z^{\varepsilon}_{s})  \Xi_{s, t}^{1} \ra = \sigma (X^{\varepsilon}_{s})  B_{s, t}^{1}.
\label{eq.1101-1}
\end{equation}
Note that the $\R^m$-component of $ \Sigma_\ve (z)$ equals 
$\sigma (\rho_1 \la z\ra) \circ \pi_1 = \sigma (x) \circ \pi_1$.
In particular, $\nabla_y \rho_1\la \Sigma_\ve (z) \ra$ vanishes.
By standard calculation for block matrices, we see that
\begin{align}  
 \rho_1 \bigl\la 
 \nabla_z \Sigma_{\varepsilon} (Z^{\varepsilon}_s)
\la  \Sigma_{\varepsilon} (Z^{\varepsilon}_s) \bullet, \star \ra
 \bigr\ra 
 &=
  \nabla_z [\rho_1 \la \Sigma_{\varepsilon} (Z^{\varepsilon}_s)\ra ]
\la  \Sigma_{\varepsilon} (Z^{\varepsilon}_s) \bullet, \star \ra
\nn\\
&=
 ( \nabla_x  \sigma) (X^{\varepsilon}_s)
 \la \rho_1\la \Sigma_{\varepsilon} (Z^{\varepsilon}_s) \bullet \ra, \pi_1 \star \ra
 \nn\\
&=
 ( \nabla_x  \sigma) (X^{\varepsilon}_s)
 \la  \sigma (X^{\varepsilon}_s) \pi_1 \bullet, \, \pi_1 \star \ra
 \nn
\end{align}
and therefore
\begin{equation}
\rho_1 \la \{  \Sigma_{\varepsilon} (Z^{\varepsilon}) \}^{\dagger}_s \Xi_{s, t}^{2}\ra 
=
 ( \nabla_x  \sigma) (X^{\varepsilon}_s)
 \la  \sigma (X^{\varepsilon}_s) \bullet', \, \star'  \ra\vert_{( \bullet', \star')=B^2_{s,t}}
 = 
  ( \nabla_x  \sigma\cdot  \sigma) (X^{\varepsilon}_s) \la B^2_{s,t}\ra .
\label{eq.1101-2}
\end{equation}
Set $(X^\ve)^{\dagger} =  \sigma (X^{\varepsilon})$ and 
$(X^\ve)^{\dagger\dagger} =  ( \nabla_x  \sigma\cdot  \sigma) (X^{\varepsilon})$.
Then, we can easily see from \eqref{eq.1031-1}-- \eqref{eq.1101-2} that
$X^\ve_t - X^\ve_s = (X^\ve)^{\dagger}_s B^1_{s,t}
+(X^\ve)^{\dagger\dagger}_s B^2_{s,t} + O((t-s)^{3\beta})$.
We can also easily check that
$(X^\ve)^{\dagger}_t -(X^\ve)^{\dagger}_s= 
(X^\ve)^{\dagger\dagger}_s \la  B^1_{s,t}, * \ra+ O((t-s)^{2\beta})$.
Hence, 
$(X^\ve, (X^\ve)^{\dagger}, (X^\ve)^{\dagger\dagger})$ is a 
CP with respect to  $B=(B^1, B^2, B^3)$.
It should be noted that $ (X^\ve)^{\dagger\dagger} = \{\sigma (X^\ve)\}^\dagger$.

By cumbersome but similar calculations for block matrices as above,
we also have 
\begin{align}
\rho_1 \la \{  \Sigma_{\varepsilon} (Z^{\varepsilon}) \}^{\dagger\dagger}_s \Xi_{s, t}^{3}\ra 
&=
  \nabla_x \sigma (X^{\varepsilon}_s)
  \bigl\la
   (\nabla_x \sigma \cdot \sigma ) (X^{\varepsilon}_s) 
   \la  \bullet' , \star' \ra, \,*'
   \bigr\ra\vert_{( \bullet', \star',*')=B^3_{s,t}}
\nn\\
&\qquad 
+
\nabla_x^2  \sigma (X^{\varepsilon}_s)
\la \sigma (X^{\varepsilon}_s)\bullet', 
\sigma (X^{\varepsilon}_s)\star', \,*'\ra\vert_{( \bullet', \star',*')=B^3_{s,t}}
\nn\\
&=
 \{\sigma (X^\ve)\}^{\dagger\dagger}_s \la B^3_{s,t}\ra. 
\label{eq.1101-3}
\end{align}
Hence, we have 
\[
\rho_1 \la K_{s, t}\ra
 =
 \sigma (X^\ve_s) B^1_{s,t}+
  \{\sigma (X^\ve)\}^{\dagger}_s \la B^2_{s,t}\ra
 +
  \{\sigma (X^\ve)\}^{\dagger\dagger}_s \la B^3_{s,t}\ra
\]
and therefore
\[
\rho_1 \Bigl\la  \int_{s}^{t} \Sigma_{\varepsilon} (Z^{\varepsilon}_{u}) d \Xi_{u}
\Bigr\ra
 =
 \int_{s}^{t} \sigma (X^{\varepsilon}_{u}) dB_{u}.
\]
This completes the proof of the lemma.
\end{proof}

As for the slow component of $K_{s,t}$, 
we can easily see the following lemma, 
in which  $\nabla_x h (x,y)\cdot \sigma (x)$ is a shorthand for 
 the linear map from $\cV_1 \otimes \cV_2$ to $\R^n$ defined by
$\xi \otimes \eta \mapsto 
\nabla_x h (x,y) \la \sigma (x)\xi , \eta \ra$.
Note that $2\alpha +\gamma >1$ and that
the third term on the right hand side of \eqref{eq.1204-1}
is (formally) the same as the It\^o-Stratonovich 
correction term.
(In the probabilistic part of this paper, 
$W^2$ and $\overline{W}^2$ will be the second level of 
the Stratonovich-type and the It\^o-type  Brownian RP, respectively.)
\begin{lemma}\label{lem.1204hiru}
Let the assumptions be the same as in Lemma \ref{lem.rde.Xve} above. Then, we have
\begin{align}  
\rho_2 \la K_{s,t} \ra
&=
\varepsilon^{-1/2} h(X^\varepsilon_s, Y^\varepsilon_s) W^1_{s,t} 
+
\varepsilon^{-1/2} \nabla_x h(X^\varepsilon_s, Y^\varepsilon_s)\cdot
\sigma( X^\varepsilon_s ) 
\langle I[B,W]_{s,t}  \rangle
\nn\\
&\qquad +
\varepsilon^{-1} (\nabla_y h\cdot h)(X^\varepsilon_s, Y^\varepsilon_s)
\langle W^2_{s,t}  \rangle
  +O ((t-s)^{2\alpha +\gamma})
\nn\\
&=
\varepsilon^{-1/2} h(X^\varepsilon_s, Y^\varepsilon_s) W^1_{s,t} 
+
\varepsilon^{-1/2} 
\nabla_x h(X^\varepsilon_s, Y^\varepsilon_s)\cdot
\sigma( X^\varepsilon_s ) 
\langle I[B,W]_{s,t}  \rangle
\nn\\
&\qquad +
\varepsilon^{-1} 
(\nabla_y h \cdot h)
(X^\varepsilon_s, Y^\varepsilon_s)
\langle \overline{W}^2_{s,t}  \rangle
 \nn\\
 &\qquad 
 +\tfrac12 \varepsilon^{-1}  {\rm Trace}\bigl[
(\nabla_y h\cdot h)(X^\varepsilon_s, Y^\varepsilon_s)
\langle \bullet, \star   \rangle \bigr] (t-s)
  +O ((t-s)^{2\alpha +\gamma})
\label{eq.1204-1}
\end{align}
for all $(s,t)\in \triangle_{[0, \tau]}$.
Here, we set $ \overline{W}^2_{s,t}:= 
W^2_{s,t} -\tfrac{t-s}{2} \sum_{k=1}^e \mathbf{e}_k \otimes \mathbf{e}_k$ for the canonical orthonormal basis 
$\{\mathbf{e}_k\}_{k=1}^e$ of $\cV_2 =\R^e$.
\end{lemma}

\begin{proof}
We only check the first equality 
since the second one is obvious.
The components of $\Xi^3$  involved in $\rho_2 \la K_{s,t} \ra$ 
are $I[B, B, W]_{s,t}$, $I[B, W, W]_{s,t}$, $ I[W, B, W]_{s,t}$ and $W^3_{s,t}$, all of which are $O ((t-s)^{2\alpha +\gamma})$.
The rest is trivial.
\end{proof}

\subsection{Probabilistic aspects}\label{subsec.probab}

In what follows we work under ${\bf (A)}$ and
 assume that $1/4< \beta< \alpha <\alpha_0 (\le 1/3)$.
For the rest of this section,
$\Xi$ is as in Definition \ref{def.random_uRP}.  
The precise meaning of the random RDE in our main theorem is
 RDE \eqref{rde.SF} driven by this $\Xi$.

We extend the time interval of the filtration 
$\{\mathcal{F}_t\}$
by setting $\mathcal{F}_{t}=\mathcal{F}_{t\wedge T}$ for $t \ge 0$.
Denote by $\hat{\R}^{m+n} := \R^{m+n}\cup \{\infty\}$
the one-point compactification of $\R^{m+n}$.
If a global solution $(Z^{\ve}_t)_{t \in [0,T]}$ exists, 
then we set $Z^{\ve}_t =Z^{\ve}_{t\wedge T}$ for $t \ge 0$.
Otherwise, denote by
$(Z^{\ve}_t)_{t \in [0,u^\ve)}$, $0< u^\ve \le T$, 
be a maximal local solution
and set $Z^{\ve}_t =\infty$ for $t \in [u^\ve, \infty)$.
Either way, $(Z^{\ve}_t)$ is constant in $t$ on $[T,\infty)$, a.s.

Define
$\tau^\ve_N =\inf\{ t \ge 0 \mid  |Z^{\ve}_t|\ge N\}$
for each $N \in \N$
and $\tau^\ve_\infty = \lim_{N\to \infty}\tau^\ve_N$.
(As usual $\inf \emptyset :=\infty$.)
These are $\{\mathcal{F}_t\}$-stopping times. 
Then, the following are equivalent:

\begin{itemize} 
\item
A global solution $(Z^{\ve}_t)_{t \in [0,T]}$ of RDE \eqref{rde.SF} exists.
\item
$(Z^{\ve}_t)_{t \in [0,T)}$ defined as above is bounded in $\R^{m+n}$.
\item
$\tau^\ve_N =\infty$ for some $N$.
\item
$\tau^\ve_\infty > T$.
\end{itemize}
It should be noted that while a solution of RDE 
\eqref{rde.SF} moves in a bounded set, its trajectory is uniformly 
continuous in $t$ (because its H\"older norm is bounded). 
Hence, if $(Z^{\ve}_t)_{t \in [0,s)}$, $0<s\le T$, is bounded, 
then $(Z^{\ve}_t)_{t \in [0,s]}$ solves  RDE \eqref{rde.SF}.

On the other hand, if no global solution exists, 
then we have $u^\ve = \tau^\ve_\infty \in (0,T]$ and 
$\limsup_{t \nearrow \tau^\ve_\infty} |Z_t| =\infty$.
Moreover, $\lim_{t \nearrow \tau^\ve_\infty} |Z_t| =\infty$
because of the uniform continuity mentioned above.
Therefore, $(Z^{\ve}_t)_{t \ge 0}$ is a continuous process 
that takes values in $\hat{\R}^{m+n}$.


\begin{proposition} \label{pr.0520}
Let the notation be as above
and assume ${\bf (A)}$. Then, for every $\ve \in (0,1]$,
$Y^{\ve}$ satisfies
the following It\^o SDE up to the explosion time of 
$Z^{\ve}=(X^{\ve}, Y^{\ve})$:
\[
Y^{\ve}_t = y_0 + 
\frac{1}{\varepsilon}\int_0^{t \wedge T}
 \tilde{g}(X^\varepsilon_s, Y^\varepsilon_s) ds
 +
\frac{1}{\sqrt{\varepsilon}}\int_0^{t \wedge T}
 h(X^\varepsilon_s, Y^\varepsilon_s) d^{{\rm I}}w_s,
 \qquad
 0 \le t < \tau^{\ve}_{\infty}.
\]
Recall that $\tilde{g}(x,y)$ was defined by \eqref{def.tilde_g}, that is,
\[
\tilde{g}(x,y) = g (x,y) +\tfrac12   {\rm Trace}\bigl[
(\nabla_y h\cdot h)(x,y)
\langle \bullet, \star   \rangle \bigr].
\]
\end{proposition}

\begin{proof} 
The proof is very similar to the corresponding one 
in the case $\alpha_0 \in (1/3, 1/2]$ in \cite[Proposition 4.7]{ina_thk}.
Hence, our argument here is a little bit sketchy.

Let $\mathcal{P}=\{ 0=t_0 <t_1 <\cdots < t_K =t\}$ be a partition of 
$[0,t]$ for $0<t \le T$.
Recall that the summand of the Riemann sum for the RP integral in
RDE \eqref{rde.SF} was calculated in the previous subsection.

First, we prove the lemma 
when $h, \sigma$ are of $C^4_{{\rm b}}$
and $f, g$ are bounded and globally Lipschitz continuous.
In this case the solution never explodes, i.e. 
$\tau_{\infty}^{\varepsilon} =\infty$, a.s.
It is easy to see that 
\begin{equation}\nn
\lim_{|\mathcal{P} | \searrow 0}
\sum_{i=1}^K
h(X^\varepsilon_{t_{i-1}}, Y^\varepsilon_{t_{i-1}}) W^1_{t_{i-1},t_i}
=
\int_0^t  h(X^\varepsilon_s, Y^\varepsilon_s)  d^{{\rm I}}w_s  
\qquad
\mbox{in $L^2 ({\mathbb P})$}
\end{equation}
and 
\begin{eqnarray*}
\lim_{|\mathcal{P} | \searrow 0}
\sum_{i=1}^K
{\rm Trace}\bigl[
(\nabla_y h\cdot h)(X^\varepsilon_{t_{i-1}}, Y^\varepsilon_{t_{i-1}})
\langle \bullet, \star   \rangle \bigr] (t_{i}-t_{i-1})
\\
=
\int_0^t 
{\rm Trace}\bigl[
(\nabla_y h\cdot h)(X^\varepsilon_s, Y^\varepsilon_s)
\langle \bullet, \star   \rangle \bigr] ds, \qquad \mbox{a.s.}
\end{eqnarray*}
Note that $\tfrac12   {\rm Trace}\bigl[
(\nabla_y h\cdot h)(x,y)
\langle \bullet, \star   \rangle \bigr] = \tilde{g}(x,y)- g (x,y) $.

If $A_{s,t} = O ((t-s)^{2\alpha +\gamma})$, one can easily see that 
$
\lim_{|\mathcal{P} | \searrow 0}
\sum_{i=1}^KA_{t_{i-1},t_i}=0,\mbox{a.s.}
$
since $2\alpha +\gamma >1$.
By the same argument as in \cite{ina_thk}, we can prove
\begin{eqnarray}
\lim_{|\mathcal{P} | \searrow 0}
\sum_{i=1}^K 
(\nabla_y h \cdot h) (X^\varepsilon_{t_{i-1}}, Y^\varepsilon_{t_{i-1}})
\langle \overline{W}^2_{t_{i-1},t_i}  \rangle=0 \qquad
\mbox{in $L^2 ({\mathbb P})$,}
\nn\\
\lim_{|\mathcal{P} | \searrow 0}
\sum_{i=1}^K 
\nabla_y h(X^\varepsilon_{t_{i-1}}, Y^\varepsilon_{t_{i-1}}) \cdot 
\sigma (X^\varepsilon_{t_{i-1}}, Y^\varepsilon_{t_{i-1}})
\langle I[B,W]_{t_{i-1},t_i}  \rangle=0 \qquad
\mbox{in $L^2 ({\mathbb P})$.}
\nn
\end{eqnarray}
Note that $\overline{W}^2$ is the second level of 
Brownian RP of It\^o type.
Combining these all, we have shown that 
\begin{equation}\label{eq.0531-5}
Y^{\ve}_t = y_0 + 
\frac{1}{\varepsilon}\int_0^{t \wedge T}
 \tilde{g}(X^\varepsilon_s, Y^\varepsilon_s) ds
 +
\frac{1}{\sqrt{\varepsilon}}\int_0^{t \wedge T}
 h(X^\varepsilon_s, Y^\varepsilon_s) d^{{\rm I}}w_s,
 \qquad
t \ge 0
\end{equation}
holds a.s. in this case.

From here only the standing assumption is assumed 
on the coefficients $h, \sigma, f, g$.
Take any sufficiently large $N$.  
Let $\phi_N \colon \R^{m+n}\to [0,1]$ be a smooth function 
with compact support such that $\phi_N \equiv 1$ on 
the ball $\{ z\in \R^{m+n}\mid |z| \le N \}$
and set $\hat{h} := h \phi_N$. Also,
$\hat{\sigma}, \hat{f}, \hat{g}$ are defined in the same way.
We replace the coefficients of 
RDE \eqref{rde.SF} by these corresponding data with ``hat" and 
denote a unique solution by $\hat{Z}^\ve =(\hat{X}^\ve, \hat{Y}^\ve)$.
Then, \eqref{eq.0531-5} holds with 
$\hat{X}^\ve, \hat{Y}^\ve,  \hat{h}, \hat{g}$ 
in place of $X^\ve, Y^\ve, h, g$.
By the uniqueness of the RDE, it holds that 
$\hat{Z}^{\ve}_{t \wedge \tau^\ve_N} = Z^\ve_{t \wedge \tau^\ve_N}$
for all $0 \le t \le T$.
Therefore, we almost surely have 
\begin{align} 
Y^\ve_{t \wedge \tau^\ve_N  \wedge T}
&=
\hat{Y}^\ve_{t \wedge \tau^\ve_N  \wedge T}
\nn\\
&=
 y_0 + 
\frac{1}{\varepsilon}\int_0^{t \wedge \tau^\ve_N \wedge T}
 (\hat{g})^{\sim} (\hat{X}^\varepsilon_s, \hat{Y}^\varepsilon_s) ds
 +
\frac{1}{\sqrt{\varepsilon}}\int_0^{t \wedge \tau^\ve_N \wedge T}
 \hat{h} (\hat{X}^\varepsilon_s, \hat{Y}^\varepsilon_s) d^{{\rm I}}w_s,
\nn\\
&=
 y_0 + 
\frac{1}{\varepsilon}
\int_0^{t \wedge \tau^\ve_N \wedge T}
\tilde{g} (X^\varepsilon_s, Y^\varepsilon_s) ds
 +
\frac{1}{\sqrt{\varepsilon}}
\int_0^{t \wedge \tau^\ve_N\wedge T}
 h(X^\varepsilon_s, Y^\varepsilon_s) d^{{\rm I}}w_s,
  \quad
t \ge 0.
\label{eq.0601-1}
\end{align}
Since $\tau^\ve_N \nearrow \tau^\ve_\infty$ as $N \to \infty$ a.s. 
on  the set $\{  \tau^\ve_\infty \le T\}$, 
we  finish the proof by letting $N \to \infty$.
\end{proof}


In the same way as in  the author's previous work \cite{ina_thk},
we can show non-explosion of the solution under the assumptions 
on the coefficients.
Note that our assumptions on the coefficients are stronger 
than the corresponding ones in \cite{ina_thk}.
\begin{proposition} \label{prop.0601}
Assume ${\bf (A)}$ and ${\bf (H1)}$--${\bf (H5)}$.
Then, the probability that 
$Z^{\ve}=(X^{\ve}, Y^{\ve})$ explodes on $[0,T]$ is zero. 
Moreover, we have 
\begin{align}
\sup_{0<\ve \le 1}
\E [\|X^{\ve} \|_{\beta, [0,T]}^p] &<\infty, \qquad 1\le p<\infty,
\label{in.0602-1}
\\
\sup_{0<\ve \le 1}\sup_{0\le t \le T}\E [ |Y^{\ve}_t|^2]  &<\infty.
\label{in.0602-2}
\end{align}
\end{proposition}

\begin{proof}
Thanks to Proposition \ref{prop.0429},
Lemma \ref{lem.rde.Xve} and Proposition \ref{pr.0520},
the same proof as in \cite[Proposition 4.8]{ina_thk} still works.
\end{proof}


Now that Proposition \ref{prop.0601} has been obtained,
our arguments in what follows are very similar to 
the level 2 case in the corresponding part of \cite{ina_thk}.
Therefore, in order to avoid repetition, 
our proofs will be sketchy from here.

Now, we introduce a new parameter $\delta$ with $0<\ve <\delta \le 1$.
(In spirit, $0<\ve \ll \delta \ll 1$.
Later, we will set $\delta := \ve^{1/(6\beta)} \log \ve^{-1}$.)
We divide $[0,T]$ into subintervals of equal length $\delta$
(except perhaps the last subinterval).
For $s\ge 0$,  set $s(\delta) := \lfloor s/\delta\rfloor \delta$,
which is the nearest breaking point preceding or equal to $s$.

We set
\begin{equation}\label{def.0603}
\hat{Y}^{\ve}_t = y_0 + 
\frac{1}{\varepsilon}\int_0^{t }
 \tilde{g} (X^\varepsilon_{s(\delta)}, \hat{Y}^\varepsilon_s) ds
 +
\frac{1}{\sqrt{\varepsilon}}\int_0^{t}
 h(X^\varepsilon_{s(\delta)}, \hat{Y}^\varepsilon_s) d^{{\rm I}}w_s,
 \qquad
t \in [0,T].
\end{equation}
Note that $\hat{Y}^{\ve}$'s dependence on $\delta$ 
is suppressed in the notation.
This approximation process satisfies the following two estimates.
(The next two lemmas are basically the same as 
or easier than \cite[Lemma 5.1 and Lemma 5.2]{ina_thk}.)

\begin{lemma} \label{lem.0606-1}
Under the same assumptions as in Proposition \ref{prop.0601},
we have the following: 
For every $\delta$ and $\ve$ with $0<\ve <\delta \le 1$, 
the above process $\hat{Y}^{\ve}$ does not explode 
and satisfies
\begin{equation}\label{in.0707-1}
\sup_{0<\ve <\delta \le 1}\sup_{0\le t \le T}
\E [ |\hat{Y}^{\ve}_t|^{2}]  <\infty.
\end{equation}
\end{lemma}

\begin{proof}
The proof is essentially the same as that of Proposition \ref{prop.0601}.
(In fact, this one is easier because we already know 
$X^\varepsilon_{s(\delta)}$ exists and satisfies the estimate 
\eqref{in.0602-1}.)
\end{proof}

\begin{lemma} \label{lem.0603-1}
Assume ${\bf (A)}$ and ${\bf (H1)}$--${\bf (H6)}$.
Then, there exists a positive constant $C$ independent of $\delta$
such that
\[
\sup_{\ve\in (0,\delta)}\sup_{0\le t \le T} 
\E [|Y^\varepsilon_t -  \hat{Y}^\varepsilon_t|^2 ] 
\le C \delta^{2\beta}.
\]
\end{lemma}

\begin{proof} 
We omit the proof since 
it is essentially the same as in \cite[Lemma 5.2]{ina_thk}. 
\end{proof} 

%


It is easy to see that, if we define 
\begin{align}
M_t &=
           \int_{0}^{t} \{
             f (X_{s}^{\varepsilon}, Y_{s}^{\varepsilon} )
                -f (X_{s(\delta)}^{\varepsilon}, Y_{s}^{\varepsilon}) 
           \} d s
             +  \int_{0}^{t} \{f(X_{s(\delta)}^{\varepsilon}, 
                      Y_{s}^{\varepsilon})
              -f(X_{s(\delta)}^{\varepsilon}, \hat{Y}_{s}^{\varepsilon})
                \} d s
                       &\quad 
                       \label{def.230428}\\
                       &\qquad +
         \int_{0}^{t}  \{ 
            f(X_{s(\delta)}^{\varepsilon}, \hat{Y}_{s}^{\varepsilon})
              -\bar{f} (X_{s(\delta)}^{\varepsilon} )  \} d s
               +  
                 \int_{0}^{t} \{\bar{f} (X_{s(\delta)}^{\ve} )-\bar{f} (X_{s}^{\ve})) 
                 \} ds,
                 \nn
\end{align}
then 
\begin{align} 
X_{t}^{\varepsilon}-\bar{X}_{t} 
&=
\int_{0}^{t} f (X_{s}^{\varepsilon}, Y_{s}^{\varepsilon}) d s
   +\int_{0}^{t} \sigma (X_{s}^{\varepsilon} ) d B_{s} 
    - 
       \int_{0}^{t} \bar{f} (\bar{X}_{s}) d s
          -\int_{0}^{t} \sigma (\bar{X}_{s}) d B_{s}
          \label{eq.0603-1}\\ 
                   \nn\\
                   &= M_t + \Bigl( 
                       \int_{0}^{t} \{ \bar{f} (X_{s}^{\ve})
                           - \bar{f}(\bar{X}_{s})\}ds 
                                                + \int_{0}^{t} \{ \sigma (X_{s}^{\ve}) 
                               - \sigma (\bar{X}_{s}) \}d B_{s} \Bigr)
                            \nn
                                \end{align}
holds as an equality of CPs with respect to $B$.
We will later apply Proposition \ref{prop.0506}  to \eqref{eq.0603-1}
after estimating $\|M\|_{3\beta}$.

\begin{lemma} \label{lem.0606-2}
Assume the same condition as in Lemma \ref{lem.0603-1}
and let $1/2 <\gamma <1$.
Then, there exists a positive constant $C=C_\gamma$ 
independent of $\ve, \delta$ such that
\begin{align*}
\E \Bigl[
\Bigl\| \int_{0}^{\cdot} \{
             f (X_{s}^{\varepsilon}, Y_{s}^{\varepsilon} )
                -f (X_{s(\delta)}^{\varepsilon}, Y_{s}^{\varepsilon}) 
           \} d s \Bigr\|_{1}^2
           +
           \Bigl\| 
                \int_{0}^{\cdot} \{\bar{f} (X_{s(\delta)}^{\ve} )-\bar{f} (X_{s}^{\ve})) 
                 \} ds
                                \Bigr\|_{1}^2 
                                \Bigr]
                                \le C\delta^{2\beta},
\\  
\E \Bigl[\Bigl\| 
 \int_{0}^{\cdot} \{f(X_{s(\delta)}^{\varepsilon}, 
                      Y_{s}^{\varepsilon})
              -f(X_{s(\delta)}^{\varepsilon}, \hat{Y}_{s}^{\varepsilon})
                \} d s
 \Bigr\|_{\gamma}^2
 \Bigr]   \le        
                          C\delta^{2\beta (1-\gamma)}
\end{align*}
for all $0 < \varepsilon < \delta \le 1$.
Here, $\|\cdot \|_{\gamma}$ stands  for the $\gamma$-H\"older (i.e. Lipschitz) norm.
\end{lemma}

\begin{proof} 
The proof is easy and essentially the same as in 
\cite[Lemma 5.3]{ina_thk} (and its corrigendum).
Roughly speaking, we combine
 Proposition  \ref{prop.0601}, Lemma \ref{lem.0603-1}, 
the boundedness and global Lipschitz property of $f$.
We also use the quantitative estimate of local Lipschitz property of $\bar{f}$,
which is provided in \cite[Proposition C.5]{ina_thk} for instance.
\end{proof}

We can estimate the most difficult one among the four terms in the definition of $M$ as follows.
\begin{lemma} \label{lem.0606-3}
Assume the same condition as in Lemma \ref{lem.0603-1}
and let $0< \gamma <1$.
Then, there exists a positive constant $C$ independent of $\ve, \delta$
such that
\begin{align*}
\E \Bigl[
           \Bigl\| 
                \int_{0}^{\cdot} \{f (X_{s(\delta)}^{\ve}, \hat{Y}_{s}^{\varepsilon})-\bar{f} (X_{s(\delta)}^{\ve}) 
                 \} ds
                                \Bigr\|_{\gamma}^2 
 \Bigr]   \le        
                          C(\delta^{2(1-\gamma)} 
+ \delta^{-2\gamma} \ve )
\end{align*}
for all $0 < \varepsilon < \delta \le 1$.
\end{lemma}

\begin{proof} 
Essentially the same proof as in \cite[Lemma 5.4]{ina_thk}
works in this case, too.
However, it should be noted that the proof is not easy.
Both a careful approximation on each subinterval and
the Markov property of the frozen SDE must be used.
\end{proof} 

%
%

Now we prove our main theorem.
\begin{proof}[Proof of Theorem \ref{thm.main}]
Applying
Proposition \ref{prop.0506}  to \eqref{def.230428} and \eqref{eq.0603-1}, 
we obtain that
\begin{align} 
 \|X^{\ve} - \bar{X} \|_{\beta} &\le 
C \exp ( C \vertiii{B}_{\alpha}^\nu )  \,  \| M\|_{3\beta}
  \nn
 \end{align}
  for certain positive constant $C$ and $\nu$ which are 
  independent of $\ve, \delta$. (Below, $C$ 
     and $\nu$ may vary from line to line.)
By Lemmas \ref{lem.0606-2} and \ref{lem.0606-3}, 
there exists a positive constant $\kappa =\kappa_\beta >0$ 
 independent of $\ve, \delta$ such that
\[
\E [ \| M\|_{3\beta}^2 ] \le C (\delta^{\kappa} 
+ \delta^{-6\beta} \ve),
\qquad 0 <  \ve< \delta \le 1.
\]
Therefore, if we set $\delta := \ve^{1/(6\beta)} \log \ve^{-1}$ for example,
then $\| M\|_{3\beta}$ converges to $0$ in $L^2$-sense as $\ve\searrow 0$.
It immediately follows that $\|X^{\ve} - \bar{X} \|_{\beta}^p$ 
converges to $0$ in probability as $\ve\searrow 0$
for every $p \in [1,\infty)$.

On the other hand, 
we see from Proposition \ref{prop.0429} that
\[
\sup_{0<\ve \le 1}
\E
[ \|X^{\ve} - \bar{X}\|_{\beta}^p]  \le 
 2^{p-1}\sup_{0<\ve \le 1}
\E
[ \|X^{\ve}\|_{\beta}^p] + 2^{p-1}\E [\|\bar{X}\|_{\beta}^p]
\le C  (\E[\vertiii{B}_{\alpha}^{\nu p} ]+1) <\infty
\]
for every $p \in [1, \infty)$.
This implies that $\{ \|X^{\ve} - \bar{X}\|_{\beta}^p \}_{0<\ve \le 1}$
are uniformly integrable for each fixed $p$.
Hence, we have 
$\E [\|X^{\ve} - \bar{X}\|_{\beta}^p] \to 0$ as $\ve\searrow 0$.
This completes the proof of the main theorem.
\end{proof}

%

\medskip
\noindent
{\bf Acknowledgement:}~
The author is supported by JSPS KAKENHI Grant No. 23K20216.

%

\bigskip
\bigskip
\begin{flushleft}
  \begin{tabular}{ll}
    Yuzuru \textsc{Inahama}
    \\
    Faculty of Mathematics,
    \\
    Kyushu University,
    \\
    744 Motooka, Nishi-ku, Fukuoka, 819-0395, JAPAN.
    \\
    Email: {\tt inahama@math.kyushu-u.ac.jp}
  \end{tabular}
\end{flushleft}

\end{document}